\documentclass[final,onefignum,onetabnum]{siamart190516}

\usepackage{amsfonts}
\usepackage{graphicx}
\usepackage{epstopdf}
\usepackage{algorithm,algorithmic}
\usepackage{mathtools}
\usepackage{commath,nicefrac}
\usepackage{comment}
\usepackage{setspace}
\usepackage{eucal}

\usepackage{subfiles}

\usepackage[normalem]{ulem}

\usepackage{tikz}
\usepackage{bm}
\ifpdf%
  \DeclareGraphicsExtensions{.eps,.pdf,.png,.jpg}
\else
  \DeclareGraphicsExtensions{.eps}
\fi

\DeclareMathOperator{\rank}{rank}
\DeclareMathOperator{\nz}{nz}

\DeclareMathOperator{\diag}{diag}

\def\ddefloop#1{\ifx\ddefloop#1\else\ddef{#1}\expandafter\ddefloop\fi}
\def\ddef#1{\expandafter\def\csname bf#1\endcsname{\ensuremath{\mathbf{#1}}}}
\ddefloop ABCDEFGHIJKLMNOPQRSTUVWXYZabcdefghijklmnopqrstuvwxyz\ddefloop

\def\ddef#1{\expandafter\def\csname bb#1\endcsname{\ensuremath{\mathbb{#1}}}}
\ddefloop ABCDEFGHIJKLMNOPQRSTUVWXYZ\ddefloop

\def\ddef#1{\expandafter\def\csname cal#1\endcsname{\ensuremath{\CMcal{#1}}}}
\ddefloop ABCDEFGHIJKLMNOPQRSTUVWXYZ\ddefloop

\DeclareMathAlphabet{\mathpzc}{OT1}{pzc}{m}{it}

\def\ddef#1{\expandafter\def\csname pz#1\endcsname{\ensuremath{\mathpzc{#1}}}}
\ddefloop ABCDEFGHIJKLMNOPQRSTUVWXYZabcdefghijklmnopqrstuvwxyz\ddefloop

\def\ddef#1{\expandafter\def\csname rm#1\endcsname{\ensuremath{\mathrm{#1}}}}
\ddefloop ABCDEFGHIJKLMNOPQRSTUVWXYZabcdefghijklmnopqrstuvwxyz\ddefloop

\makeatletter
\renewcommand*\env@matrix[1][*\c@MaxMatrixCols c]{%
  \hskip -\arraycolsep
  \let\@ifnextchar\new@ifnextchar
  \array{#1}}
\makeatother

\usepackage{enumitem}
\usepackage{booktabs}

\newsiamremark{remark}{Remark}
\newsiamremark{hypothesis}{Hypothesis}
\crefname{hypothesis}{Hypothesis}{Hypotheses}
\newsiamthm{claim}{Claim}

\headers{Preconditioned Inexact Interior Point Method}
{S.~Karim, and E.~Solomonik}

\title{%
Efficient Preconditioners for Interior Point Methods via a new Schur Complement-Based Strategy
  \thanks{Submitted to the editors April 29, 2021.}
}

\author{Samah Karim\thanks{Department of Computer Science, University of Illinois at Urbana-Champaign, Urbana, IL 61801 
  (\email{swkarim2@illinois.edu}).}
\and Edgar Solomonik\thanks{Department of Computer Science, University of Illinois at Urbana-Champaign, Urbana, IL 61801 
  (\email{solomon2@illinois.edu}, \url{http://solomonik.cs.illinois.edu}).}
}

\usepackage{amsopn}

\ifpdf%
\hypersetup{%
  pdftitle={single-fact-ipm},
  pdfauthor={S.~Karim, and E.~Solomonik}
}
\fi

\usepackage{cleveref}

\begin{document}
\maketitle

%\onehalfspacing %remove this later
\begin{abstract}
%One of the main benefits of our formulation is the ability to re-use the factorization of the invariant equality-constrained part of KKT system, which is computed only once (per quadratic program). Further, two preconditioners are suggested specifically for the iterative solution of the new reduced system. In the high number of equalities regime, the low-degree-of-freedom preconditioner is shown to bounds the number of unique eigenvalues, and subsequently the number of conjugate gradient iteration count. The same goes for the high-degree-of-freedom preconditioner in the low number of equalities regime. 
We propose a novel preconditioned inexact primal-dual interior point method for constrained convex quadratic programming problems.  The algorithm we describe invokes the preconditioned conjugate gradient method on a new reduced Schur complement KKT system, in implicit form.
%that is different from the reduced system employed in prior approaches.
%prior of than strategy produces an original symmetric positive definite inequality-constrained reduced system.
%as well as two preconditioners designed specifically for the iterative solution of this system. 
%One of the main benefits of our formulation is the ability to
In contrast to standard approaches, the Schur complement formulation we consider enables reuse of the factorization of the KKT matrix with rows and columns corresponding to inequality constraints excluded, across all interior point iterations. %\samahnew{This fixed subsystem can be thought of as the KKT system of another quadratic program which has the same equality constraints, but no inequality constraints.}
%of the Hessian and equality constraints.
%the factorization of the invariant equality-constrained part of KKT system, which is computed only once (per quadratic program).
%Further, two preconditioners are suggested specifically for the iterative solution of the new reduced system.
Further, two new preconditioners are presented for the resulting reduced system, that alleviate the ill-conditioning associated with slack variables in primal-dual interior point methods.
%\samahnew{, especially near optimality.
Each of the preconditioners we propose also provably reduces the number of unique eigenvalues for the coefficient matrix, and thus the CG iteration count.
One preconditioner is efficient when the number of equality constraints is small, while the other is efficient when the number of remaining degrees of freedom is small.
%We suggest two preconditioners for this Schur complement system.
%In the high number of equalities regime, the low-degree-of-freedom preconditioner is shown to bounds the number of unique eigenvalues, and subsequently the number of conjugate gradient iteration count. The same goes for the high-degree-of-freedom preconditioner in the low number of equalities regime. 
%Similarities and differences to the conventional Schur complement augmented system formulation are discussed, in terms of popular preconditioning techniques, spectral properties and performance models. 
%Numerical experiments on synthetic test sets validate our theoretical conjugate gradient iteration bounds for both preconditioners.
Numerical experiments with synthetic problems and problems from the Maros-M\'esz\'aros QP collection show that our preconditioned inexact interior point solvers are effective at improving conditioning and reducing cost.
%approaches outperform other preconditioners.
%Our preconditioners are highly ef reduced system is typically much better conditioned than preconditioned systems.
Across all test problems for which the direct method is not fastest, our preconditioned methods achieve a reduction in cost by a geometric mean of $1.432$ relative to the best alternative preconditioned method for each problem. 
\end{abstract}

% REQUIRED

% REQUIRED

% Column width: \the\columnwidth
%
% Line width: \the\textwidth
%
% Text width: \the\linewidth

% Notes:
%  - labelname = fig:figurename
%  - labels should be attached to the caption
%  - use ~ to attach cite and ref
%  - use % to attach a figure after the citing paragraph

\section{Introduction} \label{sec:intro}
%Start with motivation. 
%Nocedal: QP are important in their own right
Constrained quadratic programming problems are important standalone problems, and they also arise as sub-problems when solving nonlinear programming problems with a sequential quadratic programming (SQP) approach~\cite{NoceWrig06}. %For the past few decades, interior point methods (IPM) have gained wide appreciation due to their remarkable success in solving linear and non-linear optimization problems; The readers are referred to \cite{wright1992interior} and   \cite{forsgren2002interior} for excellent comprehensive surveys.
%In this work, we propose a new strategy where we perform a single factorization in a setup phase of IPM, and then we solve a new reduced symmetric positive definite linear system at every IPM iteration using an iterative method with preconditioners specifically designed for this system.
In this paper, we 
 consider interior point methods for the solution of the following constrained convex quadratic programming (QP) problem,
\begin{equation}
\label{qp}
	\begin{aligned}
\underset{x \in \mathbb{R}^n}{\text{minimize}} \quad & \frac{1}{2}x^{T}Hx+x^{T}c\\
\textrm{subject to} \quad & Ax=b,\\
  &Cx\geq d,    \\
\end{aligned}
\end{equation}
where $c, x \in \mathbb{R}^n$, $b \in \mathbb{R}^{m_1}$, $d \in \mathbb{R}^{m_2}$, 
$H \in \bbR^{n\times n}$ is a symmetric matrix corresponding to the Hessian of the objective function, $A \in \bbR^{m_1\times n}$ encodes equality constraints, and $C \in \bbR^{m_2\times n}$ encodes inequality constraints.
Often these inequality constraints take the form of {\it simple bounds}, for instance, non-negativity constraints $(x \geq 0)$, in which case $C = I$, or upper and lower bounds $(l \leq x \leq u)$, where $C = [I~-I]^T$. We assume that $0 < m_1 \leq n$ and $A$ is of full rank, as $A$ can be pre-processed to remove linearly dependent rows if necessary~\cite{NoceWrig06}.
{We also assume throughout this paper that $H$ is positive-definite.
%, without great loss of generality; %. This assumption may appear to be rather restrictive at first glance, %since $H$ is singular in many applications~\cite{haber2001preconditioned} @ADD CITATIONS. 
While such strictly convex optimization problems are of independent interest~\cite{bergamaschi2004preconditioning,ali2018iterative}, in general, regularization techniques can be used to push the eigenvalues of the Hessian away from zero. Regularization is commonly used as it guarantees that an optimal solution exists, is bounded and is unique~\cite{gill1991solving,gill1993solving,saunders1996cholesky,gondzio2012matrix,altman1999regularized,pougkakiotis2019dynamic,pougkakiotis2021interior}. We discuss pathways for extension of our method to handle semi-definite $H$ in the conclusion (\cref{sec:conc}).
} %The reader is referred to~\cite{gondzio2012matrix,altman1999regularized,pougkakiotis2019dynamic,pougkakiotis2021interior} for the benefits of using regularization within IPMs.

For the past few decades, interior point methods (IPM) have gained wide appreciation due to their remarkable success in solving linear and non-linear optimization problems~\cite{wright1992interior,forsgren2002interior}.
These iterative methods are dominated in cost by the solution of the Karush-Kuhn-Tucker (KKT) linear system of equations at every iteration. 
The reduced KKT matrix arising at the $k$th interior point iteration may be written as
\begin{equation}
\label{augmented-intro}
K^{(k)}= \begin{bmatrix}
-H & A^T & C^T\\ 
A & 0 & 0\\
C & 0 & D^{(k)}
\end{bmatrix},
\end{equation}
where $D^{(k)}$ is a positive diagonal matrix dependent on the values of Lagrange multipliers and slack variables (interior point method parameters) at the $k$th iteration.
By %\sout{Schur complementation} 
block elimination of the inequality constraints, the linear system is typically transformed
into a more compact form  known as the {\it augmented system}.
At the $k$th iteration this linear system involves the following matrix,
\begin{equation}
\label{augmented-2-intro}
K_C^{(k)} =
	\begin{bmatrix}
-\left(H+C^T\bigl(D^{(k)}\bigr)^{-1}C\right) & A^T \\
A & 0 
\end{bmatrix}.
%\begin{pmatrix}
%\Delta x\\ 
%\Delta \lambda
%\end{pmatrix} =  -\begin{pmatrix}
%r_u\\ 
%r_e
%\end{pmatrix}.
\end{equation}
%\sout{where $\lambda \in \mathbb{R}^{m_1}$, and $\nu \in \mathbb{R}^{m_2}$ are Lagrange multipliers,}\sout{$D = V^{-1}S \in \mathbb{R}^{m_2 \times m_2}$ is a positive diagonal matrix, $V$ and $S \in \mathbb{R}^{m_2 \times m_2}$ are diagonal matrices containing the elements of $\nu$ and $s$ on their diagonals respectively.}
% ... are some IPM parameters.
%\sout{A further Schur complementation results in a reduced system referred to as the {\it normal equations},}

Sparse direct methods have been the most popular approaches for solving the symmetric indefinite augmented system~\eqref{augmented-2-intro}~\cite{zhang1998solving,gertz2003object,wachter2006implementation}. 
However, for certain QP problems these direct approaches may suffer from 
poor performance due to significant amount of fill~\cite{al2008preconditioning}, in addition to poor scaling~\cite{rees2007preconditioner} and memory limitations as the problem size increases.  
For these reasons, iterative solvers and in particular Krylov-subspace methods have emerged as a viable alternative for the solution of these linear systems.  
An IPM in which the KKT system is solved approximately, using an iterative method for instance, is sometimes referred to as an {\it inexact} or truncated IPM~\cite{bellavia1998inexact,freund1999convergence}.

A principal challenge for iterative solvers of these linear systems is the poor conditioning of the KKT matrix in late steps of the IPM method~\cite{forsgren2007iterative,benzi2005numerical}.
The systems become ill-conditioned due to highly variable magnitude of the elements in $D^{(k)}$ near the solution~\cite{NoceWrig06}.
Therefore, it is essential to use an effective preconditioner
to accelerate the convergence of iterative solvers.
A number of previous works have studied preconditioners for saddle point problem $K_C^{(k)}$, 
%for the saddle point matrix in \eqref{augmented-2-intro}, 
most notably constraint preconditioners~\cite{keller2000constraint,gould2001solution,bergamaschi2004preconditioning,durazzi2003indefinitely,dollar2007using,bellavia2015updating,lukvsan1998indefinitely,perugia2000block, rozloznik2002krylov,benzi2008some,benzi2005numerical} and augmented Lagrangian preconditioners~\cite{golub2003solving,greif2006preconditioners,rees2007preconditioner,cao2008augmentation,pearson2017fast,bergamaschi2021new,morini2016spectral,shen2012augmentation} that rely on the block structure of the augmented system.
Such sophisticated preconditioners typically require the factorization of a matrix approximating $K_C^{(k)}$ at each iteration.
%, as reviewed in section~\ref{sec:other-prec}. 

%Edgar: call normal eq: eq-constrained reduced system vs. ours ineq constr reduced system 
 
In this work, instead of the Schur complementation of the inequality constraints leading to $K_C^{(k)}$, we instead consider a Schur complement reduction step that results in a new inequality-constraint reduced system.
By maximally delaying the factorization of the block that changes at each interior point iteration ($D^{(k)}$), our method requires only a 
single factorization of the 
invariant matrix $F =  \bigl[ \begin{smallmatrix} -H & A^T\\
A &0 \end{smallmatrix}\bigr]$. %\samahnew{, which we denote as the {\it KKT subsystem}.}
This matrix %, which we refer to as the  {\it KKT subsystem}, 
can be thought of as the KKT system of another QP with the same equality constraints, but no inequality constraints.
 %,$F = \Big[\begin{matrix} -H & A^T\\A &0\end{matrix}\Big]$.
Our {\it single factorization inexact IPM algorithm} iteratively solves a system with the symmetric positive-definite (SPD) matrix,
%a new Schur complementation strategy for the solution of KKT systems arising in primal-dual interior point methods is proposed, resulting in a new reduced linear system given by,     
\begin{equation}
\label{new-reduced-intro}
\begin{aligned}
    K_F^{(k)} = D^{(k)}  - 
    \begin{bmatrix}
C & 0
\end{bmatrix} F^{-1} \begin{bmatrix}
C^T\\ 
0
\end{bmatrix},
%\Delta \nu = -r_a + \begin{bmatrix}
%C & 0
%\end{bmatrix} F^{-1}\begin{pmatrix}
%r_g\\ 
%r_e
%\end{pmatrix}.
\end{aligned}
\end{equation}
leveraging the factorization of $F$ to perform products with $K_F^{(k)}$.
Effective preconditioners are proposed for the reduced matrix $K_F^{(k)}$  that directly alleviate the ill-conditioning caused by the diagonal matrix $D^{(k)}$.
These preconditioners reduce the number of unique eigenvalues, which also bounds the iteration count of the conjugate gradient method~\cite{hestenes1952methods} in exact arithmetic.
Consequently, these preconditioners achieve improvements in asymptotic cost complexity.
We first consider the preconditioner $P_L = D^{(k)}$, which is especially effective when the number of degrees of freedom $n-m_1$ is low, as $ 
    \begin{bmatrix}
C & 0
\end{bmatrix} F^{-1} \begin{bmatrix}
C  &
0
\end{bmatrix}^T$ is of low rank in this case.
When the number of degrees of freedom is high, we instead propose the preconditioner
$P_H = D^{(k)} + CH^{-1}C^T$. {In both cases, rank analysis provides an upper bound on the number of non-unit eigenvalues of the preconditioned matrices
$P_L^{-1/2}K_F^{(k)}P_{L}^{-1/2}$
and $P_H^{-1/2}K_F^{(k)}P_{H}^{-1/2}$.}

\Cref{tab:comp_methods} compares these preconditioned methods to two state-of-the-art preconditioning techniques for $K_C^{(k)}$ (constraint preconditioning, given by $P_{CP}$, and augmented Lagrangian preconditioning, given by $P_{RG}$, both described in~\cref{sec:other-prec}).
The trade-off between the methods is clearly nontrivial as they factorize different matrices.
However, the new preconditioned  methods require fewer factorizations, while reducing the theoretical iteration count in comparison to alternative known preconditioned iterative solvers for this system.
\begin{table}
\begin{tabular}{c|c|c|c|c}
Preconditioned matrix & $n_d$ &  Factorized matrix & $n_f$ & $n_{eig}$\\ \hline
$K_C^{(k)}$ & $n + m_1$ & - & $0$ & $n+m_1$\\ \hline
$\bigl(K_C^{(k)}\bigr)^{-1}K_C^{(k)}$ & $n + m_1$ & $AG^{-1}A^T$& $N_I$ & $0$\\ \hline
$P_{CP}^{-1}K_C^{(k)}$  & $n + m_1$ & $AE^{-1}A^T$ & $N_I$ &  $n-m_1$  \\ \hline
$P_{RG}^{-1}K_C^{(k)}$  & $n + m_1$ & $G+\gamma A^TA$ & $N_I$ &  $m_1$  \\ \hline
$K_F^{(k)}$  & $m_2$ & $AH^{-1}A^T$ & $1$ &  $m_2$  \\ \hline
$P_L^{-1}K_F^{(k)}$  & $m_2$ & $AH^{-1}A^T$&$1$ &  $n-m_1$\\ \hline
$P_H^{-1}K_F^{(k)}$  & $m_2$ & $AH^{-1}A^T$&$1$ &  $m_1$\\ \hline
\end{tabular}
\caption{The first two rows represent (1) unpreconditioned iterative and (2) direct solve methods. $G=H+C^T(D^{(k)})^{-1}C$ is assumed to be positive-definite, $E$ is chosen as an approximation of $H$, $n$ is the number of variables, $m_1$ is the number of equality constraints, $N_I$ is the number of interior point iterations, $n_d$ is the dimension of the linear system, $n_f$ is the number of factorizations, $n_{eig}$ is the number of non-unit eigenvalues.}\label{tab:comp_methods}
\end{table}

In~\cref{sec:numerics}, we provide experimental evidence of the efficacy of the new single-factorization inexact IPM, with preconditioner $P_L$ or $P_H$. 
We verify our spectral analysis via experiments on synthetic test problems.
Based on analytic cost models derived in~\cref{sec:perf-analysis}, we compare the performance of the proposed methods to IPM with direct solve as well as unpreconditioned, constraint-preconditioned, and augmented-Lagrangian-preconditioned iterative methods.
We measure the condition numbers and predicted costs obtained by the preconditioned iterative strategies on test problems from the Maros-M\'esz\'aros QP collection.
%We summarize improvements in both metrics by considering the geometric mean across these datests.
The condition number achieved by our preconditioners is on average\footnote{Geometric mean across test problems of ratio of improvement in the geometric mean of the condition numbers of linear systems arising across all interior point iterations for that problem.} lower by a factor of 240 relative to the best among alternative preconditioners for each problem.
%base-10 logarithm of the geometric mean of the condition number across IPM iterations for each test problem achieved by the preconditioned methods is less than those of the alternative preconditioners by a mean (geometric mean across test problems) factor of \samahnew{2.37}.
%Our preconditioners achieve lowest cost among both direct and iterative variants for 44\% of the Maros-M\'esz\'aros test problems.
Our single factorization inexact IPM method achieves the lowest cost among iterative IPM variants for %44\% 
2/3 of the Maros-M\'esz\'aros test problems, when using either $P_L$ or $P_H$ as a preconditioner.
Relative to the best alternative preconditioner, they reduce cost by a mean factor of 1.432, where the former refers to the geometric mean across test problems for which one of the iterative methods was the fastest.% \sout{(geometric mean across test problems for which the direct method is not fastest) factor of 1.432.} 

\section{Background on interior point methods}\label{sec:qp}
As a preliminary step to solving~\eqref{qp}, the {\it Lagrangian} function is defined as 
\begin{equation}
	\label{lagrangian}
	L(x,\lambda, \nu) = \frac{1}{2}x^THx + x^Tc -\lambda^T (Ax-b) - \nu^T(Cx - d),
\end{equation}
where $\lambda \in \mathbb{R}^{m_1}$ is the vector of Lagrange multipliers corresponding to the equality constraints, and $\nu \in \mathbb{R}^{m_2}$ is the vector of Lagrange multipliers corresponding to the inequality constraints.
%, also known as the dual variables.
For convex QPs, the first order optimality conditions, also known as the Karush-Kuhn-Tucker (KKT) conditions, are both necessary and sufficient to guarantee the global optimality of the solution~\cite{wright1997primal},
\begin{subequations}
	\begin{align} 
\nabla_x L(x, \lambda, \nu)  &=  0,\\
Ax-b &= 0,\\
Cx - d &\geq 0,\\
\nu^T(Cx - d) &= 0,\\
\nu &\geq 0.
\end{align}
\end{subequations}
By introducing a vector of slack variables $s \geq 0$, we can rewrite these conditions as

\begin{subequations}
	\begin{align} 
Hx+c-A^T\lambda - C^T \nu &=  0,\\
Ax-b &= 0 \label{primal-eq},\\
Cx - d - s &= 0 \label{primal-ineq},\\
SVe &= 0,\label{complem}\\
s, \nu &\geq 0 \label{nonlinear_ineq},
\end{align}
\end{subequations}
where
 \begin{align*}
 	V = \diag(\nu_1,\nu_2,\ldots,\nu_{m_2}), &&S = \diag(s_1,s_2,\ldots,s_{m_2}), && e =[1,\ldots1]^T\in \mathbb{R}^{m_2}.
 \end{align*}

Primal-dual interior point methods solve a perturbed version of the KKT conditions where the inequalities~\eqref{nonlinear_ineq} are {\it strictly} satisfied and the average value of the pairwise products $s_i\nu_i$, for $i=1,\ldots,m_2$ in~\eqref{complem} is reduced at every iteration.  In particular, a primal-dual IPM solves the system of equations, 

\begin{subequations}
\begin{align}\label{nonlinear}
 \begin{bmatrix}Hx+c-A^T\lambda - C^T \nu\\
Ax-b\\
Cx - d - s \\
SVe - \sigma \mu e \end{bmatrix} &=0,\\ \label{strict-pos}
s, \nu &> 0,
\end{align}
\end{subequations}
where $\sigma \in [0,1]$ is a centering parameter, and the {\it duality measure} $\mu$ is defined as
\begin{equation}
	\label{duality}
	\mu = \frac{1}{m_2} \sum_{i=1}^{m_2} s_i \nu_i = \frac{s^T\nu}{m_2}.
\end{equation}
The perturbed KKT conditions~\eqref{nonlinear} are solved using Newton's method for nonlinear equations in order to generate search directions $(\Delta x^{(k)},\Delta \lambda^{(k)},\Delta \nu^{(k)},\Delta s^{(k)})$, where $k$ is the iteration number. The resulting system of linear equations  to be solved at every interior point iteration will be
\begin{equation}
\label{nonsym-indef-KKT}
\begin{bmatrix}
H & -A^T & -C^T & 0\\ 
A & 0 & 0 & 0\\ 
C & 0 & 0 & -I\\
0 &0 & S^{(k)}  & V^{(k)}
\end{bmatrix}\begin{pmatrix}
\Delta x^{(k)}\\ 
\Delta \lambda^{(k)}\\ 
\Delta \nu^{(k)} \\
\Delta s^{(k)}
\end{pmatrix} =  - \begin{pmatrix}Hx^{(k)}+c-A^T\lambda^{(k)} - C^T \nu^{(k)}\\ Ax^{(k)}-b\\ Cx^{(k)} -d -s^{(k)}\\S^{(k)}V^{(k)}e - \sigma \mu e\end{pmatrix}.
\end{equation}
This KKT system~\eqref{nonsym-indef-KKT} is not symmetric but a multiplication of the first equation by $-I$
%\samahnew{in anticipation of later symmetrization,} will
results in the following symmetric indefinite system,
\begin{equation}
\label{sym-indef-KKT}
	\begin{bmatrix}
-H & A^T & C^T & 0\\ 
A & 0 & 0 & 0\\ 
C & 0 & 0 & -I\\
0 &0 & S^{(k)}  & V^{(k)}
\end{bmatrix}\begin{pmatrix}
\Delta x^{(k)}\\ 
\Delta \lambda^{(k)}\\ 
\Delta \nu^{(k)} \\
\Delta s^{(k)}
\end{pmatrix} =  -\begin{pmatrix}
r_g^{(k)}\\ 
r_e^{(k)}\\ 
r_i^{(k)}\\
r_c^{(k)}
\end{pmatrix},
\end{equation}
where
\begin{equation}
\label{rhs-KKT4}
	\begin{aligned}
		r_g &= -Hx^{(k)}-c+A^T\lambda^{(k)}+C^T \nu^{(k)},\\
		r_e &= Ax^{(k)}-b,\\
		r_i &= Cx^{(k)}-d-s^{(k)},\\
		r_c &=S^{(k)}V^{(k)}e - \sigma \mu e.
	\end{aligned}
\end{equation}

\begin{comment}
\begin{equation}
\label{sym-indef-KKT}
	\begin{bmatrix}
-H & A^T & C^T & 0\\ 
A & 0 & 0 & 0\\ 
C & 0 & 0 & -I\\
0 &0 & -I  & -\bigl(S^{(k)}\bigr)^{-1}V^{(k)}
\end{bmatrix}\begin{pmatrix}
\Delta x^{(k)}\\ 
\Delta \lambda^{(k)}\\ 
\Delta \nu^{(k)} \\
\Delta s^{(k)}
\end{pmatrix} =  -\begin{pmatrix}
r_g^{(k)}\\ 
r_e^{(k)}\\ 
r_i^{(k)}\\
r_c^{(k)}
\end{pmatrix},
\end{equation}

where
\begin{equation}
%\label{rhs-KKT4}
	\begin{aligned}
		r_g &= -Hx^{(k)}-c+A^T\lambda^{(k)}+C^T \nu^{(k)},\\
		r_e &= Ax^{(k)}-b,\\
		r_i &= Cx^{(k)}-d-s^{(k)},\\
		r_c &=-\nu^{(k)}+\sigma \mu \bigl(S^{(k)}\bigr)^{-1} e .
	\end{aligned}
\end{equation}
\end{comment}

\subsection{Standard reduced systems}\label{sec:redsystem}
Alternative to the formulation in~\eqref{sym-indef-KKT}, it is very common to reduce the KKT system into a smaller symmetric indefinite system that is cheaper to factorize. This is achieved by elimination of
\begin{equation}
    \label{delz}
	\Delta s^{(k)} %= -\bigl(V^{(k)}\bigr)^{-1}\bigl(r_c^{(k)} + S^{(k)}\Delta \nu^{(k)}\bigr) 
	= -\bigl(V^{(k)}\bigr)^{-1}
	r_c^{(k)} - D^{(k)}\Delta \nu^{(k)}
\end{equation}
from the third equation in~\eqref{sym-indef-KKT}, where we have introduced the notation
\begin{equation}
\label{diag-d}
    D^{(k)} = \bigl(V^{(k)}\bigr)^{-1}S^{(k)}.
\end{equation}
The resulting block $3\times 3$ system will be referred to as the {\it first augmented system} throughout this paper,

\begin{equation}
\label{augmented}
\underbrace{	\begin{bmatrix}
-H & A^T & C^T\\ 
A & 0 & 0\\
C & 0 & D^{(k)}
\end{bmatrix}}_{K^{(k)}}
\begin{pmatrix}
\Delta x^{(k)}\\ 
\Delta \lambda^{(k)}\\
\Delta \nu^{(k)}
\end{pmatrix}%}_{\Delta_3^{(k)}} 
=  -%\underbrace{
\begin{pmatrix}
r_g^{(k)}\\ 
r_e^{(k)}\\
r_a^{(k)}
\end{pmatrix},
%}_{\bfr_1^{(k)}},
\end{equation}
where $r_a^{(k)} = r_i^{(k)} + \bigl(V^{(k)}\bigr)^{-1} r_c^{(k)}$, and $K^{(k)}$ $\in \bbR^{(n+m_1+m_2)\times(n+m_1+m_2)}$ denotes the coefficient matrix of this system. {Such a problem is typically addressed in the literature as a double (or twofold) saddle point system, since it has a saddle point structure and contains a submatrix that is itself a saddle point matrix~\cite{benzi2019uzawa, ramage2013preconditioned,gatica2000dual,ali2018iterative}.}
One can go a step further by eliminating $\Delta \nu^{(k)} = -\bigl(D^{(k)}\bigr)^{-1}\bigl(r_a^{(k)} + C \Delta x^{(k)}\bigr) \label{delnu}$ from the first equation in~\eqref{augmented} to get a smaller linear system, known as the {\it augmented system},  
%can say block $2 \times 2$
\begin{equation}
\label{augmented-2}
	\underbrace{\begin{bmatrix}
-\left(H+C^T\bigl(D^{(k)}\bigr)^{-1}C\right) & A^T \\A & 0 
\end{bmatrix}}_{K_C^{(k)}}%\underbrace{
\begin{pmatrix}
\Delta x^{(k)}\\ 
\Delta \lambda^{(k)}
\end{pmatrix}%}_{\Delta_2^{(k)}} 
=  -%\underbrace{
\begin{pmatrix}
r_u^{(k)}\\ 
r_e^{(k)}
\end{pmatrix},%}_{\bfr_2^{(k)}},
\end{equation}
%\edgar{drop $\bfr_2$}
\noindent where $r_u^{(k)} = r_g^{(k)} - C^T\bigl(D^{(k)}\bigr)^{-1} r_a^{(k)}$, and $K_C^{(k)}\in \bbR^{(n+m_1)\times(n+m_1)}$.  {This block $2\times2$ linear system is sometimes referred to as a classical saddle point problem~\cite{sogn2019schur}.}
A further Schur complement reduction step results in an equality-constraint reduced system $\in \bbR^{m_1 \times m_1}$,
known as the {\it normal equations},
%Even further reduction is possible by eliminating $\Delta x^{(k)}$ from (\ref{augmented-2}) in order to get the following system $\in \bbR^{m_1 \times m_1}$, known as the {\it normal equations}, 
{
\begin{equation}
\label{del_lambda}
	\underbrace{
	\left[A\left(H+C^T\bigl(D^{(k)}\bigr)^{-1}C\right)^{-1} A^T\right]}_{K_H^{(k)}} 
	\Delta \lambda^{(k)} 
	= -%\underbrace{
	\left[ A\left(H+C^T\bigl(D^{(k)}\bigr)^{-1}C\right)^{-1}r_u^{(k)} + r_e^{(k)}\right].
	%}_{\bfr_1^{(k)}}.
\end{equation}}
After solving for $\Delta \lambda^{(k)}$, the search direction in the primal variable can be recovered from
$\left(H+C^T\bigl(D^{(k)}\bigr)^{-1}C\right)\Delta x^{(k)} =  A^T \Delta \lambda^{(k)} + r_u^{(k)}.$
%It is worth noting that the coefficient matrix in the normal equations (\ref{del_lambda}), is the Schur complement of $K_C^{(k)}$, the coefficient matrix in the augmented system (\ref{augmented-2}). 
%%%%%%%%%%%%%%%%%
In practice, a plethora of factors can make the block $3 \times 3$ first augmented system~\eqref{augmented}, the block $2 \times 2$ augmented system~\eqref{augmented-2} or the normal equations~\eqref{del_lambda} a more favorable  formulation in terms of total computational cost and numerical accuracy.
The coefficient matrices of the $3\times3$, and the $2\times2$ systems, $K^{(k)}$, and $K_C^{(k)}$, are both symmetric indefinite but the former may be more sparse than the latter. On the other hand, the coefficient matrix for the normal equations~\eqref{del_lambda}, $K_H^{(k)}$, is symmetric positive-definite and tends to be even less sparse.
But all KKT matrices $K^{(k)}$, $K_C^{(k)}$, and $K_H^{(k)}$, preserve their initial sparsity structures throughout the interior point algorithm as the diagonal matrix $D^{(k)}$ is the only component that changes at every IPM iteration.  

%Whether one choose to solve unreduced primal-dual KKT system, the block $2 \times 2$ augmented system (\ref{augmented-2-intro}) or the normal equations (\ref{normal-eqs-intro}), it is widely known that the numerical solution of that linear system is the most computationally expensive step at every primal-dual IPM iteration \cite{benzi2005numerical}.  In practice a plethora of factors can make one formulation more favorable than the rest in terms of total computational cost as well as numerical accuracy including but not restricted to, the sparsity of the Hessian $H$ and the constraint matrices $A$ and $C$ as well as the type of linear system solver whether it be direct or iterative.

Many IPM software packages use sparse direct methods to compute the Newton search directions. For instance LIPSOL~\cite{zhang1998solving} uses a Cholesky-based factorization to solve the normal equations~\eqref{del_lambda}, while OOQP~\cite{gertz2003object} and IPOPT~\cite{wachter2006implementation} use a symmetric indefinite $LDL^T$ factorization to solve the augmented system~\eqref{augmented-2}, where $L$ is unit lower triangular and $D$ is diagonal with $2\times2$ or $1\times1$ pivots~\cite{arioli1989augmented,duff2017direct}.
 However, the $L$ factor given by any of the direct solvers often suffers from a significant amount of fill-in vis-a-vis the original matrix for both~\eqref{augmented-2}, and~\eqref{del_lambda}~\cite{al2008preconditioning}.   
 Moreover, the coefficient matrices for the KKT systems have to be formed and factorized at every IPM iteration due to the change in the diagonal elements of the matrix $D^{(k)}$. 
 Both the computation and the factorization procedures may get prohibitively expensive as the QP size increases~\cite{wang2000adaptive}. 
 
 On the other hand, the iterative solution of the linear system at every IPM iteration 
 allows for significant savings in memory and computational effort as the coefficient matrix need not be explicitly formed. When solving any of the KKT systems in~\eqref{augmented},~\eqref{augmented-2}, or~\eqref{del_lambda}, all the information that an iterative linear solver needs is the action of the coefficient matrix on a vector. 
 In addition, 
 %\sout{iterative solvers for instance Krylov-subspace methods, unlike direct solvers,allow the use of approximate search directions \cite{bellavia1998inexact} and that is especially useful in early iterations of the inexact IPM when the iterate is far from the optimal solution.}
 iterative solvers provide tunable accuracy, and high accuracy may only be needed in late iterations of IPM~\cite{bellavia1998inexact}.
But, the matrix $D^{(k)}$ presents a particular challenge for any type of iterative method in this context. It exhibits increasing ill-conditioning, as some of its diagonal elements become very large and others  very tiny as the interior point method approaches optimality~\cite{NoceWrig06}. 
And since this behavior often translates into the ill-conditioning of KKT matrices $K^{(k)}$, $K_C^{(k)}$, and $K_H^{(k)}$, at late IPM iterations, Krylov-subspace methods are known to be accurate only when used with an effective preconditioner~\cite{forsgren2007iterative}. %\sout{A preconditioner is necessary in improving the convergence properties of the iterative method by reducing the condition number and/or achieving good clustering of the eigenvalues of the preconditioned matrix~\cite{demmel1997applied,golub2013matrix}.}
 {Prior approaches to circumvent the ill-conditioning of the $3\times3$ KKT matrix include proposing equivalent systems that do not include the matrix $D^{(k)}$~\cite{forsgren2002inertia}. By premultiplying the third equation in~\eqref{augmented} by $V^{1/2}$, and scaling the variable $\Delta \nu^{(k)}$ by $V^{-1/2}\Delta \nu^{(k)}$, one can get a system with coefficient matrix,
 \begin{equation}
\label{scaled-augmented}
\hat{K}^{(k)}	= \begin{bmatrix}
-H & A^T & C^T(V^{(k)})^{\frac{1}{2}}\\ 
A & 0 & 0\\
(V^{(k)})^{\frac{1}{2}}C & 0 & S^{(k)}
\end{bmatrix}.%}_{K^{(k)}}
%\begin{pmatrix}\Delta x^{(k)}\\ \Delta \lambda^{(k)}\\(V^{(k)})^{-\frac{1}{2}}\Delta \nu^{(k)}\end{pmatrix} =  -\begin{pmatrix}r_g^{(k)}\\ r_e^{(k)}\\(V^{(k)})^{\frac{1}{2}}r_a^{(k)}\end{pmatrix}.
\end{equation}
Spectral analysis in~\cite{doi:10.1137/120890600,morini2016spectral} shows that this matrix may be well-conditioned throughout the IPM iterations, and is at the very least analytically nonsingular.
However, to recover the search direction in the dual variables $\Delta \nu^{(k)}$ one needs to apply another ill-conditioned transformation, $V^{-1/2}\Delta \nu^{(k)}$~\cite{forsgren2002inertia}. 
%Moreover, a study done in~\cite{morini2017comparison} shows that 
Further, when using effective preconditioning, iterative solvers exhibit similar convergence behavior with any of the $3\times3$ formulations~\eqref{augmented}, and~\eqref{scaled-augmented}, as well as the $2\times2$ formulation, despite the apparent superior spectral properties of~\eqref{scaled-augmented}, when unpreconditioned~\cite{morini2017comparison}.
Our approach may be applied to~\eqref{scaled-augmented} via a change of variables, but this also results in equivalent or similar methods after preconditioning is applied.
%Thus based on their theoretical and experimental analysis they conclude that, it is more advantageous to use the $2\times2$ augmented system~\eqref{augmented-2}, as opposed to the $3\times3$ system~\eqref{scaled-augmented} in terms of computational time. 

The augmented system is also oftentimes the system of choice, as opposed to the normal equations, as the coefficient matrix $K_H^{(k)}$ corresponding to~\eqref{delz} is likely to become more ill-conditioned than $K_C^{(k)}$ corresponding to~\eqref{augmented-2}, as IPM approaches an optimal solution~\cite{arioli1989augmented,NoceWrig06}. 
} {Although many preconditioners have been formulated for the normal equations~\cite{bocanegra2013improving,dravzic2015sparsity,CASACIO2017129}, they
 have generally not resulted in uniformly improved results~\cite{lustig1992implementing}.} 
 {
 %An additional reason why 
 In our experimental results, we consider preconditioned solvers for the $2\times2$ block system as a baseline, as such solvers have been shown to perform best in prior work~\cite{morini2017comparison}.
 %formulation will be the baseline that we will compare to later on, is the abundance of preconditioners for this system.
 This has been explained in part by an analysis presented in~\cite{oliveira2005new}, where  multiple different preconditioners for the augmented system were shown to yield an identical preconditioner for the normal equations.} In Section~\ref{sec:other-prec}, we give an overview of some popular preconditioning techniques for the $2\times2$ saddle point system in~\eqref{augmented-2}.

%Nocedal: General purpose sparse Cholesky can be applied to normal eqs but modifications are needed because it may be ill-conditioned or singular especially during the final stages of primal-dual algorithm, when the elements of D take on both huge and tiny values. Cholesky may encounter diagonal elements that are very small, zero, or because of roundoff error slightly negative.
 
%In \cite{morini2017comparison}, Simoncini, Morini and Tani conducted a thorough theoretical and experimental analysis in order to assess the desirability of the $3 \times 3$ block vis-à-vis the $2 \times 2$ block {\it augmented system} formulations when using iterative linear solvers with an appropriate preconditioning strategy.  

%Also see: Matrix-free interior point method Jacek Gondzio for more discussion on choice of system for iterative solver

\section{Preconditioners for the augmented system}\label{sec:other-prec}
For saddle point problems such as~\eqref{augmented-2}, the quality of a preconditioner varies depending on the problem at hand. Hence there is no `best' preconditioner for linear systems that have the same $2 \times 2$ block structure as the augmented system~\cite{benzi2005numerical}. But several classes of preconditioners have been successful in the iterative solution of such linear systems in optimization, most notably constraint preconditioners and augmented Lagrangian preconditioners which we will review in the next subsections. {Interested readers may refer to~\cite{axelsson2003preconditioning,pearson2020preconditioners} for comprehensive surveys on preconditioning techniques in constrained optimization.  }  
\subsection{Constraint preconditioners}

Constraint preconditioning is a popular technique that has been used in the solution of saddle point linear systems stemming from a mixed finite element formulation of elliptic partial differential equations in~\cite{bank1989class,ewing1990preconditioning,mihajlovic2004efficient,perugia2000block,rozloznik2002krylov,tong1998iterative,axelsson2003preconditioning}, and references therein. This technique has also been utilised in constrained optimization problems~\cite{benzi2005numerical} for the solution of block $2 \times 2$ KKT saddle point systems of the form,
\begin{equation}
    \label{saddle-point-system}
    \underbrace{\begin{bmatrix}
G^{(k)} & A^T \\ 
A & 0 
\end{bmatrix}}_{\CMcal{A}}
\begin{pmatrix}
x\\y
\end{pmatrix}=\begin{pmatrix}
b\\c
\end{pmatrix},
\end{equation}
%\samah{it's a system so I think the above equation is better than this one:}{\color{red} I agree, please integrate.}
%\begin{equation*}\label{saddle-point-system}\CMcal{A} = \begin{bmatrix}G^{(k)} & A^T \\ A & 0 \end{bmatrix},\end{equation*}
where $G^{(k)} \in \bbR^{n \times n}$ is symmetric positive semi-definite,  and $A \in \bbR^{m_1 \times n}$ is of 
full rank.
%\sout{a full rank ($=m_1 \leq n$) matrix as above.}
For a quadratic programming problem within a primal-dual interior point  formulation, $G^{(k)} = -\bigl(H+C^T\bigl(D^{(k)}\bigr)^{-1}C\bigr)$ ; see ~\eqref{augmented-2}. 
%and $C$ is symmetric but is equal to zero in many cases, such as in (\ref{augmented-2}). 
Constraint preconditioners (CPs) are indefinite saddle point matrices that have the same $2 \times 2$ block structure as the original saddle point system~\eqref{saddle-point-system}. The constraints in the $(1,2)$ and $(2,1)$ blocks of $\CMcal{A}$ %\sout{the saddle point matrix} 
are kept intact, and the matrix $G^{(k)}$ in the $(1,1)$ block is approximated by another $n \times n$ matrix $E^{(k)}$ that is much easier to invert, 
\begin{equation}
    \label{const-prec}
    P_{CP} =  \begin{bmatrix}
E^{(k)} & A^T \\ 
A & 0 
\end{bmatrix}.
\end{equation}
The nomenclature {\it constraint preconditioning} comes from the fact that the preconditioner %\sout{matrix}
$P_\text{CP}$ is the augmented KKT matrix for a QP with the same equality constraints. 
Multiple choices are possible for $E^{(k)}$ that insure the invertibility of the preconditioner, the simplest of which being $E^{(k)} = I$~\cite{gould2001solution},~\cite{benzi2005numerical}. When $G^{(k)}$ has positive diagonal entries, a common alternative is to take $E^{(k)} = \diag(G^{(k)})$~\cite{bergamaschi2004preconditioning,bellavia2015updating,benzi2005numerical}. 

CPs have been used to solve equality constrained nonlinear optimization problems (NLPs) in~\cite{lukvsan1998indefinitely}, as well as inequality constrained NLPs in~\cite{bergamaschi2004preconditioning,dollar2007using}. They have also been used to solve KKT systems coming from equality constrained QPs in~\cite{keller2000constraint,gould2001solution}. While in~\cite{bergamaschi2004preconditioning,durazzi2003indefinitely,dollar2007using}, CPs are used specifically for augmented KKT systems given by an interior point method formulation for an equality-constrained QP with simple bounds.  

Several authors have studied the eigenvalues and the corresponding eigenvectors of the preconditioned matrix $P_{CP}^{-1}\CMcal{A}$; see, e.g.,~\cite{lukvsan1998indefinitely,perugia2000block, rozloznik2002krylov,durazzi2003indefinitely}. Under the assumptions that $E^{(k)} = (E^{(k)})^T \neq G^{(k)}$, and $P_\text{CP}$ in~\eqref{const-prec} is nonsingular, the preconditioned matrix $P_{CP}^{-1}\CMcal{A}$ has eigenvalue $1$ with algebraic multiplicity $2m_1$, and $n-m_1$ eigenvalues defined by the generalized eigenproblem $Z^TG^{(k)}Z x_z = \lambda Z^T E^{(k)} Z x_z$, where $Z \in \bbR^{n \times (n-m_1)}$ is a matrix whose columns form a basis for the nullspace of $A$~\cite{keller2000constraint}. Generally, the better $E^{(k)}$ approximates $G^{(k)}$, the tighter is the cluster of eigenvalues around $1$.  

The application of a constraint preconditioner to a linear system involves computing its factorization either through a direct $LDL^T$ decomposition of $P_\text{CP}$ as has been done in~\cite{keller2000constraint}. But it is not always the case that the symmetric indefinite preconditioner in~\eqref{const-prec} is significantly cheaper to factorize than the original symmetric indefinite saddle point system in~\eqref{saddle-point-system}~\cite{benzi2005numerical,benzi2008some}.
Alternatively the factorization of $P_\text{CP}$ can be computed through the following block decomposition~\cite{benzi2005numerical},
\begin{equation*}
    \label{cp-factor}
    \footnotesize{
 \begin{bmatrix}
E^{(k)} & A^T \\ 
A & 0 
\end{bmatrix}^{-1} =  \Bigg[\begin{matrix}
I & -\bigl(E^{(k)}\bigr)^{-1}A^T \\ 
0 & I 
\end{matrix}\Bigg] \Bigg[\begin{matrix}
\bigl(E^{(k)}\bigr)^{-1} & 0 \\ 
0 & -\left(A\bigl(E^{(k)}\bigr)^{-1}A^T\right)^{-1} 
\end{matrix}\Bigg] \Bigg[\begin{matrix}
I & 0 \\ 
-A\bigl(E^{(k)}\bigr)^{-1} & I 
\end{matrix}\Bigg]. }
\end{equation*}
%\sout{Using} 
Per this identity, the cost of factorizing the constraint preconditioner $P_\text{CP}$ depends on that of factorizing $E^{(k)}$ and $A\bigl(E^{(k)}\bigr)^{-1}A^T$. If $E^{(k)}$ is sufficiently simple, for instance diagonal, then its %\sout{the}
factorization %\sout{of $E^{(k)}$} 
is trivial to compute.
%Furthermore, the factorization of $A^{(k)}(E^{(k)})^{-1}A^T$ will not necessarily come at low cost.  

Constraint preconditioning requires the factorization of the preconditioner $P_\text{CP}$ at each IPM step. 
In recent years, some attempts have been made to reduce the cost of factorizing preconditioners for $2 \times 2$ block matrices similar to~\eqref{saddle-point-system} through incomplete factorizations or sparse approximate inverse techniques in~\cite{chow1997approximate,perugia2000block,toh2004block,haws2002preconditioning}.  
%@Samah: may say something about $AE^{-1}A^T$ having fill, but so does our method!
The indefiniteness of the preconditioner also prohibits using some Krylov subspace methods such as the minimal residual method (MINRES)~\cite{paige1975solution} which require a positive-definite preconditioner~\cite{benzi2008some}. 
Consequently CP is often used in conjunction with popular iterative methods for indefinite systems such as the biconjugate gradient stabilized method (BiCGSTAB)~\cite{van1992bi}, the generalized minimal residual method (GMRES)~\cite{saad1986gmres}, or a simplified variant of the quasi-minimal residual method (QMR)~\cite{freund1995software}. 

%\sout{Moreover within the context of primal-dual interior point methods, the (1,1) block of $\CMcal{A}$ in (\ref{saddle-point-system}) namely $G^{(k)}$ will change at every IPM iteration as previously seen in (\ref{augmented-2}). Further, since the diagonal of  $G^{(k)}$ changes, generally the matrix $E^{(k)}$ approximating $G^{(k)}$ in the constraint preconditioner (\ref{const-prec}) will also change, prompting the factorization of $P_\text{CP}$ at every iteration.  Thus using CPs to precondition KKT systems within interior point method may have a very high overall computational cost.}

\subsection{{Augmented Lagragngian preconditioners}}

{Block diagonal preconditioners have been used in the solution of saddle point linear systems resulting from a discretization of partial differential equations~\cite{silvester1994fast,wathen1995convergence,perugia2000block,krzyzanowski2001block,powell2003optimal,toh2004block,mardal2004uniform,greif2006preconditioners,greif2007preconditioners}. They have also been used in the solution of KKT systems arising in constrained linear programming problems (LPs)~\cite{gill1992preconditioners,bergamaschi2021new}, QPs~\cite{pearson2017fast,bergamaschi2021new}, as well as NLPs~\cite{lukvsan1998indefinitely}.
A particular family of block diagonal preconditioners that has been popular in optimization especially in more recent years are {\it augmented Lagrangian preconditioners},
\begin{equation}
\label{rees-greif}
    P_\text{RG} = \begin{bmatrix}
G^{(k)}+A^TW^{-1}A & 0 \\ 
0 & W 
\end{bmatrix}.
\end{equation}
%\sout{In~\cite{rees2007preconditioner}, Rees and Greif propose block preconditioners for augmented KKT systems~\eqref{augmented-2} given by an IPM formulation for a QP of the form~\eqref{qp} with both equality and inequality constraints.} 
This type of preconditioner has been originally proposed in \cite{greif2006preconditioners} and later extended in \cite{rees2007preconditioner,cao2008augmentation}. The main approach is to augment the $(1,1)$ block in the saddle point matrix $\CMcal{A}$~\eqref{saddle-point-system} using a symmetric positive definite (SPD) weight matrix $W \in \mathbb{R}^{m_1\times m_1}$, resulting in an SPD preconditioner $P_\text{AL}$. A common choice for the weight matrix is a a weighted identity $W = \gamma I$~\cite{golub2003solving,rees2007preconditioner}, where $\gamma$ is chosen so that the norm of the augmented term $\gamma^{-1}A^TA$ is comparable to that of $G^{(k)}$; for instance $\gamma = ||A||^2/||G^{(k)}||$ for QPs. In~\cite{rees2007preconditioner}, an augmented Lagrangian preconditioner of the form shown in~\eqref{rees-greif} has been applied to augmented KKT systems~\eqref{augmented-2} given by an IPM formulation for QPs with both equality and inequality constraints. Using spectral analysis, they show that the preconditioned matrix $P_\text{RG}^{-1}\CMcal{A}$ has eigenvalues $\lambda = 1$ with multiplicity $n$ and $\lambda = -1$ with multiplicity $p$, while the rest of the eigenvalues lie in the interval $(-1,0)$. 
These eigenvalues are more tightly clustered as the $(1,1)$ block of $\CMcal{A}$ becomes more ill-conditioned, which is known to happen when the interior point iterate approaches optimality~\cite{wright1997primal}.

The name of this preconditioner stems from augmented Lagrangian approaches that have been used in optimization when the $(1,1)$ block $G^{(k)}$ in the saddle point matrix $\CMcal{A}$ is ill-conditioned, or even singular~\cite{glowinski1989augmented,fortin2000augmented,golub2003solving,greif2004augmented,fletcher2013practical,benzi2008some}. The original saddle point system~\eqref{saddle-point-system} is replaced by an equivalent one, with the same solution, but having better numerical properties~\cite{golub2003solving},
\begin{equation*}
\begin{bmatrix}
G^{(k)}+A^TW^{-1}A & A^T \\ 
A & 0 
\end{bmatrix}
\begin{pmatrix}
x\\y
\end{pmatrix}=\begin{pmatrix}
b+A^TW^{-1}c\\c
\end{pmatrix}.
\end{equation*}
Further, augmented Lagrangian preconditioners have been applied to regularized saddle point systems where the $(2,2)$ block of $\CMcal{A}$ is non-zero~\cite{morini2016spectral,shen2012augmentation}, and similarly to nonsymmetric saddle point systems~\cite{farrell2019augmented,cao2008augmentation}. 
Block triangular preconditioners based on the augmented Lagrangian idea have also been proposed to solve saddle point systems in~\cite{benzi2008some,benzi2006augmented,rees2007preconditioner,cao2008augmentation,shen2012augmentation}.}

%In \cite{morini2016spectral}, they addressed the use of positive deﬁnite augmented preconditioners that preserve the symmetry of the coefﬁcient matrix, while coping with the possible high singularity of the diagonal blocks. General spectral bounds for the preconditioned matrix were derived; they signiﬁcantly expand results in the literature, covering the case of all nonzero diagonal blocks.

%@Samah: Maybe add this: But this ill-conditioning will increase the computational cost of computing the (1,1) block of P_\text{RG}.   

%@Samah: Maybe add this: This strong clustering of the eigenvalues will often promote the fast convergence of many Krylov subspace methods for symmetric, and in general normal matrices.  

%\sout{The $(1,1)$ block of the block diagonal preconditioner in (\ref{rees-greif}) will change at every primal-dual interior point iteration $k$, hence a factorization of the matrix $G^{(k)}+A^TW^{-1}A$ is required at each IPM iteration.}
%per iteration, which might lead to a high overall cost.

\section{Our approach: single factorization inexact interior point method}\label{sec:approach}
The aforementioned state-of-the-art preconditioning techniques for the augmented KKT system all require the factorization of the corresponding preconditioner at each IPM step. This step may make the inexact IPM algorithm prohibitively expensive given the relatively large size of the augmented matrix.
We propose a different reduced system, which reuses the same matrix factorization across all IPM steps. %\samah{Not sure about this. It says the other approaches have to factorize the preconditioners at every iteration while we reuse a factorization for the reduced system (not the preconditioners).}

%\sout{
%Evidently in the formulation \eqref{new-reduced-intro}, $F^{-1}$ refers to the application of the forward and backward substitution of the factors of $F$ rather than the explicit inverse. And since $F$ is a constant matrix for a given QP, then we only need to compute the Bunch-Parlett factorization of that matrix once prior to starting the IPM iterations. Then at every interior point iteration, forward and backward substitution algorithms can be applied to get the new reduced system \eqref{new-reduced-intro} in an implicit form. This approach has the advantage of reducing the computational cost per IPM iteration by reformulating the system so that the constant portion of the KKT system is factorized beforehand.% The details will be described in Section \ref{sec:approach}.       
%}

\subsection{New reduced system}

Consider the following splitting of $K^{(k)}$, the coefficient matrix of the first augmented KKT system in~\eqref{augmented}. %\sout{, where w}
We isolate the bottom right block, $D^{(k)}$, which will be the only sub-block changing at every interior point iteration,
\begin{equation}
\label{augmented-split}
%\underbrace{
\begin{bmatrix}[cc|c]
-H & A^T & C^T\\ 
A & 0 & 0\\\cmidrule[0.3pt]{1-3}
C & 0 & D^{(k)}
\end{bmatrix}%}_{K^{(k)}}
\begin{pmatrix}
\Delta x^{(k)}\\ 
\Delta \lambda^{(k)}\\
\Delta \nu^{(k)}
\end{pmatrix} =  -%\underbrace{
\begin{pmatrix}
r_g^{(k)}\\ 
r_e^{(k)}\\
r_a^{(k)}
\end{pmatrix}%}_{\bfr_1^{(k)}}
.
\end{equation}
%\sout{We denote the coefficient matrix of \eqref{augmented-split} by $K^{(k)}$.}\samah{I removed the underbrace $K^{(k)}$ and the previous sentence because we introduced this notation first in eq (2.10) }
%, $\bfr_1^{(k)}$ respectively. 
Let $F \in \bbR^{(n+m_1)\times(n+m_1)}$ be the top left block of $K^{(k)}$, %\samahnew{which we refer to as the KKT subsystem,}
\begin{equation}
	\label{f-def}
	F = \begin{bmatrix}
		-H & A^T\\
		A & 0
	\end{bmatrix}.
\end{equation}
%Then we can write the system \eqref{augmented-split} as 2 separate equations:
%\begin{subequations}
%\begin{align}	
%\label{sep-1}
%	F\begin{pmatrix}\Delta x^{(k)}\\ 
%\Delta \lambda^{(k)}\end{pmatrix} +  \begin{bmatrix}
%C^T\\ 0\end{bmatrix} \Delta \nu^{(k)} &=  -\begin{pmatrix}
%r_g^{(k)}\\ r_e^{(k)}\end{pmatrix},\\ \label{sep-2}
%\begin{bmatrix}C & 0\end{bmatrix} \begin{pmatrix}
%\Delta x^{(k)}\\ \Delta \lambda^{(k)}\end{pmatrix} + D^{(k)} \Delta \nu^{(k)} &= -r_a^{(k)}.
%\end{align}
%\end{subequations}
%\begin{equation}\label{schur-R}
%	K_F^{(k)} = D^{(k)}  - \begin{bmatrix}
%C & 0\end{bmatrix} F^{-1} \begin{bmatrix}
%C^T\\ 0\end{bmatrix}.
%\end{equation}
$F$ is %\sout{also a saddle point matrix} \samahnew{
non-singular, and its  invertibility can be readily established for a convex QP such as in~\eqref{qp} assuming $H$ is symmetric definite and $\rank(A)=m_1$~\cite[Corollary 3.1]{gansterer2003mathematical}. 
%Even though this result is well-known, we will state it here for completeness. If $H$ is symmetric definite and $\rank(A)=m_1$, then $F$ defined in \eqref{f-def} is nonsingular \cite[Corollary 3.1]{gansterer2003mathematical}.
%While if $H$ is symmetric semidefinite then $F$ defined in \eqref{f-def} is nonsingular if and only if $\rank(A)=m_1$, and $\rank\begin{pmatrix} H\\ A\end{pmatrix}=n$  \cite[Lemma 3.2]{gansterer2003mathematical}.
%if -H is symmetric positive definite and Null($H$) $\cap$ Null($A$) = $\{0\}$ \cite{benzi2005numerical} where Null() is the nullspace of a matrix?  
%\sout{It is straightforward to show that }
%\edgarnew{
Note that $F$, also a saddle point matrix, is symmetric and indefinite.
Assuming $H$ is positive-definite, by the Haynesworth inertia additivity formula~\cite{haynsworth1968inertia}, the inertia of $F$ is given by the inertia of $-H$ ($n$ negative eigenvalues) combined with the inertia of the Schur complement $AH^{-1}A^T$ ($m_1$ positive eigenvalues).
%}
%\sout{since $F = F^T$,and $e_{n+1}^TF e_{n+1}=0$ where $e_i$ is the $i$th unit vector of $\mathbb{R}^{n+m_1}$. In fact, by Sylvester's law of inertia \cite{horn2012matrix}, $F$ has $m_1$ positive and $n$ negative eigenvalues considering that the inertia of $F$ is equal to the sum of the inertias of $-H$ and the schur complement of $-H$ in $F$. }

 Using the splitting introduced in~\eqref{augmented-split}, we can reduce the first augmented system into a smaller %\sout{$m_2 \times m_2$}
 linear system of the form,

\begin{equation}
\label{ineq-const-reduced}
\begin{aligned}
  \underbrace{  \left[D^{(k)}  - \begin{bmatrix}
C & 0
\end{bmatrix} F^{-1} \begin{bmatrix}
C^T\\ 
0
\end{bmatrix} \right]}_{K_F^{(k)}} \Delta \nu^{(k)} = \underbrace{-r_a^{(k)} + \begin{bmatrix}
C & 0
\end{bmatrix} F^{-1}\begin{pmatrix}
r_g^{(k)}\\ 
r_e^{(k)}
\end{pmatrix}}_{r_\nu^{(k)}}.%{\beta^{(k)}}.
\end{aligned}
\end{equation}
where $K_F^{(k)}$ $\in \bbR^{m_2\times m_2}$ is the Schur complement of $F$ in $K^{(k)}$. We will denote this new reduced linear system~\eqref{ineq-const-reduced} as the {\it inequality-constraint reduced} system, %\samahnew{
as it solves for the step in the Lagrange multipliers corresponding to the inequality constraints. %}.
%It is straightforward to show that $K_F^{(k)}$ is indeed symmetric positive definite, and this allows the use of Cholesky factorisation if one chooses to use a direct solver or the Conjugate Gradient method \cite{hestenes1952methods} if one chooses to use an iterative solver. But 
The main advantage in adopting the new reduced system~\eqref{ineq-const-reduced} is the ability to factorize $F$ using a symmetric indefinite $LDL^T$ factorization in the {\it setup} phase of the interior point method, as  %\samahnew{$F$} 
this saddle point matrix does not change with iteration number $k$, and then use that factorization at every iteration of the IPM {\it solve} phase to compute the coefficient matrix, $K_F^{(k)}$, and the right-hand-side vector, $r_\nu^{(k)}$, of the new reduced system either implicitly or explicitly.  
%The right-hand side of the new reduced system \eqref{ineq-const-reduced} can be expressed as
%\begin{equation}
%\label{rhs-schur-nu}
%\begin{aligned}
%\beta^{(k)} &= -r_a^{(k)} + \begin{bmatrix}
%C & 0\end{bmatrix} F^{-1}\begin{pmatrix}
%r_g^{(k)}\\ r_e^{(k)}\end{pmatrix}
%= -r_a^{(k)} + C %\left(\tilde{F}_{11}r_g^{(k)} + %\tilde{F}_{12}r_e^{(k)}\right)\\
%	&= -r_a^{(k)} + C %\left[(-H^{-1}+H^{-1}A^TF_H^{-1}AH^{-1})r_g^{(k)} + %(H^{-1}A^TF_H^{-1}) r_e^{(k)}\right].
%\end{aligned}
%\end{equation}

Once the inequality-constraint reduced system system $K_F^{(k)} \Delta \nu^{(k)} = r_\nu^{(k)}$ is solved for the search direction $\Delta \nu^{(k)}$, the other variables in the first augmented KKT system~\eqref{augmented-split} can be recovered in the following manner,
\begin{equation}
\label{schur-delta-x}
\begin{aligned}
	F\begin{pmatrix}
\Delta x^{(k)}\\ 
\Delta \lambda^{(k)}
\end{pmatrix}  %&= -F^{-1} \left( \begin{pmatrix}
%r_g^{(k)}\\ r_e^{(k)}
%\end{pmatrix} + \begin{bmatrix}
%C^T\\ 0\end{bmatrix} \Delta \nu^{(k)}\right) \\
&=  - \begin{pmatrix}
r_g^{(k)}+C^T \Delta \nu^{(k)}\\ 
r_e^{(k)}
\end{pmatrix}.
%&=  -\begin{pmatrix}
%\tilde{F}_{11}(r_g^{(k)}+C^T \Delta \nu^{(k)}) + \tilde{F}_{12} r_e^{(k)}\\ 
%\tilde{F}_{12}^T(r_g^{(k)}+C^T \Delta \nu^{(k)}) + \tilde{F}_{22} r_e^{(k)}
%\end{pmatrix}. 
\end{aligned}
\end{equation}
Using the $LDL^T$ factorization of $F$ computed in the setup phase of IPM, a forward and backward substitution can be performed on the right-hand side vector of~\eqref{schur-delta-x} at every interior point iteration $k$ in order to obtain the search directions $(\Delta x^{(k)},\Delta \lambda^{(k)})$.

%Therefore solving the original KKT system \eqref{sym-indef-KKT} reduces to solving the Schur complemented system \eqref{schur-delta-s} for $\Delta s$ and then performing some matrix-vector multiplications using \eqref{schur-delta-x} to get $\Delta x$ and $\Delta \lambda$. Since the cost of solving a linear system is an order of magnitude greater than that of performing a mat-vec product, we will focus on making the factorization of $R$ more efficient. 

Although the inverse of $F$ in~\eqref{f-def} is not computed explicitly, we can further simplify the expression for $K_F^{(k)}$ in the inequality-constraint reduced system~\eqref{ineq-const-reduced} by using the following explicit expression for the inverse of a block $2 \times 2$ matrix,
\begin{equation}
\begin{aligned}
\label{f-inv}
	F^{-1} &= \Biggl[ \, \begin{matrix}
		%\overbrace{
		-H^{-1}+H^{-1}A^TF_H^{-1}AH^{-1}%}^{\tilde{F}_{11}} 
		& %\overbrace{
		H^{-1}A^TF_H^{-1}%}^{\tilde{F}_{12}}
		\\%\underbrace{
		F_H^{-1}AH^{-1}%}_{\tilde{F}_{12}^T}
		& %\underbrace{
		F_H^{-1}%}_{\tilde{F}_{22}}
	\end{matrix} \, \Biggr]
% = \begin{bmatrix}	\tilde{F}_{11} & \tilde{F}_{12}\\
%		\tilde{F}_{12}^T & \tilde{F}_{22}	\end{bmatrix}.
	\end{aligned}
\end{equation}
where $F_H = AH^{-1}A^T$, the Schur complement of $-H$ in $F$, is symmetric positive-definite~\cite{benzi2005numerical}. Consequently,
\begin{equation}
%\label{schur-R}
\label{schur-R}
\begin{aligned}
	K_F^{(k)} &= D^{(k)} - \begin{bmatrix}
		C & 0
	\end{bmatrix} F^{-1} \begin{bmatrix}
		C^T\\
		0
	\end{bmatrix}\\
	%= D^{(k)} - C \tilde{F}_{11} C^T\\
	&= D^{(k)} + C (H^{-1}-H^{-1}A^TF_H^{-1}AH^{-1}) C^T\\
	&= D^{(k)} + C H^{-1}(H-\underbrace{A^TF_H^{-1}A}_{\bar{H}_A})H^{-1} C^T.
	\end{aligned}
\end{equation}
$K_F^{(k)}$ is symmetric positive-definite, because it is the sum of a diagonal positive matrix ($D^{(k)}$) and a symmetric semi-definite matrix ($CH^{-1}(H-\bar{H}_A)H^{-1}C^T$). This allows for the use of the Cholesky factorization %\sout{if one chooses to use} \samahnew{
when opting for a direct solver or the conjugate gradient method (CG)~\cite{hestenes1952methods} when opting for an iterative solver. However, iterative methods have better scaling and produce less fill than their direct counterparts. In addition, an iterative solver, such as CG, %\sout{is}
would be especially useful in the solution of system~\eqref{ineq-const-reduced} since it does not require the explicit formation of the coefficient matrix $K_F^{(k)}$. 
To accelerate the convergence of CG, we introduce two preconditioners for the inequality-constraint reduced system in the next subsection.

\subsection{Preconditioners}\label{sec:pre}
In order to find a good preconditioner for the new reduced system~\eqref{ineq-const-reduced}, we would like to 
find an easily invertible transformation that improves the conditioning of $K_F^{(k)}$.
%approximate $K_F^{(k)}$ by a simpler form that is inexpensive to construct and ``invert''.
It is well-known that the diagonal matrix $D^{(k)}$ in~\eqref{ineq-const-reduced}, becomes ill-conditioned as $k$ increases~\cite{oliveira2005new}. Further, due to the structure of the IPM algorithm, the entries of $D^{(k)}$ are necessarily strictly positive (nonsingular) at every iteration. Therefore, we use the matrix $D^{(k)}$ in order to design suitable preconditioners for the new reduced system.  
\begin{comment}
Although the inverse of $F$ introduced in~\eqref{f-def} is not computed explicitly, we can simplify the expression for $K_F^{(k)}$ in the inequality-constraint reduced system~\eqref{ineq-const-reduced} by using the following explicit expression for the inverse of a block $2 \times 2$ matrix,
\begin{equation}
\begin{aligned}\label{f-inv}
	F^{-1} &= \Biggl[ \, \begin{matrix}
		-H^{-1}+H^{-1}A^TF_H^{-1}AH^{-1}
		& H^{-1}A^TF_H^{-1}\\
		F_H^{-1}AH^{-1}
		& F_H^{-1}
	\end{matrix} \, \Biggr]
	\end{aligned}
\end{equation}
where $F_H = AH^{-1}A^T$, the Schur complement of $-H$ in $F$, is symmetric positive-definite~\cite{benzi2005numerical}. Consequently,
\begin{equation}\label{schur-R}
\begin{aligned}
	K_F^{(k)} &= D^{(k)} - \begin{bmatrix}
		C & 0
	\end{bmatrix} F^{-1} \begin{bmatrix}
		C^T\\0
	\end{bmatrix}\\
	&= D^{(k)} + C (H^{-1}-H^{-1}A^TF_H^{-1}AH^{-1} C^T\\
	&= D^{(k)} + C H^{-1}(H-\underbrace{A^TF_H^{-1}A}_{\bar{H}_A})H^{-1} C^T.
	\end{aligned}
\end{equation}
\end{comment}

 We define the quantity $n_d = n - m_1 \geq 0$ as the {\it number of degrees of freedom} since $n \geq m_1$. 
 We differentiate between a {\it low-degree-of-freedom} case $(n_d \approx 0)$ when the number of equality constraints in the QP \eqref{qp} is approximately equal to that of the primal variables $(m_1 \approx n)$, and a {\it high-degree-of-freedom} case $(n_d \approx n)$ when the number of equality constraints is much less than that of the primal variables $(m_1 \ll n)$.
Since $F_H = AH^{-1}A^T$, we use the following approximations,
 \begin{equation}
 \begin{aligned}
 	\bar{H}_A &= A^TF_H^{-1}A = A^T(AH^{-1}A^T)^{-1}A,\\
 	&\approx \left\{\begin{matrix}
 H,& & n_d \approx 0,\\ % m_1 \approx n.\\  
 0,& &n_d \approx n.%m_1 << n.
\end{matrix}\right. 	\end{aligned}
 \end{equation}
Based on the approximation $\bar{H}_A = H$ when $n_d \approx 0$, we propose the following low-degree-of-freedom (low-d.o.f.) preconditioner for the inequality-constraint reduced matrix $K_F^{(k)}$~\eqref{schur-R},
\begin{equation}
\label{preconditioner-1}
\begin{aligned}
	P_L = D^{(k)}. %$m_1 \approx n$
\end{aligned}
\end{equation}
$P_L$ %\sout{$\in \bbR^{m_2 \times m_2}$} \samah{(removing size, it's redundant since we know the sice of $K_F$)} 
is a diagonal matrix with positive diagonal entries due to the strict positivity requirement~\eqref{strict-pos} of primal-dual interior point methods ($s,\nu > 0$). As a result, the cost of using this low-d.o.f. preconditioner within the conjugate gradient method is negligible. 
  
Additionally, we propose a high-degree-of-freedom (high-d.o.f.) preconditioner, 
\begin{equation}
\label{preconditioner-2}
\begin{aligned}
	P_H = D^{(k)} + CH^{-1}C^T, %$m_1 << n$
\end{aligned}
\end{equation}
for the new reduced matrix $K_F^{(k)}$~\eqref{schur-R}.
This preconditioner can be derived using the approximation $\bar{H}_A = 0$ when $n_d \approx n$. $P_H$ %\sout{$\in \bbR^{m_2 \times m_2}$} 
is a symmetric positive-definite matrix since it is equal to the sum of a positive diagonal matrix ($D^{(k)}$) and a symmetric semi-definite matrix ($CH^{-1}C^T$).  

The cost of applying this high-d.o.f. preconditioner depends on the cost of factorizing $H$, and that of performing a forward/backward substitution on the columns of $C$. In practice, Hessian matrices $H$ 
(or approximations thereof)
often have special sparsity structures, for instance diagonal or block diagonal, hence the cost of factorizing $H$ is tractable. In other cases, the preconditioner can employ an approximation of $H$ by a positive diagonal matrix, as we do in Section~\ref{sec:numerics}. Further, many constrained QP problems only have simple bounds on the variables, in which case $C = I$. The high-d.o.f. preconditioner can then be expressed as $P_H = D^{(k)} + H^{-1}$.   
%And a preconditioner for \eqref{ineq-const-reduced} is given by:   
%\begin{equation}
%	M^{-1} = (V^{-1}S + CH^{-1}C^T)^{-1}
%\end{equation}

\subsection{Spectral analysis of preconditioned matrices}

We now characterize the spectrum of the preconditioned matrices, which enable bounds on CG iteration counts.
To achieve this, we first provide two lemmas that characterize the ranks of the Schur complement terms $\bar{H}_A = A^TF_H^{-1}A$ and $H-\bar{H}_A$, which arise in our inexact IPM approach.
In the remainder of this section, we omit the superscript $(k)$, which denotes the IPM iteration number, for readability purposes.

\begin{lemma}\label{lem_s_a}
Assume $H \in \mathbb{R}^{n\times n}$ is symmetric positive-definite and $A \in \mathbb{R}^{m_1\times n}$ is full rank with $m_1 \leq n$.
%$\rank=m_1 \leq n$.
Define %\sout{$\bar{H}_A = A^T(AH^{-1}A^T)^{-1}A \in \mathbb{R}^{n\times n}$}
%\sout{$F_H = AH^{-1}A^T$, and $\bar{H}_A = A^TF_H^{-1}A$,}
%Given symmetric positive definite $H \in \mathbb{R}^{n\times n}$, full rank $A \in \mathbb{R}^{m_1\times n}$ with $m_1 \leq n$, define $\bar{H}_A = A^T(AH^{-1}A^T)^{-1}A \in \mathbb{R}^{n\times n}$ 
%is rank-deficient, with \sout{then the $\rank(\bar{H}_A)$ is less than or equal to $m_1$.}
$\bar{H}_A = A^T(AH^{-1}A^T)^{-1}A$, then $\rank(\bar{H}_A)\leq m_1$.
\end{lemma}  
{
\begin{proof}
By the properties of the rank of a product of matrices, we have 
\begin{equation}
\begin{aligned}
\rank(\bar{H}_A) = \rank\bigl(A^T(AH^{-1}A^T)^{-1}A\bigr) &\leq \min\bigl(\rank(A), \rank(AH^{-1}A^T)\bigr)\\
 	&= \min(m_1, n)
 	= m_1,%\;{\rm since }\;m_1 \leq n.
\end{aligned}
\end{equation}
 since $m_1 \leq n$, and $AH^{-1}A^T$ is SPD.
\end{proof}

}

\begin{lemma}\label{lem_h-s_a}
Assume $H \in \mathbb{R}^{n\times n}$ is symmetric positive-definite and $A \in \mathbb{R}^{m_1\times n}$ is full rank with $m_1 \leq n$.
%$\rank=m_1 \leq n$.
Define $\bar{H}_A = A^T(AH^{-1}A^T)^{-1}A$ 
then $\rank(H - \bar{H}_A)\leq n-m_1 $.
\end{lemma}  
\begin{proof}
To derive the rank upper-bound, we show that the columns of $H^{-1}A\in \mathbb{R}^{n\times m_1}$ form a basis for the nullspace of $H-\bar{H}_A$. Let $y$ be an arbitrary vector in the $\textrm{span}\{H^{-1}A\}$, $y = H^{-1}A z$. The following shows that $y \in \textrm{ker}(H-\bar{H}_A)$,
\begin{equation*}
\begin{aligned}
    (H-\bar{H}_A)y &= \Big(H - A^T\bigl(AH^{-1}A^T\bigr)^{-1}A\Big)y\\
    &= \Big(H - A^T\bigl(AH^{-1}A^T\bigr)^{-1}A\Big)H^{-1}A z\\
    &= (A-A)z = 0.
\end{aligned}
\end{equation*}
Thus we can deduce that the dimension of the nullspace of $(H-\bar{H}_A)$ is at least equal to $m_1$, the number of columns in $H^{-1}A$,
\begin{equation} \label{dim-ker-h-a}
    \textrm{dim}\left(\textrm{ker}(H-\bar{H}_A)\right) \geq m_1.
\end{equation}
The bound in the Lemma follows from the rank-nullity theorem, $\rank(H-\bar{H}_A) +  \textrm{dim}\left(\textrm{ker}(H-\bar{H}_A)\right) = n$.
%), and the bound on the nullity of $(H-\bar{H}_A)$ in~\ref{dim-ker-h-a} we can now show that 
%\begin{equation*} 
%\begin{aligned}
%    \textrm{dim}\left(\textrm{ker}(H-\bar{H}_A)\right)+\rank(H-\bar{H}_A) &\geq m_1 + \rank(H-\bar{H}_A),\\
%   % n &\geq m_1 + \rank(H-\bar{H}_A)\\
%   % n-m_1 &\geq \rank(H-\bar{H}_A)
%   \rank(H-\bar{H}_A)) \leq n - m_1.
%   \end{aligned}
%\end{equation*}
\end{proof}

 We want to describe the spectrum of the preconditioned matrix $ P^{-\frac{1}{2}}K_F P^{-\frac{1}{2}}$ to understand the convergence of CG. 
 We consider this split preconditioning as opposed to left-preconditioning ($P^{-1}K_F$) in order to preserve symmetry and simplify analysis.

\begin{theorem}\label{thm1}
  Consider any symmetric positive-definite matrix $H \in \mathbb{R}^{n\times n}$, any full rank matrix $A \in \mathbb{R}^{m_1\times n}$ with $m_1\leq n$,  %row $\rank (= m_1 \leq n)$,
  and $C \in \mathbb{R}^{m_2\times n}$. %where $m_2 \geq n$. 
  Let $D \in \mathbb{R}^{m_2\times m_2}$ be a diagonal matrix, $D = \diag(d_1,\ldots,d_{m_2})$, where $d_i \neq 0$ for $i=1,\ldots,m_2$. Further, define $K_F \in \mathbb{R}^{m_2\times m_2}$  as
  \begin{equation*}K_F = D - \begin{bmatrix}C & 0
	\end{bmatrix} F^{-1}\begin{bmatrix}
		C^T\\0\end{bmatrix},\end{equation*}
where $F = \begin{bmatrix}
		-H & A^T\\A & 0\end{bmatrix}$ %\sout{$\in \mathbb{R}^{(n+m_1)\times (n+m_1)}$} 
	is nonsingular. Preconditioning $K_F$ by a diagonal matrix of the form \begin{equation*}P_L = D,
	\end{equation*}
	%when $m_1 \approx n$, 
	implies that the preconditioned matrix $P_L^{-\frac{1}{2}}K_F P_L^{-\frac{1}{2}}$ has
\begin{enumerate}[label=(\arabic*)]
\item at most $n-m_1$ eigenvalues of the form $1 + \lambda_i$, where $\lambda_i$   is an eigenvalue of  $X = P_L^{-\frac{1}{2}}CH^{-1} \left( H - \bar{H}_A\right)H^{-1}C^TP_L^{-\frac{1}{2}}$, with $\bar{H}_A = A^T(AH^{-1}A^T)^{-1}A$, for $i = 1,\ldots,n-m_1$, 
and 
	\item 
	eigenvalue $1$ with multiplicity equal to $\max(0,m_2 - (n-m_1) )$. 
	%the remainder are unit eigenvalues. %an eigenvalue $1$ with multiplicity $m_2 - \min (n, 2(n-m_1) )$ 
\end{enumerate}
\end{theorem}
\begin{proof} 
Using~\eqref{schur-R} % and~\eqref{S_A_expression}, we can write \begin{align*}K_F = D + C H^{-1}\left( H - \bar{H}_A\right)H^{-1} C^T.\end{align*}
and since $P_L=D$, we can write the preconditioned matrix as
\begin{equation}\label{prec-mat-1}
\begin{aligned}
	P_L^{-\frac{1}{2}}K_F P_L^{-\frac{1}{2}} 
	&= P_L^{-\frac{1}{2}}\left[D + C H^{-1}\left( H - \bar{H}_A\right)H^{-1} C^T\right] P_L^{-\frac{1}{2}}\\
	&= I + \underbrace{P_L^{-\frac{1}{2}}CH^{-1} \left( H - \bar{H}_A\right) H^{-1}C^TP_L^{-\frac{1}{2}}}_{X},
	\end{aligned}
\end{equation}
where $\bar{H}_A = A^T(AH^{-1}A^T)^{-1}A^T$. The preconditioned matrix $P_L^{-\frac{1}{2}}K_F P_L^{-\frac{1}{2}}$ is a sum of the identity matrix and  a second term, $X= P_L^{-\frac{1}{2}}CH^{-1} \left( H - \bar{H}_A\right) H^{-1}C^TP_L^{-\frac{1}{2}}$. Using~\cref{lem_h-s_a}, we find an upper bound on the rank of the matrix $X$
 \begin{align*}
 	\rank(X) &\leq \min\bigl( \rank(P_L), \rank(C), \rank (H), \rank(H-\bar{H}_A)\bigr)\\ 
 	&= \min(m_2, n, n-m_1)\\
 	&\leq (n-m_1).%,\;\;\;\;\;\;\;\;{\rm since }\;n \leq m_2.
 \end{align*}
 % Using the assumption $m_1 \approx n$, we find that $(n-m_1)$ is a small integer approximately equal to zero. Hence $X$ is a low rank matrix, and $\rank(X) \leq 2(n-m_1) \approx 0$. Thus the preconditioned matrix \eqref{prec-mat-1}) has an eigenvalue 1 with multiplicity $m_2-2(n-m_1) \approx m_2$, and $2(n-m_1) \approx 0$ eigenvalues of the form $1+\lambda_i$, where $\lambda_i$ is an eigenvalue of $X = P_L^{-\frac{1}{2}}CH^{-1} ( H - QQ^THQQ^T)H^{-1}C^TP_L^{-\frac{1}{2}}$, for $i = 1, \ldots, 2(n-m_1)$. 
Therefore, $X$ has at most $(n-m_1)$ non-zero eigenvalues $\lambda_i$.
%And since by \eqref{prec-mat-1}, the preconditioned matrix $P_L^{-\frac{1}{2}}K_F P_L^{-\frac{1}{2}}$ reduces to the sum $(I+X)$, 
Consequently, the matrix $P_L^{-\frac{1}{2}}K_F P_L^{-\frac{1}{2}}=I+X$ has at most $(n-m_1)$ non-unit eigenvalues of the form $1 + \lambda_i$, while the remainder of its eigenvalues are equal to $1$.
\end{proof}
% It is well-known that the eigenvalues of the sum of the identity matrix and an arbitrary matrix $B$ is given by: $\lambda(I + B) = 1 + \lambda(B)$. Therefore we know that the preconditioned matrix \eqref{prec-mat-1} has $m_2$ eigenvalues equal to $1 + \lambda(X)$.
 \begin{corollary}\label{cor:cg1} The conjugate gradient method applied to~\eqref{ineq-const-reduced} with preconditioner $P_L$~\eqref{prec-mat-1} 
converges %\sout{to the solution when executed in exact arithmetic} 
in at most $(n-m_1)+1$ iterations, when executed in exact arithmetic. 
%\edgar{Abbreviation 'PCG' has not been introduced I think, also should say what preconditioner is used here.}
\end{corollary}
\begin{proof}
By~\cref{thm1}, %\sout{the preconditioned matrix}
$P_L^{-1/2}K_FP_L^{-1/2}$ has an eigenvalue 1 with multiplicity $k = \max(0,m_2 -(n-m_1) )$.
Consequently, the preconditioned matrix has at most $m_2-(k-1)=(n-m_1)+1$ unique eigenvalues.
Given %\sout{a zero vector as} 
an initial guess, conjugate gradient converges after a number of iterations equal to the number of unique eigenvalues of the preconditioned matrix~\cite{golub2013matrix}.
\end{proof}

\begin{theorem}\label{thm2}
 Consider any symmetric positive-definite matrix $H \in \mathbb{R}^{n\times n}$, any full rank matrix $A \in \mathbb{R}^{m_1\times n}$ with $m_1\leq n$, and
 %row $\rank (= m_1 \leq n)$ with $m_1\leqn$,
 $C \in \mathbb{R}^{m_2\times n}$. %where $m_2 \geq n$. 
 Let $D \in \mathbb{R}^{m_2\times m_2}$ be a diagonal matrix, $D = \diag(d_1,\ldots,d_{m_2})$, where $d_i \neq 0$ for $i=1,\ldots,m_2$. Further, define  $K_F \in \mathbb{R}^{m_2\times m_2}$ as
  \begin{equation*}
      K_F = D - \begin{bmatrix}C & 0
	\end{bmatrix} F^{-1} \begin{bmatrix}C^T\\0\end{bmatrix},
  \end{equation*}
where $F = \begin{bmatrix}-H & A^T\\A & 0
	\end{bmatrix}$ %\sout{$\in \mathbb{R}^{(n+m_1)\times (n+m_1)}$} 
	is nonsingular. Preconditioning $K_F$ by a symmetric positive-definite matrix of the form, 
	\begin{equation*}P_H = D + CH^{-1}C^T,
	\end{equation*}
	%when $m_1 << n$, 
	implies that the preconditioned matrix $P_H^{-\frac{1}{2}}K_F P_H^{-\frac{1}{2}}$ has
\begin{enumerate}[label=(\arabic*)]
\item at most $m_1$ eigenvalues of the form $1 - \lambda_i$, where
	$\lambda_i$ is an eigenvalue of $Y = P_H^{-\frac{1}{2}}CH^{-1} \bar{H}_AH^{-1}C^TP_H^{-\frac{1}{2}}$, with $\bar{H}_A=A^T(AH^{-1}A^T)^{-1}A$, for $i = 1,\ldots,m_1$, and 
	\item 
	%the remainder are unit eigenvalues. %an \edgarnew{
	eigenvalue $1$ with multiplicity $\max(0,m_2 - m_1)$. 
\end{enumerate}
\end{theorem}

\begin{proof} 
Using~\eqref{schur-R} and $P_H = D + CH^{-1}C^T$, we can write the preconditioned matrix as %$P_H^{-\frac{1}{2}}K_F P_H^{-\frac{1}{2}}$ as
\begin{equation}\label{prec-mat-2}
\begin{aligned}
	P_H^{-\frac{1}{2}}K_F P_H^{-\frac{1}{2}} 
	    =& P_H^{-\frac{1}{2}}\left[D + C H^{-1}(H-\bar{H}_A)H^{-1} C^T\right] P_H^{-\frac{1}{2}}\\
%		=& P_H^{-\frac{1}{2}}\left[(D + C H^{-1}C^T)-C(H^{-1}\bar{H}_AH^{-1})C^T\right] P_H^{-\frac{1}{2}}\\
%		=& P_H^{-\frac{1}{2}}(D + C H^{-1}C^T)P_H^{-\frac{1}{2}} 	- P_H^{-\frac{1}{2}}CH^{-1} \bar{H}_AH^{-1}C^TP_H^{-\frac{1}{2}} \\
		=& I - \underbrace{P_H^{-\frac{1}{2}}CH^{-1} \bar{H}_AH^{-1}C^TP_H^{-\frac{1}{2}}}_{Y} .\\
	\end{aligned}
\end{equation}
The preconditioned matrix $P_H^{-\frac{1}{2}}K_F P_H^{-\frac{1}{2}}$ is equal to the difference between the identity matrix and a second term, $Y=P_H^{-\frac{1}{2}}CH^{-1} \bar{H}_AH^{-1}C^TP_H^{-\frac{1}{2}}$.
Using~\cref{lem_s_a}, we find an upper bound on the rank of the matrix $Y$
 \begin{align*}
 	\rank(Y) &\leq \min\bigl( \rank(P_H), \rank(C), \rank (H), \rank(\bar{H}_A)\bigr)\\ 
 	&= \min(m_2, n, m_1) \leq m_1.
 	%m_1,\;\;\;\;\;\;\;\;{\rm since }\;m_1 \leq n.% \leq m_2.\\
 	%&= 2(n-m),\;\;\;\;\;\;\;\;{\rm since }\;m \approx n
 \end{align*}
 Therefore, $Y$ has $m_1$ %$\leq m_2$ 
 non-zero eigenvalues $\lambda_i$. 
 Consequently, the preconditioned matrix $P_H^{-\frac{1}{2}}K_F P_H^{-\frac{1}{2}}=I-Y$ has at most $m_1$ non-unit eigenvalues of the form $1 - \lambda_i$, while %\sout{the remainder of its} \samahnew{
 its remaining eigenvalues are equal to $1$.
\end{proof}
 %Let us consider two cases:
 %\begin{enumerate}[label=(\roman{*})]
  %   \item $\min(m_2, m_1) = m_1$, i.e. $m_1 \leq m_2$, then $Y$ is rank deficient and its $\rank(Y)$ is bounded by $m_1$. The preconditioned matrix \eqref{prec-mat-2} has an eigenvalue 1 with multiplicity $m_2-m_1$, and $m_1$ eigenvalues of the form $1-\lambda(Y)$.
     %\item $\min(m_2, m_1) = m_2$, then $Y$ is a full rank matrix. The preconditioned matrix \eqref{prec-mat-2} has $m_2$ eigenvalues of the form $1-\lambda(Y)$. 
 %\end{enumerate}
%It is well-know that the eigenvalues of $\lambda(I - B) = 1 - \lambda(B)$ for an arbitrary matrix $B$. 
\begin{corollary}\label{cor:cg2}
%\sout{PCG} 
%\edgar{same comment as above} \samahnew{The Preconditioned Conjugate Gradient (PCG) 
The conjugate gradient method applied to~\eqref{ineq-const-reduced} with preconditioner $P_H$~\eqref{prec-mat-2},
 converges %\sout{to the solution when executed in exact arithmetic} 
 after at most $m_1+1$ iterations, when executed in exact arithmetic.
\end{corollary}
The correctness of~\cref{cor:cg2} follows by the same argument as the proof of~\cref{cor:cg1}.
%\samahnew{
Note that these upper bounds on iteration count also apply when using conjugate gradient with one-sided preconditioning ($P^{-1}K_F$), since it is equivalent to executing CG on $P^{-1/2}K_FP^{-1/2}$~\cite{saad2003iterative}.%} \sout{These iteration count bounds also apply when using CG with one-sided preconditioning ($P^{-1}K_F$), since it is equivalent to executing CG on $P^{-1/2}K_FP^{-1/2}$~\cite{saad2003iterative}.}
%\begin{proof}
%Make an argument here about the number of CG iterations being bounded by the number of distinct eigenvalues of the preconditioned matrix. See Thm 10.2.5 in "Matrix Computations" by Gene Golub for proof of this.
%\end{proof} 
%Using the assumption $m_1 \ll n$, we find that $Y$ is a low rank matrix, and $\rank(Y) \leq m_1 \approx 1$. Thus the preconditioned matrix \eqref{prec-mat-2} has an eigenvalue 1 with multiplicity $m_2-m_1$, and $m_1 \approx 1$ eigenvalues of the form $1-\lambda_i$, where $\lambda_i$ is an eigenvalue of $Y = P_H^{-\frac{1}{2}}CH^{-1} (QQ^THQQ^T)H^{-1}C^TP_H^{-\frac{1}{2}}$, for $i = 1, \ldots, m_1$. 

\section{Performance models}

%\sout{In this section, we develop performance models for IPM with our new inequality-constraint reduced system formulation solved with an iterative solver preconditioned with our low and high degree of freedom preconditioners and we compare those to other state-of-the art IPM solvers with an augmented KKT system formulation along with multiple preconditioners.}\edgar{run-on!}
%\samahnew{
In this section, we develop performance models for our single factorization inexact IPM algorithm, which uses the preconditioned conjugate gradient method (PCG) with one of the preconditioners introduced in~\cref{sec:pre}.
We then compare these predicted costs to other popular inexact IPM algorithms that solve the augmented system~\eqref{augmented-2} using an iterative solver with multiple state-of-the art preconditioners.
Our models approximate the total number of floating-point operations (flops) for 
%\edgarnew{
sparse matrix operations in %}
the interior point algorithm to leading order in  %\sout{the} 
input parameters $n, m_1, m_2$.    

Consider a QP of the form~\eqref{qp}, where the matrices $H$, $A$ and $C$ are sparse. We consider six variants of the IPM algorithm that differ in the KKT linear system solution step. First, we consider an IPM algorithm where we solve the augmented system~\eqref{augmented-2} either with a direct incomplete $LDL^T$ factorization, or with an iterative solver, namely BiCGSTAB, without any preconditioning. 
%\edgar{Introduce BiCGSTAB abbreviation, and give reference if possible.} Done in constraint preconditioning subsection
Additionally, we consider versions where the augmented system is solved with BiCGSTAB preconditioned using a constraint preconditioner such as in~\eqref{const-prec} as well as with a block diagonal preconditioner such as in~\eqref{rees-greif}. We compare the previous variants to an IPM algorithm where the inequality-constraint reduced system~\eqref{ineq-const-reduced} is solved using PCG, with our low-d.o.f. preconditioner~\eqref{prec-mat-1} and our high-d.o.f. preconditioner~\eqref{prec-mat-2}. In~\cref{ipm-variants}, we summarize the different variants and preconditioners used for our models.
%Justification for BiCGSTAB: If no good symmetric positive definite preconditioner is known, or if a very good nonsymmetric preconditioner is available, the possible advantage of A being symmetric is lost, and one usually must use a solution method for nonsymmetric matrices. See Simoncini (2004a)%
%\sout{First, let us start with some definitions and notations:}
Our cost analysis employs the following notation.
%\begin{itemize}
 %   \item[-] $\nz(X)$: number of non zeros in sparse matrix $X \in \bbR^{m\times n}$,
 %   \item[-] $X[:,j]$: the $j$-th column of $X$,
 %   \item[-] $L_X$: lower triangular factor obtained from a Cholesky or an $LDL^T$ factorization of $X$ if $X$ is SPD or symmetric indefinite respectively,
 %   \item[-] $n_\text{kr}$: number of iterations of the Krylov solver (BiCGSTAB/CG)
 %   \item[-] $N_I$: number of IPM iterations
 %   \item[-] $f(X) = 2 \sum\limits_{j=1}^{n}  \big[nz(X[:,j])\big]^2$
 %   \item[-] $t(X) = 2 \nz (X)$
 %   \item[-] $g(X) = 2 \sum\limits_{i=1}^{m}  \big[nz(X[i,:])\big]^2$
%\end{itemize}
\begin{equation*}
\begin{array}{ll}
nz(X) & \mbox{number of non-zeros in sparse matrix } X \in \bbR^{m\times n} \\ 
 X[:,j]& \mbox{the } j\mbox{-th column of }X\\ 
 L_X& 
 \mbox{lower triangular factor from Cholesky or (pivoted) }LDL^T \mbox{ factorization if}\\
 &X\mbox{ is SPD or symmetric indefinite respectively}\\ 
n_\text{kr} & \mbox{number of iterations of the Krylov solver (BiCGSTAB/CG)}\\ 
N_I & \mbox{number of IPM iterations}%\\ 
%f(X) & 2 \sum\limits_{j=1}^{n} \big[nz(X[:,j])\big]^2\\ 
%t(X) & 2 \nz (X)\\ 
%g(X) & 2 \sum\limits_{i=1}^{m} \big[nz(X[i,:])\big]^2
\end{array}
\end{equation*}
\begin{table}[]
\begin{tabular}{llll}
\hline
Variant & Solver & Preconditioner & Linear System \\ \hline
    D-KC     & $LDL^T$ & -  & Augmented            \\ \hline
   U-KC   &  BiCGSTAB & None & Augmented            \\ \hline
      CP-KC  &    BiCGSTAB & Constraint  & Augmented           \\ \hline
      RG-KC  &    BiCGSTAB & Augmented Lagrangian  & Augmented            \\ \hline
      U-KF &   CG & None & Inequality-constraint reduced           \\
      \hline
   PL-KF &   CG & Low-d.o.f.  & Inequality-constraint reduced          \\ \hline
    PH-KF &  CG & High-d.o.f.  & Inequality-constraint reduced           \\ \hline
\end{tabular}
\caption{IPM variants with direct and iterative (preconditioned) solvers considered in the paper, and the Schur complement of $K^{(k)}$ they operate on.}\label{ipm-variants}
\end{table}

In our models, we use five basic computational kernels, a factorization (fact) kernel whether $LDL^T$ or Cholesky, a triangular solve (trsv) kernel, sparse matrix vector product (spmv) kernel, sparse matrix times sparse matrix multiplication (spmm) kernel and a Frobenius norm (frob) kernel. We ignore vector operations and diagonal matrix scaling operations. 
We now define the costs of these kernels, given a sparse symmetric matrix $\mathpzc{A}\in \bbR^{n\times n}$
with triangular factor (from Cholesky or pivoted $LDL^T$ factorization) $L_\mathpzc{A}$, a sparse matrix $\mathpzc{B} \in \bbR^{m\times n}$, and diagonal matrices $\mathpzc{D_1} \in \bbR^{m\times m}$, $\mathpzc{D_2} \in \bbR^{n\times n}$,
\begin{align}
    &c_\text{fact}(\mathpzc{A}) = %2
    \sum\limits_{j=1}^{n}  \big[nz(L_{\mathpzc{A}}[:,j])\big]^2,\\% = f(L_{A}),\\
    &c_\text{trsv}(L_\mathpzc{A}) = 2 \nz(L_\mathpzc{A}),\\% = t(L_A),\\
    &c_\text{spmv}(\mathpzc{B}) = 2 \nz (\mathpzc{B}),\\% = t(L_A),\\
    &
    %\left\{\begin{matrix}
c_\text{spmm}(\mathpzc{B}^T,\mathpzc{D}_1\mathpzc{B})= %2
\sum \limits_{i=1}^{m}\big[nz(\mathpzc{B}[i,:])\big]^2,\\% = g(A),\\ 
&
c_\text{spmm}(\mathpzc{B},\mathpzc{D}_2\mathpzc{B}^T)= %2
\sum \limits_{j=1}^{n}\big[nz(\mathpzc{B}[:,j])\big]^2, \\% = f(A), \\
%\end{matrix}\right.\\
    &c_{norm}(\mathpzc{A}) = %2
    \nz (\mathpzc{A}).% = t(L_A).
\end{align}
For all the variants, at every IPM iteration we need to compute $[r_g,r_e,r_i,r_c]^T$ previously defined in~\eqref{rhs-KKT4}. We can do that with a cost of
\begin{equation}
\label{c-rhs-r4}
c_{rhs} = N_I[c_\text{spmv}(H)+2c_\text{spmv}(A)+2c_\text{spmv}(C)].
\end{equation}
Further, a common final step in all the variants is solving for the search direction $\Delta s$ using~\eqref{delz}, but we  ignore that cost since $D$ is a diagonal matrix.
The cost %\sout{of} \samahnew{
breakdown of all the variants is  summarized in~\cref{tab:costs}.
These cost expressions are derived in the following subsections.
%The cost~\eqref{c-rhs-r4} needs to be added to all of the costs of the multiple variants.

\subsection{Augmented system}\label{sec:aug-cost}

For the first four variants, in all of which the augmented~\eqref{augmented-2} system is solved, at every IPM iteration the coefficient matrix $K_C$ and the right-hand side vector %$\bfr_2$
need to be updated. We assume $r_a = r_i - Dr_c$ has a negligible cost, while $r_u = r_g-C^TD^{-1}r_a$ has cost,
\[
N_I \times c_\text{spmv}(C).
\]
Additionally we need to compute the $(1,1)$ block of $K_C$ using the expression in~\eqref{augmented-2} at every IPM iteration at a cost similar to that of computing the matrix-matrix product $C^TDC$ since $D$ is diagonal,
\[
N_I \times c_\text{spmm}(C^T,DC).
\]
After the linear system solve for these four variants, we need to compute the Newton step $\Delta \nu = -D^{-1}(r_a+C\Delta x)$, which has cost,
\[
N_I \times c_\text{spmv}(C).
\]

\subsubsection{Direct solve (D-KC)}

For the first variant, where the direct solve is used, at every IPM iteration an $LDL^T$ factorization of the augmented system~\eqref{augmented-2} is performed, followed by two triangular solves with the right-hand side with costs of $N_I [c_\text{fact}(K_C)]$ and $N_I[2 c_\text{trsv}(L_{K_C})]$, respectively. The total cost (including terms derived in~\cref{sec:aug-cost}) is
\begin{equation}
      N_I[c_\text{fact}(K_C) + 2 c_\text{trsv}(L_{K_C}) + 2 c_\text{spmv}(C) + c_\text{spmm}(C^T,DC) ].
\end{equation}
%\begin{equation}
%  \begin{aligned}
%      c_{v1} =& N_I[c_\text{fact}(K_C) + 2 c_\text{trsv}(L_{K_C}) + 2 c_\text{spmv}(C) + c_\text{spmm}(C^T,DC) ]\\
%      =& \underbrace{N_I c_\text{fact}(K_C)}_\text{fact} + \underbrace{2 N_I c_\text{trsv}(L_{K_C})}_\text{trsv} + \underbrace{2 N_I c_\text{spmv}(C)}_\text{spmv} + \underbrace{N_I c_\text{spmm}(C^T,DC)}_\text{spmm}.
%\\    \end{aligned}
%\end{equation}
\subsubsection{Unpreconditioned iterative solve (U-KC)}
For the second variant, the augmented system~\eqref{augmented-2} is solved using BiCGSTAB without any preconditioning. In this case, the main cost of the linear system solution comes from two sparse-matrix vector multiplications with $K_C$ at every BiCGSTAB iteration for a cost of $N_I [n_\text{kr} \times 2 c_\text{spmv}(K_C)]$. The total cost of the IPM algorithm is
\begin{equation}
\begin{aligned}
    %c_{v2} =
    &N_I[n_\text{kr} \times 2 c_\text{spmv}(K_C) + 2 c_\text{spmv}(C) + c_\text{spmm}(C^T,DC)  ].
%    =& \underbrace{2N_I [n_\text{kr} c_\text{spmv}(K_C) +  c_\text{spmv}(C)]}_\text{spmv} + \underbrace{N_I c_\text{spmm}(C^T,DC)}_\text{spmm}.  
    \end{aligned}
\end{equation}
\subsubsection{Constraint preconditioner (CP-KC)}

For the third variant, the augmented system is solved using preconditioned BiCGSTAB with a constraint preconditioner such as in~\eqref{const-prec}. The symmetric indefinite preconditioner $P_{CP}$ needs to be factorized using an $LDL^T$ factorization for a cost of $N_I c_\text{fact}(P_{CP})$. Then at every iteration of the Krylov solver, the main costs are due to four triangular solves (two sets of forward/backward substitutions) and two sparse-matrix vector multiplications with $K_C$, which have costs of $N_I[n_\text{kr}\times 4 c_\text{trsv}(L_{P_{CP}})]$ and $N_I [n_\text{kr} \times 2 c_\text{spmv}(K_C)]$, respectively. Therefore the total cost of the IPM method is
\begin{equation}
\begin{aligned}
    N_I\Big[& c_\text{fact}(P_{CP}) + n_\text{kr}  \big(4 c_\text{trsv}(L_{P_{CP}})\big) + n_\text{kr}  \big(2 c_\text{spmv}(K_C)\big)\\
    &+2 c_\text{spmv}(C) + c_\text{spmm}(C^T,DC)  \Big].
%    =& \underbrace{N_I c_\text{fact}(P_{CP})}_\text{fact} + \underbrace{N_I n_\text{kr}  \big(4 c_\text{trsv}(L_{P_{CP}})\big)}_\text{trsv} + \underbrace{2N_I \big[ n_\text{kr}  c_\text{spmv}(K_C) +  c_\text{spmv}(C) \big]}_\text{spmv}\\
%    &+ \underbrace{N_I c_\text{spmm}(C^T,DC)}_\text{spmm}.  
    \end{aligned}
\end{equation}
\subsubsection{Augmented Lagrangian (block diagonal) preconditioner (RG-KC)} 
For this fourth variant, the augmented system is solved using preconditioned BiCGSTAB with a block diagonal preconditioner such as in~\eqref{rees-greif}, with $W = \gamma I$ and $\gamma = ||A||^2/||G||$, where $G=H+C^TD^{-1}C$. Since the matrix $G$ is changing at every IPM iteration, its norm needs to be computed every time for a cost of $N_I \times c_{norm}(G)$. In order to form the preconditioner $P_{RG}$, a sparse matrix multiplication $A^TW^{-1}A$ has to be performed at a cost of $N_I \times c_\text{spmm}(A^T,D_WA)$, where $D_W=W^{-1}$ is diagonal.
%since $W$ is diagonal. 
The remaining cost is identical to that of variant CP-KC, resulting in a total %\sout{cost} 
of
\begin{equation}
\begin{aligned}
     N_I\Big[&c_\text{fact}(P_{RG}) + n_\text{kr}  \big(4 c_\text{trsv}(L_{P_{RG}})\big) + n_\text{kr}  \big(2 c_\text{spmv}(K_C)\big) \\
    &+ c_\text{spmm}(A^T,D_WA) + c_{norm}(G)  + 2 c_\text{spmv}(C) + c_\text{spmm}(C^T,DC) \Big].
%    =& \underbrace{N_I c_\text{fact}(P_{RG})}_\text{fact} + \underbrace{N_I n_\text{kr}  \big(4 c_\text{trsv}(L_{P_{RG}})\big)}_\text{trsv} + \underbrace{2N_I \big[ n_\text{kr}  c_\text{spmv}(K_C) +  c_\text{spmv}(C) \big]}_\text{spmv} \\
%    &+ \underbrace{N_I [c_\text{spmm}(C^T,DC)+ c_\text{spmm}(A^T,DA)]}_\text{spmm}  + N_I c_{norm}(G).  
    \end{aligned}
\end{equation}
%\begin{equation}
%\begin{aligned}
%    c_{v4} =& N_I\Big[c_\text{fact}(P_{RG}) + n_\text{kr}  \big(4 c_\text{trsv}(L_{P_{RG}})\big) + n_\text{kr}  \big(2 c_\text{spmv}(K_C)\big) \\
%    &+ c_\text{spmm}(A^T,DA) + c_{norm}(G)  + 2 c_\text{spmv}(C) + c_\text{spmm}(C^T,DC) \Big]\\
%    =& \underbrace{N_I c_\text{fact}(P_{RG})}_\text{fact} + \underbrace{N_I n_\text{kr}  \big(4 c_\text{trsv}(L_{P_{RG}})\big)}_\text{trsv} + \underbrace{2N_I \big[ n_\text{kr}  c_\text{spmv}(K_C) +  c_\text{spmv}(C) \big]}_\text{spmv} \\
%    &+ \underbrace{N_I [c_\text{spmm}(C^T,DC)+ c_\text{spmm}(A^T,DA)]}_\text{spmm}  + N_I c_{norm}(G).  
%    \end{aligned}
%\end{equation}
\subsection{Inequality-constraint reduced system}
For the PL-SF and PH-KF variants, we separate the cost into a setup phase cost and a solve phase cost. For the setup phase, an $LDL^T$ factorization of $F$ is performed at a cost equal to $c_\text{fact}(F)$.
At every IPM iteration, the right-hand side vector of the inequality-constraint reduced system~\eqref{ineq-const-reduced} needs to be updated. In particular, we assume  the computation of $r_a = r_i - Dr_c$ has a negligible cost. We compute $r_\nu$ $= -r_a +\begin{bmatrix}
C & 0
\end{bmatrix} F^{-1}\begin{pmatrix}
r_g\\ 
r_e
\end{pmatrix}$ by doing a forward/backward solve with the %\sout{$LDL^T$} \samahnew{
$L$ and $D$ factors of $F$ on $[r_g,r_e]^T$. Combining the cost of these triangular solves with that of the matrix-vector multiplication with $C$ yields a cost of
\[
N_I \times[ 2 c_\text{trsv}(L_{F}) + c_\text{spmv}(C)].
\]
After the linear system solve for these last two variants, we need to compute $[\Delta x, \Delta \lambda]^T$ using~\eqref{schur-delta-x}. This would involve a matrix vector multiplication with $C^T$ in order to get the right-hand side of~\eqref{schur-delta-x}, and a forward/backward substitution with the factors of $F$ at a cost of
\[
N_I \times[2 c_\text{trsv}(L_{F}) + c_\text{spmv}(C)].
\]
Note that $K_F$ is not explicitly formed, since only matrix-vector products with $K_F$ are needed to implement conjugate gradient. 

\subsubsection{Low-degree-of-freedom preconditioner (PL-KF)}

For this variant, the inequality-constraint reduced system~\eqref{ineq-const-reduced} is solved using preconditioned CG with a  preconditioner $P_L = D$. The preconditioner is diagonal, so the cost of applying it at every CG iteration is negligible. Thus, the main cost at every iteration comes from computing the matrix vector product $K_F \Delta \nu$. Since $K_F$ is never formed explicitly, this matrix-vector product is computed implicitly through three steps: a sparse matrix vector product to get $C^T \Delta \nu$, a forward/backward substitution with the $LDL^T$ factors of $F$, and a sparse matrix-vector product with $C$.
The total cost of these three steps is  $N_I[n_\text{kr} (2 c_\text{trsv}(L_F) + 2 c_\text{spmv}(C))  ]$. Thus, the overall cost of this variant is
\begin{equation}
\begin{aligned}
    %c_{v5} =
    & c_\text{fact}(F) +N_I[4c_\text{trsv}(L_{F}) + 2 c_\text{spmv}(C) + 2 n_\text{kr} c_\text{trsv}(L_F) + 2 n_\text{kr}c_\text{spmv}(C)  ].%\\
    %=& \underbrace{c_\text{fact}(F)}_\text{fact} + \underbrace{N_I\big[4c_\text{trsv}(L_F)+  2n_\text{kr}c_\text{trsv}(L_F) \big]}_\text{trsv} + \underbrace{N_I\big[2 c_\text{spmv}(C)+2 n_\text{kr}c_\text{spmv}(C)\big]}_\text{spmv}\\
    %=& \underbrace{c_\text{fact}(F)}_\text{fact} + \underbrace{2N_I(2+  n_\text{kr} )c_\text{trsv}(L_F)}_\text{trsv} + \underbrace{2N_I(1 + n_\text{kr})c_\text{spmv}(C)}_\text{spmv}.\\
    \end{aligned}
\end{equation}

\subsubsection{High-degree-of-freedom preconditioner (PH-KF)}
%\subsubsection{Variant 6, PH-KF}

For this variant, the inequality-constraint reduced system~\eqref{ineq-const-reduced} is solved using preconditioned CG with a preconditioner $P_H = D +C \diag(H)^{-1}C^T$. In order to compute this preconditioner, we can compute $T = C \diag(H)^{-1}C^T$ in the IPM setup phase at a cost of $c_\text{spmm}(C,D_HC^T)$, where $D_H = \diag(H)^{-1}$.
Subsequently, at every IPM iteration we can compute $P_H = D + T$, and then perform a Cholesky factorization $P_H=LL^T$ at a cost of $N_I \times c_\text{fact}(P_H)$. At every PCG iteration, one needs to perform a matrix vector product with $K_F$.
As we described in variant PL-KF above, this matrix vector product has a cost of $N_I[n_\text{kr} (2 c_\text{trsv}(L_F) + 2 c_\text{spmv}(C))  ]$.
Further, each iteration of PCG involves a forward/backward solve with the factors of $P_H$ at a cost of $N_I[n_\text{kr}(2 c_\text{trsv}(L_{P_H}))]$. Thus, the total cost is
\begin{equation}
\begin{aligned}
    %c_{v6} =
     c_\text{fact}(F) + c_\text{spmm}(C,D_HC^T) +N_I[&4c_\text{trsv}(L_{F}) + 2 c_\text{spmv}(C) +2 n_\text{kr} c_\text{trsv}(L_F)   \\
    &+ 2 n_\text{kr}c_\text{spmv}(C)+ c_\text{fact}(P_H) + 2 n_\text{kr} c_\text{trsv}(L_{P_H})   ].
    %=& \underbrace{c_\text{fact}(F)+N_I(c_\text{fact}(P_H))}_\text{fact} + \underbrace{N_I\big[4c_\text{trsv}(L_F)+  2n_\text{kr}c_\text{trsv}(L_F)+ 2 n_\text{kr} c_\text{trsv}(L_{P_H}) \big]}_\text{trsv}  \\
    %& + \underbrace{N_I\big[2 c_\text{spmv}(C)+2 n_\text{kr}c_\text{spmv}(C)\big]}_\text{spmv}+ \underbrace{c_\text{spmm}(C,DC^T)}_\text{spmm}\\
    %=& \underbrace{c_\text{fact}(F)+N_I(c_\text{fact}(P_H))}_\text{fact} + \underbrace{2N_I(2+  n_\text{kr} )c_\text{trsv}(L_F) + 2N_I n_\text{kr} c_\text{trsv}(L_{P_H}) }_\text{trsv} \\
    %&+ \underbrace{2N_I(1 + n_\text{kr})c_\text{spmv}(C)}_\text{spmv} + \underbrace{c_\text{spmm}(C,DC^T)}_\text{spmm}. 
    \end{aligned}
\end{equation}
\begin{table}[]
\centering
\footnotesize{
\begin{tabular}{|l|l|l|l|l|}
\hline
Variant                                                                        & Factorization & TRSV & SpMV & SpMM  \\ \hline
D-KC
  &  \begin{tabular}[c]{@{}l@{}}$N_I\times$\\ $c_\text{fact}(K_C)$\end{tabular}    &  $2N_I c_\text{trsv}(L_{K_C})$    &  $2N_I c_\text{spmv}(C)$    &  \begin{tabular}[c]{@{}l@{}} $N_I\times$\\$c_\text{spmm}(C^T{,}DC)$\end{tabular}    \\ \hline

U-KC &  $0$    &  $0$    &    \begin{tabular}[c]{@{}l@{}}$2N_I\times$\\$[n_\text{kr}c_\text{spmv}(K_C)$\\$+c_\text{spmv}(C)]$\end{tabular} &    \begin{tabular}[c]{@{}l@{}} $N_I\times$\\$c_\text{spmm}(C^T{,}DC)$\end{tabular}  \\ \hline

   CP-KC                                                                            & \begin{tabular}[c]{@{}l@{}}$N_I\times$\\ $c_\text{fact} (P_{CP})$\end{tabular}     &    \begin{tabular}[c]{@{}l@{}}$4N_I n_\text{kr}\times$\\$c_\text{trsv}(L_{P_{CP}})$\end{tabular}  &  \begin{tabular}[c]{@{}l@{}}$2N_I\times$\\$[n_\text{kr}c_\text{spmv}(K_C)$\\$+c_\text{spmv}(C)]$\end{tabular}    &  \begin{tabular}[c]{@{}l@{}} $N_I\times$\\$c_\text{spmm}(C^T{,}DC)$\end{tabular}   \\ \hline
   
 RG-KC                       &  \begin{tabular}[c]{@{}l@{}}$N_I\times$\\$ c_\text{fact} (P_{RG})$\end{tabular}     &    \begin{tabular}[c]{@{}l@{}}$4N_I n_\text{kr}\times$\\$c_\text{trsv}(L_{P_{RG}})$\end{tabular}  &  \begin{tabular}[c]{@{}l@{}}$2N_I\times$\\$[n_\text{kr}c_\text{spmv}(K_C)$\\$+c_\text{spmv}(C)]$\end{tabular}    &  \begin{tabular}[c]{@{}l@{}} $N_I\times$\\$[c_\text{spmm}(C^T{,}DC)$\\$+c_\text{spmm}(A^T{,}D_WA)]$\end{tabular} 
 %& 
 %\begin{tabular}[c]{@{}l@{}}$N_I\times$\\$ c_\text{norm}(G)$\end{tabular}
 \\ \hline
 
PL-KF                                                                               &  $c_\text{fact}(F)$    & \begin{tabular}[c]{@{}l@{}}$2N_I(2{+}n_\text{kr})\times$\\$c_\text{trsv}(L_{F})$\end{tabular}      &  \begin{tabular}[c]{@{}l@{}}$2N_I(1{+}n_\text{kr})\times$\\$c_\text{spmv}(C)$\end{tabular}    &  $0$     \\ \hline

     PH-KF                                                                          & \begin{tabular}[c]{@{}l@{}}$c_\text{fact}(F)$\\$+ N_I\times$\\$ c_\text{fact}(P_H)$\end{tabular}
     
     %$c_\text{fact}(F)+N_I c_\text{fact}(P_H)$
     &\begin{tabular}[c]{@{}l@{}} $2N_I\times$\\$\big[(2{+}n_\text{kr})c_\text{trsv}(L_{F})$ \\  $+n_\text{kr}c_\text{trsv}(L_{P_H})\big]$ \end{tabular} &  \begin{tabular}[c]{@{}l@{}}$2N_I(1{+}n_\text{kr})\times$\\$c_\text{spmv}(C)$\end{tabular}   &   $c_\text{spmm}(C{,}D_HC^T)$     \\ \hline
\end{tabular}
}
\caption{
%\edgarnew{
Cost components of each IPM variant (direct, unpreconditioned, and two preconditioned variants for each of the two Schur complements, $K_C$ and $K_F$).
%}
The cost of AL-KC additionally includes the term $N_I\times c_\text{norm}(G)$.
}
\label{tab:costs}
\end{table}
%\pagebreak

\subsubsection{Comparison}
Based on Table~\ref{tab:costs}, our preconditioned variants require fewer sparse matrix factorizations all alternatives, with the exception of U-KC.
The direct method factorizes $K_C$ as opposed to $F$, which is more expensive, unless $C$ has  a simple structure.
On the other hand, the new approaches require solving many systems with the factorization of $F$, as seen in the second column.
The sparse-matrix by vector multiplication cost is evenly matched among the variants, while the sparse-matrix by sparse-matrix multiplication cost decreases substantially when using one of our two variants PL-KF or PH-KF.
The sparse matrix-matrix products needed for PL-KF and PH-KF are easy to compute when $C$ has a simple form, e.g., when the inequality constraints described by $C$ are simple bounds.
Overall, the best choice of solver and preconditioner depends nontrivially on the sparsity structure of $H$, $A$, and $C$, as well as the desired accuracy.\label{sec:perf-analysis}

\section{Numerical experiments}\label{sec:numerics}

In this section, we report experimental results that reinforce the analysis presented in the previous sections. We implement all the IPM variants previously introduced in~\Cref{sec:perf-analysis}, in addition to an unpreconditioned CG variant for the inequality-constraint reduced system, which we call U-KF. {We collect results on the condition number of coefficient matrices, as well as number of iterations of Krylov subspace methods and number of IPM iterations. Using these iteration counts, we evaluate the costs given by the models developed in~\cref{sec:perf-analysis}.}
The initial guess for all iterative methods is taken to be the zero vector and the termination condition considers whether the residual norm is less than some tolerance $\epsilon$ times the norm of the right-hand side.
%\samahnew{with a maximum iteration limit of $10^5$}. 
%\samahnew{As one of the goals of inexact IPM methods is the ability to approximately solve for the Newton search direction.
We found the value of $\epsilon=10^{-3}$ to work well, so as to converge to the optimal solution in a reasonable number of IPM iterations. 
With this tolerance, the loss of accuracy in attained IPM objective value with respect to a direct solver is below $6\cdot 10^{-7}$, although IPM takes somewhat more iterations ($22\%$ as a median). In our comparative benchmarks, to ensure the IPM iteration count and linear systems are consistent among the iterative methods used, we replace the computed inexact solution with the result from direct solve at the end of each IPM step.
%median relative difference between ipm iterative (with tol$=1e-3$) and ipm direct solve $= 6.12e-07$
%median increase in number of ipm iterations $= 22\%$
We compare the final objective value given by our IPM solver to that obtained by the CVXOPT solver~\cite{vandenberghe2010cvxopt}. 
%and find that the relative difference between the two is less than $10^{-7}$ for all the test problems.  \samah{relative error(ipm - benchmark) median = 1.9e-2, geomean=2.0e-2.relative error(cvx - benchmark) median = 2.0e-2, geomean=2.1e-2.}

%\samah{Upon further inspection, I separated regularized from unregularized problems. For unregularized problems, relative error(ipm - benchmark) median = 4.3e-2, geomean=5.2e-2.relative error(cvxopt - benchmark) median = 4.3e-2, geomean=4.5e-2. While regularized only had:relative error(ipm - benchmark) median = 2.75e-2, geomean=3.0e-2.relative error(cvx - benchmark) median = 3.56e-2, geomean=3.23e-2.Compare with CVXOPT:relative error(ipm - cvx)(unregularized) median = 2.14e-7, geomean=1.92e-7.relative error(ipm - cvx)(both regularized) median = 4.31e-8, geomean=4.69e-8.}

All the numerical experiments are carried out using Python version 3.7, SciPy library version 1.3.2. 
For the CG and BiCGSTAB methods, we employ corresponding routines in the \texttt{scipy.sparse.linalg} submodule, and use the $LDL^T$ factorization from the \texttt{scipy.linalg} submodule, which uses LAPACK's \texttt{SYTRF} routines for symmetric matrices. 
%on a Dual-Core Intel Core i5 2.7 GHz processor with 16 GB of RAM. 
In the following subsections, we present more details on the test problems, numerical experiments, and results. 

\subsection{Test problems}

We consider two main types of test problems, synthetic QP test problems and QP problems from the Maros-M\'esz\'aros QP (MMQP) benchmark test set~\cite{maros1999repository}. These problems are instances of~\eqref{qp}.

{\it Synthetic QP} (\text{SyQP}) are convex quadratic programming problems of the form,
\begin{equation}
\label{synthetic-qp}
	\begin{aligned}
\underset{x \in \mathbb{R}^n}{\text{minimize}} \quad & \frac{1}{2}x^{T}Hx+x^{T}c\\
\textrm{subject to} \quad & Ax=b,\\
  &x\geq 0.    \\
\end{aligned}
\end{equation}
%Where 
We generate $H \in \bbR^{n\times n}$ to be a symmetric positive-definite matrix of block diagonal sparsity structure. 
Each block is computed from $M^TM$, where $M$ is a square matrix with entries selected uniformly at random from $[0.0,1.0]$.
We select this block size to be 4 in all generated synthetic matrices.
%Entries in the individual blocks, $H_{ii}$, are first drawn at random from the continuous uniform distribution over the interval $[0.0, 1.0)$, then set to $H_{ii}=H_{ii}^TH_{ii}$ in order to make it SPD.
$A \in \bbR^{m_1\times n}$ is generated to be full rank and banded, where the entries in the band are populated with random samples from a uniform distribution over $[0.0, 1.0)$.
$b \in \mathbb{R}^{m_1}$ and $c \in \mathbb{R}^{n}$ are similarly drawn at random.

%\sout{We randomly}
%\edgar{how specifically? what matrix is defined with element from what distribution?}
%\sout{generate $H \in \bbR^{n\times n}$ to be a symmetric positive definite matrix of block diagonal sparsity structure, $A \in \bbR^{m_1\times n}$ to be full rank and banded, $b \in \mathbb{R}^{m_1}$, $c \in \mathbb{R}^{n}$.}  

We denote with $\text{SyQP}^n_{m_1}$ the synthetic QP problem of the form~\eqref{synthetic-qp}  with $n$ primal variables, $n$ simple bound inequality constraints ($C=I$, $d=0$), and $m_1$ equality constraints, where $m_1\leq n$.
In our experiments, for a given $n$ we have generated a collection of quadratic programming problems that have the same number of primal variables but different number of equality constraints $m_1$. Specifically, we generate a set of $n$ quadratic programming problems
%\edgar{Is this $n$ the same as the dimension of $H$? (I guess you have in mind $m_1$, should also clarify that we have an $A$ with a different number of rows, i.e., that we vary the number of equality constraints.)} \samah{Yes this $n$ is the same as the dimension of $H$ because $m_1 = 1, \ldots, n$ so there are $n$ QP problems in total. I thought the next sentence makes it clear since $\text{SyQP}^n_{m_1}$ was defined as QP problem with $n$ primal variables and $m_1$ equality constraints, so each of these ones have different $A$s. Maybe ($\text{SyQP}^n_{m_1}$) is confusing, an alternative might be $\text{SyQP}(n,m_1)$. Not sure if that's better. }
%Edgar: I must have misread this previously, reads ok to me after your explanation
$\{\text{SyQP}^n_{1}, \text{SyQP}^n_{2},\ldots, \text{SyQP}^n_{n}\}$ which share the same Hessian matrix $H$, but have different equality constraint Jacobian matrices $A$.

\begin{figure}[htbp]
\centerline{\includegraphics[width=\linewidth,height=3.5in]%[width=5in, height=3in]
{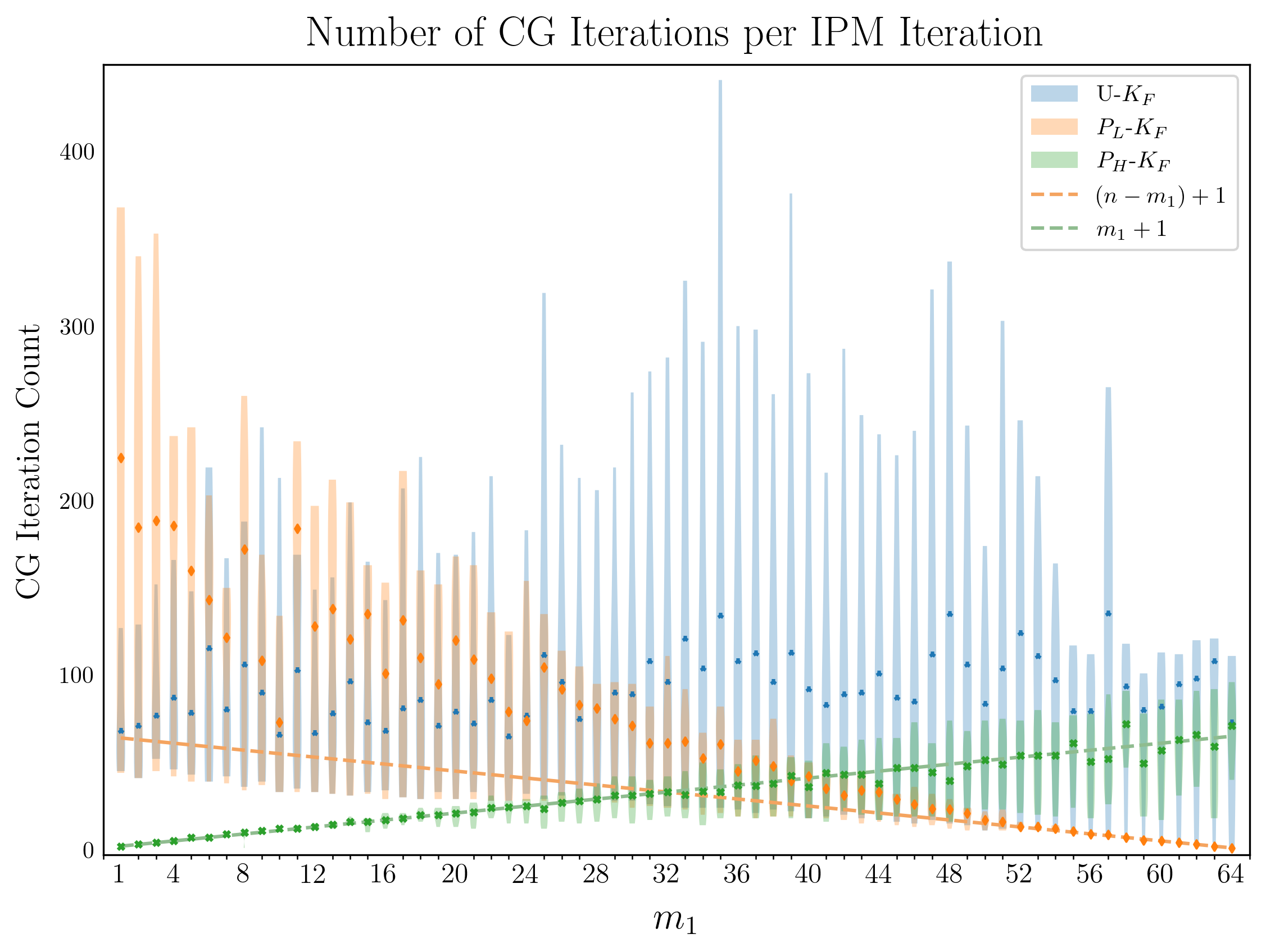}}
\caption{The range of conjugate-gradient iteration counts 
per IPM iteration,
for QPs with $n = 64$, and $m_1$  equality constraints ranging from $1$ to $n$. The CG solves were performed using preconditioners $P_L$, $P_H$ as well as without any preconditioning. The markers represent the median CG iteration counts for U-KF, PL-KF and PH-KF. The dashed lines $(\min(n,2(n-m_1))+1)$ and $(m_1+1)$ denote the theoretical upper bounds on CG iteration count using preconditioners $P_L$ and $P_H$ respectively.}
\label{violin-cg-64}
\end{figure}

We denote with MMQP the
%{\it Maros-M\'esz\'aros QP} (MMQP) are 
convex quadratic programming problems~\eqref{qp} from the {\it Maros-M\'esz\'aros} benchmark test set~\cite{maros1999repository}, which is often used to benchmark convex QP solvers.  These QP problems stem from both academic and real world applications, and can be accessed through the CUTEst testing environment~\cite{gould2015cutest}.
In our tests, we use QP problems from the MMQP test set  for which the matrix $A$ has full rank. \Cref{tab:cutest-nnz} provides a list of the MMQP problems used in our numerical experiments. %along with some of  as well as some of their important characterizing features.}
%\sout{We separate the problems used into sublists based on some of properties of the QP, one being the positive definiteness of the given Hessian matrix $H$, so that sublist $1$ has \texttt{DUAL1}, \texttt{DUAL2}, \texttt{DUAL3}, \texttt{DUAL4}, \texttt{DUALC1}, \texttt{DUALC5}, \texttt{QPCBLEND}, \texttt{QPCBOEI1}, \texttt{QPCBOEI2}, \texttt{QPCSTAIR}.}
%\edgar{We don't seem to separate them anywhere anymore.}

Most of the problems in this benchmark set have singular %\samahnew{
(positive semi-definite) Hessian matrices, contrary to the assumption we made in Section~\ref{sec:intro}.
%which can be problematic for both direct and iterative solvers of the augmented system~\eqref{augmented-2}\sout{, and does not conform with our single-factorization inexact IPM formulation}.
%, and that is why some form of regularization is often used. %One possible method for dealing with an ill-conditioned or singular $H$ is the augmented Lagrangian method, see~\cite{fortin1983augmented},~\cite{glowinski1989augmented},~\cite{golub2003solving},~\cite{golub2005algebraic} for a more detailed analysis. 
%\sout{We employ the primal regularization technique~\cite{altman1999regularized} for dealing with ill-conditioned or singular $H$.In particular, we add a multiple of the identity to the Hessian block in the first augmented system $H + \rho I$, with $\rho = 10^{-4}||H||_2$. }
{
%One approach to alleviate those issues is by introducing regularization terms to the optimization problem. First proposed in~\cite{gill1991solving,saunders1996cholesky}, regularization techniques have been extensively used in optimization to improve the numerical properties of KKT linear systems. The reader is referred to~\cite{gondzio2012matrix,altman1999regularized,pougkakiotis2019dynamic,pougkakiotis2021interior} for the benefits of using regularization in the context of IPM. One such approach is the primal regularization technique~\cite{altman1999regularized} which
Consequently, we make use of the primal regularization technique~\cite{altman1999regularized} when working with the MMQP benchmark problems. Using a primal regularization parameter $\rho = 10^{-4}||H||_2>0$, the Hessian in the (1,1) block of~\eqref{augmented} is modified into $H+\rho I$, thus rendering it positive definite.}
%With the right-hand-side unchanged, the new linear system can be interpreted as the KKT system for the {\it primal-regularized} QP~\cite{friedlander2012primal},
%
%\begin{equation}
%    \label{regular-qp}
%    \begin{aligned}
%\underset{x \in \mathbb{R}^n}{\text{minimize}} \quad & \frac{1}{2}x^{T}Hx+\frac{1}{2}\rho||x-x_k||^2+x^{T}c\\
%\textrm{subject to} \quad & Ax=b,\\
%  &Cx\geq d.    \\
%\end{aligned}
%\end{equation}
%where $x_k$ is the current iterate.
The regularization introduces a median deviation of $2.9\cdot 10^{-2}$ with respect to the objective value of the solution given by the CVXOPT solver with no regularization. While for problems with a symmetric positive-definite Hessian, the objective solution given by our IPM solver achieves a median relative accuracy of $2.9\cdot 10^{-7}$ compared to the reference solution.

\begin{table}[]%[htbp]
\centering
%\footnotesize{
\begin{tabular}{ |l|l|l|l|l|l|l|l| }
%{ |p{2.1cm}|p{0.8cm}|p{0.8cm}|p{0.8cm}|p{1.0cm}|p{1.0cm}|p{1.0cm}|p{1.11cm}| }
 \hline
 Name& $n$ & $m_1$ & $m_2$ & $\nz(A)$ & $\nz(C)$ & $\nz(H)$& $\nz(H_L)$  \\
  \hline
DUAL1 & 85 & 1 & 170 & 85 & 170 & 7031 & 3473     \\
\hline
DUAL2 & 96 & 1 & 192 & 96 & 192 & 8920 & 4412     \\
\hline
DUAL3 & 111 & 1 & 222 & 111 & 222 & 12105 & 5997     \\
\hline
DUAL4 & 75 & 1 & 150 & 75 & 150 & 5523 & 2724     \\
\hline
DUALC1 & 9 & 1 & 232 & 9 & 1944 & 81 & 36     \\
\hline
DUALC5 & 8 & 1 & 293 & 8 & 2232 & 64 & 28     \\
\hline
QPCBLEND & 83 & 43 & 114 & 298 & 276 & 83 & 0     \\
\hline
QPCBOEI1 & 384 & 9 & 971 & 168 & 4191 & 384 & 0     \\
\hline
QPCBOEI2 & 143 & 4 & 378 & 56 & 1424 & 143 & 0     \\
\hline
QPCSTAIR & 467 & 209 & 696 & 1374 & 3031 & 467 & 0     \\
\hline
QBANDM & 472 & 305 & 472 & 2494 & 472 & 504 & 16     \\
\hline
QADLITTL & 97 & 15 & 138 & 173 & 307 & 237 & 70     \\
\hline
QAFIRO & 32 & 8 & 51 & 34 & 81 & 38 & 3     \\
\hline
QBEACONF & 262 & 140 & 295 & 3309 & 328 & 280 & 9     \\
\hline
QE226 & 282 & 33 & 472 & 938 & 1922 & 2076 & 897     \\
\hline
QFFFFF80 & 854 & 350 & 1028 & 4775 & 2306 & 4130 & 1638     \\
\hline
QSC205 & 203 & 91 & 317 & 249 & 505 & 223 & 10     \\
\hline
QSCAGR25 & 500 & 300 & 671 & 1334 & 720 & 700 & 100     \\
\hline
QSCAGR7 & 140 & 84 & 185 & 362 & 198 & 174 & 17     \\
\hline
QSCFXM1 & 457 & 187 & 600 & 1467 & 1579 & 1811 & 677     \\
\hline
QSCFXM2 & 914 & 374 & 1200 & 2939 & 3158 & 3028 & 1057     \\
\hline
QSCTAP1 & 480 & 120 & 660 & 360 & 1812 & 714 & 117     \\
\hline
QSHARE1B & 225 & 89 & 253 & 891 & 485 & 267 & 21     \\
\hline
QSHARE2B & 79 & 13 & 162 & 84 & 689 & 169 & 45     \\
\hline
QSCRS8 & 1169 & 384 & 1275 & 2576 & 1775 & 1345 & 88     \\
\hline
CVXQP1\_M & 1000 & 500 & 2000 & 1498 & 2000 & 6968 & 2984     \\
\hline
CVXQP3\_M & 1000 & 750 & 2000 & 2247 & 2000 & 6968 & 2984     \\
\hline
CVXQP1\_S & 100 & 50 & 200 & 148 & 200 & 672 & 286     \\
\hline
CVXQP3\_S & 100 & 75 & 200 & 222 & 200 & 672 & 286     \\
\hline
QGROW15 & 645 & 300 & 1245 & 5620 & 1245 & 1569 & 462     \\
\hline
QGROW22 & 946 & 440 & 1826 & 8252 & 1826 & 2520 & 787     \\
\hline
QGROW7 & 301 & 140 & 581 & 2612 & 581 & 955 & 327     \\
\hline
VALUES & 202 & 1 & 404 & 202 & 404 & 7442 & 3620     \\
\hline
DUALC2 & 7 & 1 & 242 & 7 & 1610 & 49 & 21     \\
\hline
DUALC8 & 8 & 1 & 518 & 8 & 4032 & 64 & 28     \\
\hline
QETAMACR & 688 & 272 & 1033 & 1374 & 1940 & 8826 & 4069     \\
\hline
QFORPLAN & 421 & 90 & 517 & 3775 & 1234 & 1513 & 546     \\
\hline
QRECIPE & 180 & 67 & 299 & 351 & 587 & 240 & 30     \\
\hline
QSTAIR & 467 & 209 & 696 & 1374 & 3031 & 2371 & 952     \\
\hline
QSTANDAT & 1075 & 160 & 1394 & 2128 & 2098 & 2407 & 666     \\
\hline
QSEBA & 1028 & 507 & 1550 & 4330 & 1572 & 2128 & 550     \\
\hline
\end{tabular}
%}
\caption{List of the Maros-M\'esz\'aros QP test problems. $n$ is the number of variables, $m_1$ is the number of equality constraints, $m_2$ is the number of inequality constraints,
$\nz(X)$ is the number of nonzeros in $X$, and $H_L$ is strictly lower triangular part of $H$.
}\label{tab:cutest-nnz}
\end{table}

\subsection{Conjugate gradient iteration count}
In this numerical experiment, we use synthetic quadratic programming problems $\text{SyQP}^{64}_1, \text{SyQP}^{64}_2,\ldots, \text{SyQP}^{64}_{64}$ with $n=64$ primal variables, $n=64$ simple bounds, and $m_1 = 1,\ldots, 64$ equality constraints. For each QP, we execute three variants of IPM, all solving the inequality-constraint reduced system~\eqref{ineq-const-reduced}.
%\sout{They employ} 
The conjugate gradient method is employed with no preconditioning (U-KF), with our low-d.o.f. preconditioner $P_L$~\eqref{prec-mat-1} (PL-KF), and with our high-d.o.f. preconditioner $P_H$~\eqref{prec-mat-2} (PH-KF).

\Cref{violin-cg-64} shows the distribution of CG iteration count per IPM iteration for every convex quadratic program $\text{SyQP}^{64}_{m_1}$, with $ m_1 = 1,\ldots,64$. The median for each of the distributions is shown (with markers) on the violin plot, as well as the theoretical upper bounds on CG iteration count (with dashed lines) based on~\cref{cor:cg1,cor:cg2}. This experiment illustrates the effectiveness of our preconditioners in reducing CG iteration count compared to the unpreconditioned variant U-KF. $P_L$ significantly reduces the median number of CG iterations (per IPM iteration) when the number of degrees-of-freedom ($n_d=n-m_1$) is less than or equal than $n/2$. While $P_H$ is not only successful in reducing the median Krylov iteration count in the high-d.o.f. regime ($m_1 \ll n$), but also for all the QPs tested (with $n_d$ as small as 0) . %where the number of equality constraints ranged from $1$ all the way to $n$.}

Further, we observe that the median values of the CG iteration count per IPM iteration generally conform to the theoretical %expected values that we \samahnew{have} prove\samahnew{d} 
bounds proven in~\cref{thm1,thm2}. Using our %low-degree-of-freedom 
low-d.o.f. preconditioner $P_L$, the median number of CG iterations needed per IPM iteration does not exceed $2(n-m_1)+1$ for $m_1>n/2$, reaching a minimum of $1$ when $m_1 = n$. 
Using our high-d.o.f. preconditioner $P_H$, the median number of CG iterations needed per IPM iteration is less than or equal to $m_1+1$ iterations, reaching a minimum of $2$ iterations when $m_1 = 1$. 

There are deviations from the theoretical expected CG iteration counts especially due to finite precision error for PL-KF. %, and that is due to finite precision error.
Rounding error in CG is significant in later IPM iterations, as its amplification is greatest for ill-conditioned problems.
The use of the %low-degree-of-freedom
low-d.o.f. preconditioner in a %high-degree-of-freedom 
high-d.o.f. scenario (when the number of equality constraints $m_1$ is small) can result in an increase of condition number, slowing down CG due to increased round-off error.
This hypothesis is supported by the results presented in the next subsection. %For it is a common feature of interior point methods that some diagonal entries in $D$ become very large and others very small as the method nears optimality~\cite{forsgren2007iterative}. 

%The inequality constraints can take the form of {\it simple bounds} on the variables $x$, for instance, non-negativity constraints $(x \geq 0)$, in which case $C$ will be equal to an identity matrix, or upper and lower bounds on the variables $l \leq x \leq u$, where $C$ will be equal to the concatenation of an identity matrix with a negative identity matrix. We will refer to all other linear inequality constraints as {\it general inequality} constraints.

\subsection{Conditioning}
\begin{figure}[htbp]
\centerline{\includegraphics[width=\linewidth,height=3.5in]{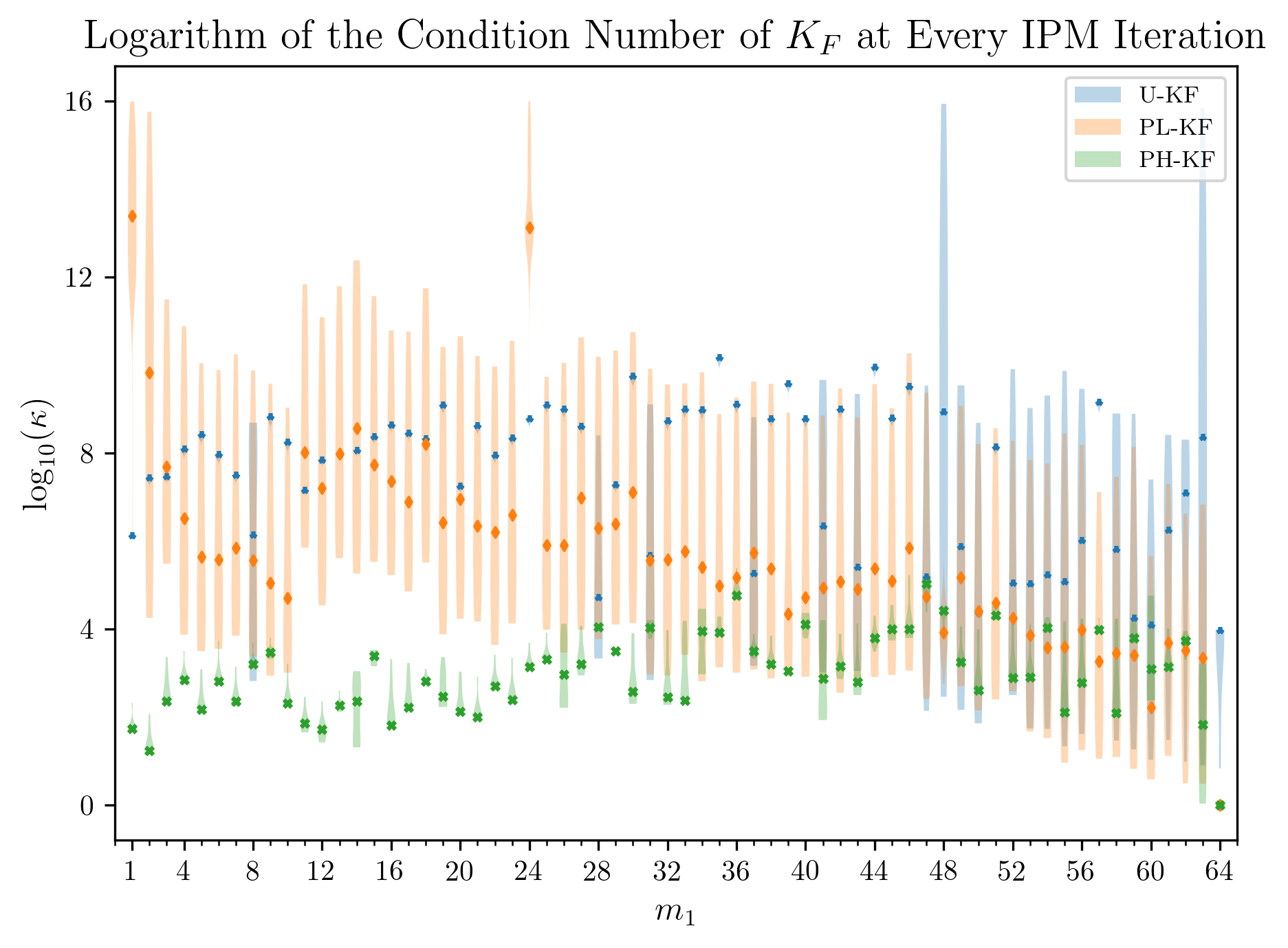}}
\caption{Condition number at every IPM iteration of the (preconditioned) inequality-constraint reduced matrix for  QPs with $n = 64$, and with $m_1$ equality constraints ranging from $1$ to $n$. The CG solves were performed using preconditioner $P_L$, $P_H$ as well as without any preconditioning.  The markers represent the median values of $\log_{10}(\kappa)$ for $K_F$, $P_L^{-1/2}K_FP_L^{-1/2}$ and $P_H^{-1/2}K_FP_H^{-1/2}$ .}
\label{violin-cond-64}
\end{figure}

\begin{figure}[htbp]
\centerline{\includegraphics[width=\linewidth]{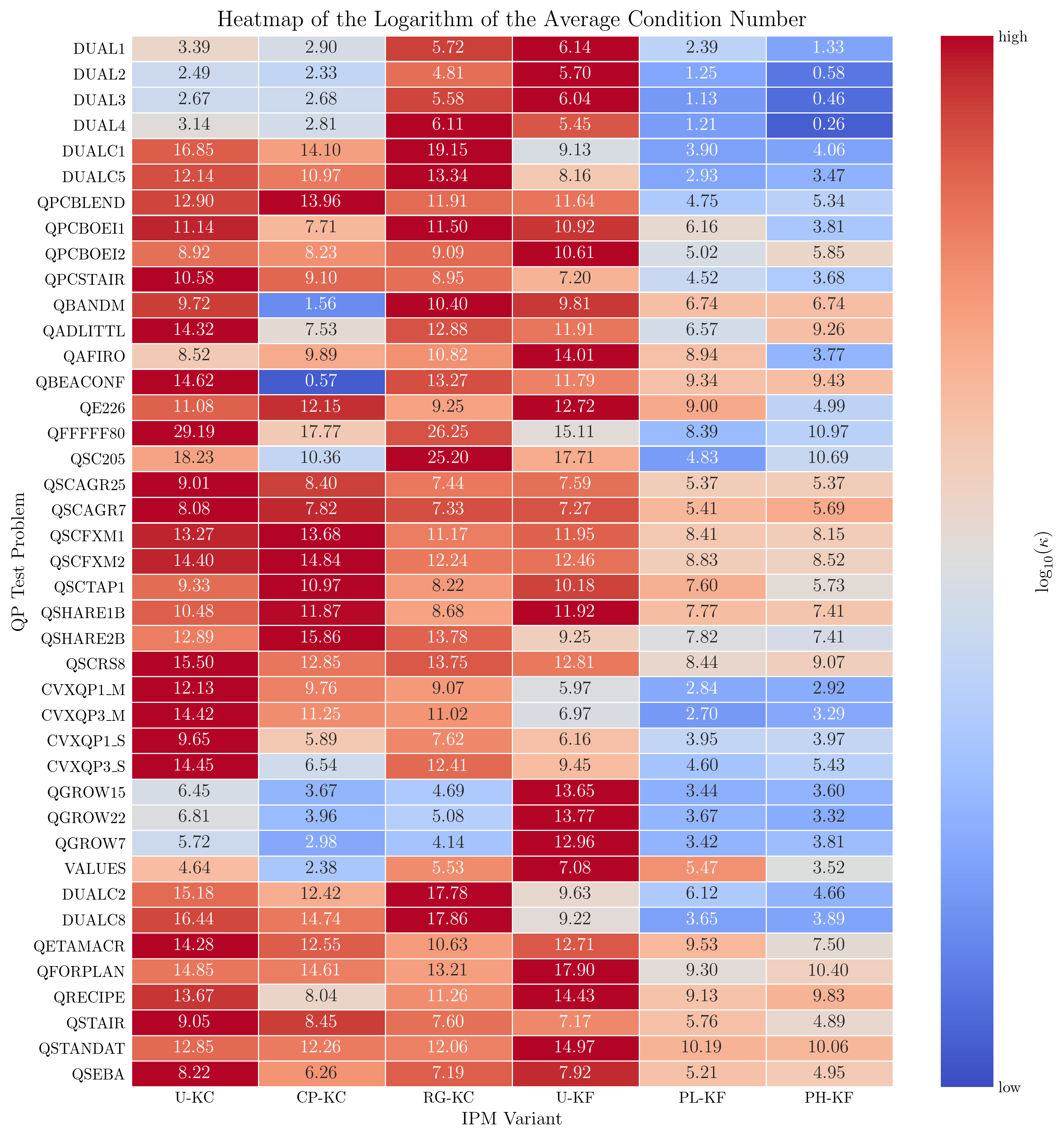}}
  %\centering
  %\includegraphics{figures/cutest-heatmap-all-cond.png}
  \caption{Heatmap of the logarithm of the average (geometric mean) condition number across IPM iterations for the augmented system both without preconditioning (U-KC) and preconditioned using a constraint preconditioner (CP-KC) and a block preconditioner (RG-KC), as well as the inequality-constraint reduced system both without preconditioning (U-KF) and preconditioned using the %low-degree-of-freedom 
  low-d.o.f. (PL-KF) and the %high-degree-of-freedom 
  high-d.o.f. (PH-KF) preconditioners.}
\label{cutest-cond}
\end{figure}

\Cref{violin-cond-64} considers the improvement to the conditioning of the inequality-constraint reduced system~\eqref{ineq-const-reduced} achieved by our low-d.o.f. preconditioner $P_L$ and high-d.o.f. preconditioner $P_H$. 
%in improving the spectral properties of the 
%We compare these improvements to the conditioning of the baseline .
The plot shows the distribution of the logarithm of the condition number across IPM iterations for $K_F$,  $P_L^{-1/2}K_F P_L^{-1/2}$, and $P_H^{-1/2}K_F P_H^{-1/2}$,  corresponding to the synthetic quadratic programming problems, $\text{SyQP}^{64}_1, \ldots, \text{SyQP}^{64}_{64}$.
%test problems were used with $n=64$ primal variables, $m_1 = 1, ...., 64$ equality constraints and $m_2 = n$ simple bounds on the variables.
$P_H$ improves the condition number of $K_F$ for all test problems with any value of $m_1$. 
Similarly,
$P_L$ improves the condition number of $K_F$ when $m_1 \geq n/2$ and in $77 \%$ of problems when $m_1 < n/2$. This reduction in condition number of $K_F$ provides the empirical confirmation of the improvements in CG convergence properties presented in~\Cref{violin-cg-64}.  %A preconditioner is necessary in improving the convergence properties of the iterative method by reducing the condition number and/or achieving good clustering of the eigenvalues of the preconditioned matrix~\cite{demmel1997applied},~\cite{golub2013matrix}. 

%Figure shows 
%describes how
%Our bounds show how cond shoud be good, 
%that our precond 

%back reference table 1 in this section, Table~\ref{precs-eigs}

Next, using the QP problems from the Maros-M\'esz\'aros test set, we execute all inexact variants of the IPM.
These include U-KC, CP-KC, RG-KC for solving the augmented KKT system~\eqref{augmented-2} with BiCGSTAB both without preconditioning, and preconditioned with a constraint preconditioner $P_{CP}$~\eqref{const-prec}, and a block diagonal preconditioner $P_{RG}$~\eqref{rees-greif}, respectively.
We also consider the variants U-KF, PL-KF, and PH-KF for solving the inequality-constraint reduced system~\eqref{ineq-const-reduced} using CG both without preconditioning, and preconditioned using our preconditioners, $P_L$ and $P_H$.

\Cref{cutest-cond} quantifies the average (across IPM iterations) of the condition numbers corresponding to unpreconditioned and preconditioned matrices $K_C$, $P_{CP}^{-1}K_C$, $P_{RG}^{-1}K_C$, $K_F$, $P_L^{-1/2}K_F P_L^{-1/2}$, and $P_H^{-1/2}K_F P_H^{-1/2}$.
These constitute the coefficient matrices of the linear systems in IPM variants U-KC, CP-KC, RG-KC, U-KF, PL-KF, and PH-KF, respectively.
The average condition number $\kappa$ is computed as the geometric mean over interior point iterations
%\edgarnew{Specifically, we compute the geometric mean of the condition number across interior point iterations for each problem, then take the arithmetic mean across test problems.}
% The colorbar in the heatmap goes from blue to red with the increase in $\log_{10}(\kappa)$, where $\kappa$ denotes the condition number. 
, and $\log_{10}(\kappa)$ for each QP is displayed in the respective cell of the heatmap.
%\sout{The number displayed in each cell is $\log_{10}(\kappa)$, where $\kappa$ is the average condition number.}

\Cref{tab:cutest-cond-table} (included at the end of the paper) provides further detail on the relative comparison of our low-d.o.f. and high-d.o.f. preconditioners to other state-of-the art preconditioners.
For every QP, the table presents a ratio of improvement in average condition number $\kappa$ for the best for our variants PL-KF and PH-KF with respect to the best alternative among U-KC, CP-KC and RG-KC. To summarize the results, we compute the geometric mean of the ratios across all MMQP problems and find that our preconditioners $P_L$ and $P_H$ achieve on average lower condition numbers for the coefficient matrices 
by a factor of $240$ relative to the best alternative iterative variant, whether preconditioned (CP-KC, RG-KC) or otherwise (U-KC, U-KF).
%Our variants PL-KF and PH-KF achieve a lower condition number for the preconditioned coefficient matrices far more often than variants CP-KC and RG-KC, as well as the unpreconditioned matrices of variants U-KC and U-KF. 
These results indicate that our preconditioners are successful in achieving their goal of improving the spectral properties of the inequality-constraint reduced system, and compare favorably to other popular preconditioners.

%The preconditioned matrices $P_L^{-1/2}K_F P_L^{-1/2}$ and $P_H^{-1/2}K_F P_H^{-1/2}$ have lower median condition number ($\tilde{\kappa}(P_L^{-1/2}K_F P_L^{-1/2}) = 2.5e+5$, $\tilde{\kappa}(P_H^{-1/2}K_F P_H^{-1/2}) \approx 2.1e+5$) compared to that of the preconditioned augmented matrices $P_{CP}^{-1}K_C$ ($\tilde{\kappa}(P_{CP}^{-1}K_C) \approx 1.2e+9$) and $P_{RG}^{-1}K_C$ ($\tilde{\kappa}(P_{RG}^{-1}K_C) \approx 4.2e+10$). And in the best case, $P_H$ and $P_L$ are able to achieve low condition numbers of $1.8$ and $15.5$, respectively.
\begin{table}[]%[htbp]
\centering
\begin{tabular}{c|c|c}
Coefficient Matrix & $\log_{10}(\tilde{\kappa})$ &  $\log_{10}(\hat{\kappa})$\\ \hline
$K_C$ & 11.14 & 11.26\\ \hline
$P_{CP}^{-1}K_C$&  9.1 & 8.95 \\ \hline
$P_{RG}^{-1}K_C$& 10.63  &  10.83 \\ \hline
$K_F$ & 10.18 & 10.43\\ \hline
$P_L^{-1/2}K_F P_L^{-1/2}$&  5.41  &  5.75 \\ \hline
$P_H^{-1/2}K_F P_H^{-1/2}$&  4.99  & 5.56 \\ \hline
\end{tabular}
\caption{ Condition numbers of (preconditioned) $K_C$ and $K_F$ averaged over all MMQP test problems.
$\tilde{\kappa}$ is the median  over QP test problems of the average condition number $\kappa$ (geometric mean over IPM iterations) , $\hat{\kappa}$ is the geometric mean over QP test problems of the average condition number $\kappa$ (geometric mean over IPM iterations).}\label{tab:cond-summary}
\end{table}

\Cref{tab:cond-summary} quantifies average condition %summarizes the condition number results \sout{in} \samahnew{of} Figure~\ref{cutest-cond}.
of the coefficient matrices $K_C$, $P_{CP}^{-1}K_C$, $P_{RG}^{-1}K_C$, $K_F$, $P_L^{-1/2}K_F P_L^{-1/2}$, and $P_H^{-1/2}K_F P_H^{-1/2}$. This table presents both the median and geometric mean of $\kappa$ over QP test problems.  
%First, we compare the conditioning of our new reduced system~\eqref{delnu} to that of the augmented system~\eqref{augmented-2}.
%, and for that we look at the average condition numbers for the MMQP problems in the U-KF and U-KC columns.
Comparing the median condition number of our inequality-constraint reduced matrix to that of the augmented matrix $K_C$, we note that $K_F$ has a lower condition number in terms of both median value and geometric mean value.  %\sout{the median value is higher for $K_C$ $(\tilde{\kappa}(K_C) \approx 1.4e+11)$ than $K_F$ $(\tilde{\kappa}(K_F) = 1.5e+10)$}. 
This means, that unlike the normal equations, our Schur complementation does not negatively impact the conditioning of the KKT matrix. In fact, even in the most ill-conditioned case, the matrix $K_F$ still has a lower condition number ($\max(\log_{10}(\kappa(K_F))) = 17.9)$) than that of $K_C$ ($\max(\log_{10}(\kappa(K_C))) \approx 29.19)$).

%unpreconditioned matrix $K_F$ for the MMQP problems considered. 

\subsection{Cost comparison}
\begin{figure}[htbp]
\centerline{\includegraphics[width=\linewidth]{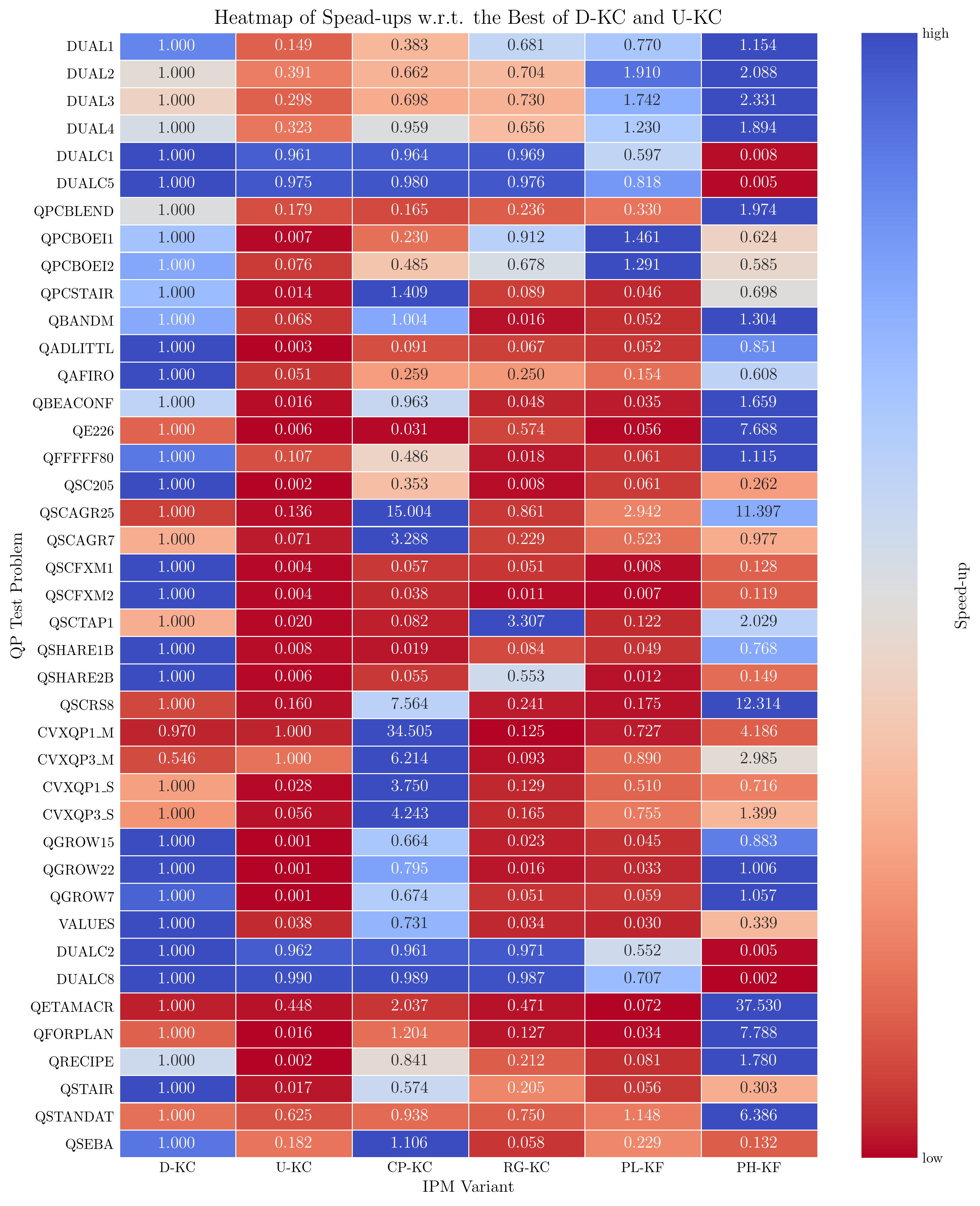}}
  %\centering
  %\includegraphics{figures/cutest-heatmap-all-cond.png}
  \caption{Heatmap of the speedup with respect to  the best of D-KC and U-KC for some QP test problems for the CP-KC, RG-KC, PL-KF, and PH-KF variants.}
\label{speedup-all-ldld-u}
\end{figure}

We now evaluate the arithmetic cost of each of the methods by applying our cost models to the specific sparse matrices arising in the QP test problems.
First, we compare the IPM variants with a preconditioned iterative solver, namely CP-KC, RG-KC, PL-KF, and PH-KF to those that either use a direct solver (D-KC) or an unpreconditioned iterative solver (U-KC).
Using the models we developed in~\cref{sec:perf-analysis}, we calculate FLOP counts for  the IPM variants formerly mentioned %\sout{D-KC, U-KC, CP-KC, RG-KC, PL-KF and PH-KF for} \samahnew{corresponding to each of} 
necessary for each of the MMQP test problems. 
We calculate speedup with respect to the best variant amongst D-KC and U-KC, according to the following formula
%\eqref{speedup-eq}

\begin{equation}
\label{speedup-eq}
 s_\text{i}^\text{D-U}(q_j)= \frac{\min\bigl(c_\text{D-KC}(q_j),c_\text{U-KC}(q_j)\bigr)}{c_{i}(q_j)},\;\;\;\;\;\left\{\begin{matrix}
i\in \{\text{CP-KC, RG-KC, PL-KF, PH-KF}\},\\ 
q_j \in \text{\{MMQP\}} .
\end{matrix}\right. 
\end{equation}
where $c_i(q_j)$ is the cost of method $i$ for QP problem $q_j$, and present the results in~\cref{speedup-all-ldld-u}.

%\sout{Comparing median speedups, $P_H$ is the only preconditioner able to achieve a speedup with respect to the direct solver variant ($\med(s_\mathrm{PH-KF})=1.01$).} %The 75th percentile is the value at which $25\%$ of the answers lie above that value, and $75\%$ lie above it. And for $75\%$ percentile of the test problems $P_H$ achieved the highest speedup of $3.31$, followed by $P_{CP}$ which achieved an speedup of $1.83$, and then $P_L$ which achieved a speedup of $1.27$, while $P_{RG}$ is still unsuccessful in achieving speedup ($s_\mathrm{RG-KC}=0.79$). 
%In terms of maximum speedup, 
Comparing maximum speedups, $P_H$ is able to achieve the highest speedup of $37.53$, followed by $P_{CP}$ with a maximum speedup of $34.5$,
$P_{RG}$ with a maximum speedup of $3.3$, and lastly
$P_L$ with a maximum speedup of $2.9$.

\begin{figure}[htbp]
\centering
  \includegraphics[width=\linewidth,height=4in]
  {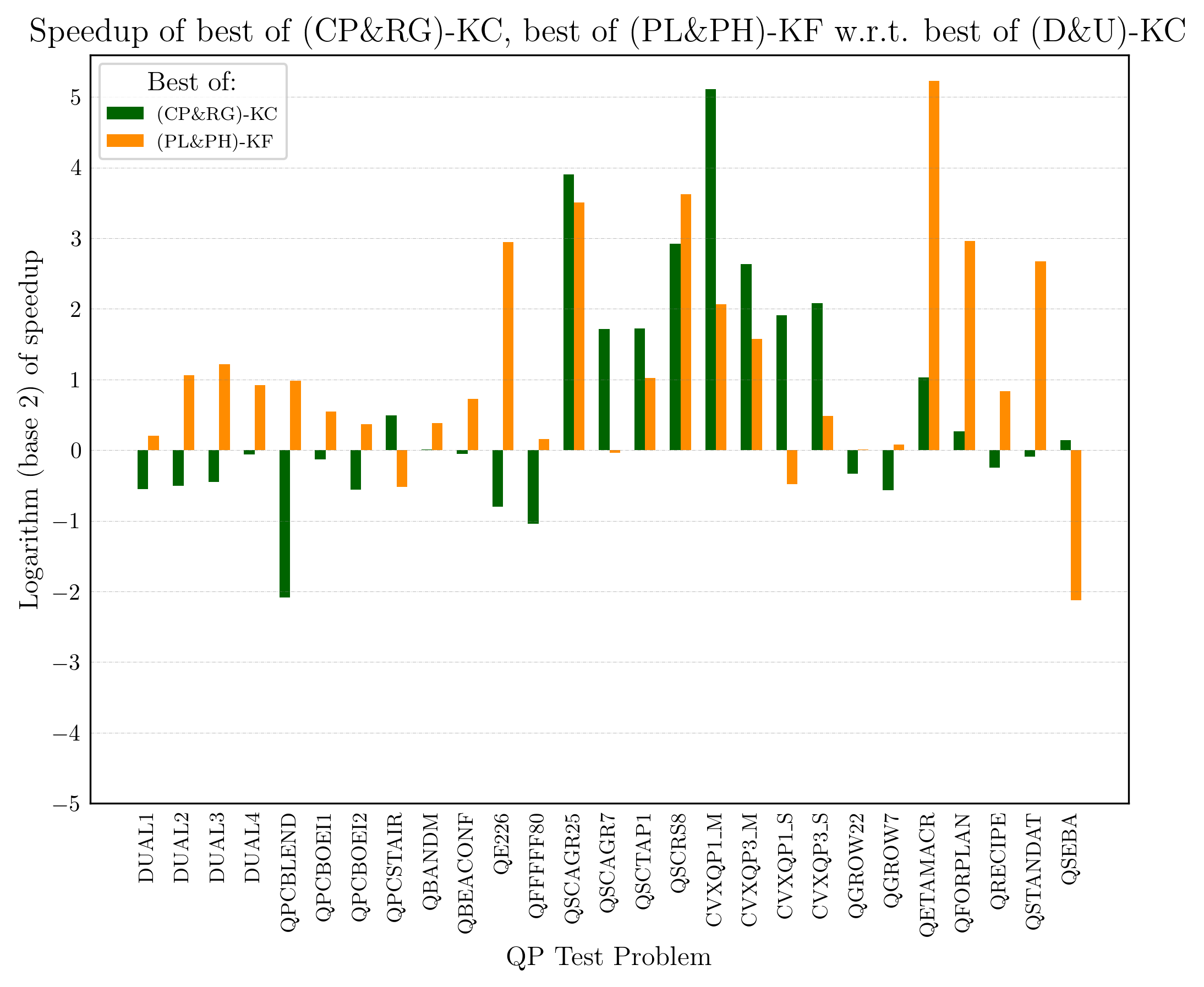}
  \caption{Speedups with respect to the fastest of the D-KC direct variant and the unpreconditioned iterative variant U-KC for the fastest of the augmented KKT system preconditioners $P_{CP}$ and $P_{RG}$, compared to the fastest of the inequality-constraint reduced system preconditioners $P_L$ and $P_H$ for QP problems from the MMQP test problems.%\samah{(CP\&RG)median=$0.96$ and geomean=$1.52$. (PL\&PH) median=$1.78$ and geomean=$2.18$}
  }
  \label{speedup-barplot_all}
\end{figure}

Excluding problems which are amenable to efficient direct solve (for which the direct solver D-KC is fastest among the variants) allows us to more clearly compare the performance among inexact IPM methods. We will denote this subset of problem as $\text{MMQP}_\text{XD}$ for the remainder of our experiments. \Cref{speedup-barplot_all} shows the maximum speedup among the (CP-KC, RG-KC) variants, $\max\bigl(s_{CP-KC}^{D-U},s_{RG-KC}^{D-U}\bigr)$, and the (PL-KF, PH-KF) variants, $\max\bigl(s_{PL-KF}^{D-U},s_{PH-KF}^{D-U}\bigr)$, for the QP test problems for which at least one of the four variants has a speedup $s_i^{D-U}$ greater than $1$.%\sout{ with respect to (D-KC, U-KC)}. 

Comparing speed-ups across the $\text{MMQP}_\text{XD}$ subset of problems, the new preconditioned variants (PL\&PH)-KF provide the lowest and most robust cost for the tested problems among alternative methods. Our preconditioned inexact IPM methods achieve a median speed-up of $1.8$ across this problem set in~\cref{speedup-barplot_all}, compared to $0.9$ for the (CP\&RG)-KC variants. \Cref{tab:speedup-summary} further shows that the PH-KF preconditioner is fastest for 59\% of these problems, and taken together, either PL-KF or PH-KF is the method of choice for 2/3 of the problems.
%\samah{(CP\&RG)median=$0.96$ and geomean=$1.52$. (PL\&PH) median=$1.78$ and geomean=$2.18$}

\begin{table}[]%[htbp]
\centering
\begin{tabular}{c|c|c|c}
Preconditioner & $p(\%)$ & $\tilde{s}^U$ & $\hat{s}^U$\\ \hline
PL-KF& 7.41 &  4.89 & 4.58 \\ \hline   %3.285 &  3.797(these on the right include all problems)\\ \hline
PH-KF  & 59.25 &  25.16 & 34.08 \\ \hline  %24.995 &  15.445  \\ \hline
CP-KC & 29.63 & 14.79 & 18.24 \\ \hline % 9.308 &  14.221  \\ \hline
RG-KC  & 3.7 &  2.97 &  3.53  \\ \hline
\end{tabular}
\caption{$p$ is the percentage of problems this preconditioner is the fastest among the four preconditioned variants for QP test problems (excluding problems where the direct method is the fastest). %$p_U$ is the percentage of problems faster than U-KC, 
$\tilde{s}^U$ is the median speedup with respect to U-KC and $\hat{s}^U$ is the geometric mean speedup with respect to U-KC, across test problems where the direct method D-KC is not the fastest.
}\label{tab:speedup-summary}
\end{table}

\begin{figure}[htbp]
\centering
  \includegraphics[width=\linewidth,height=4in]{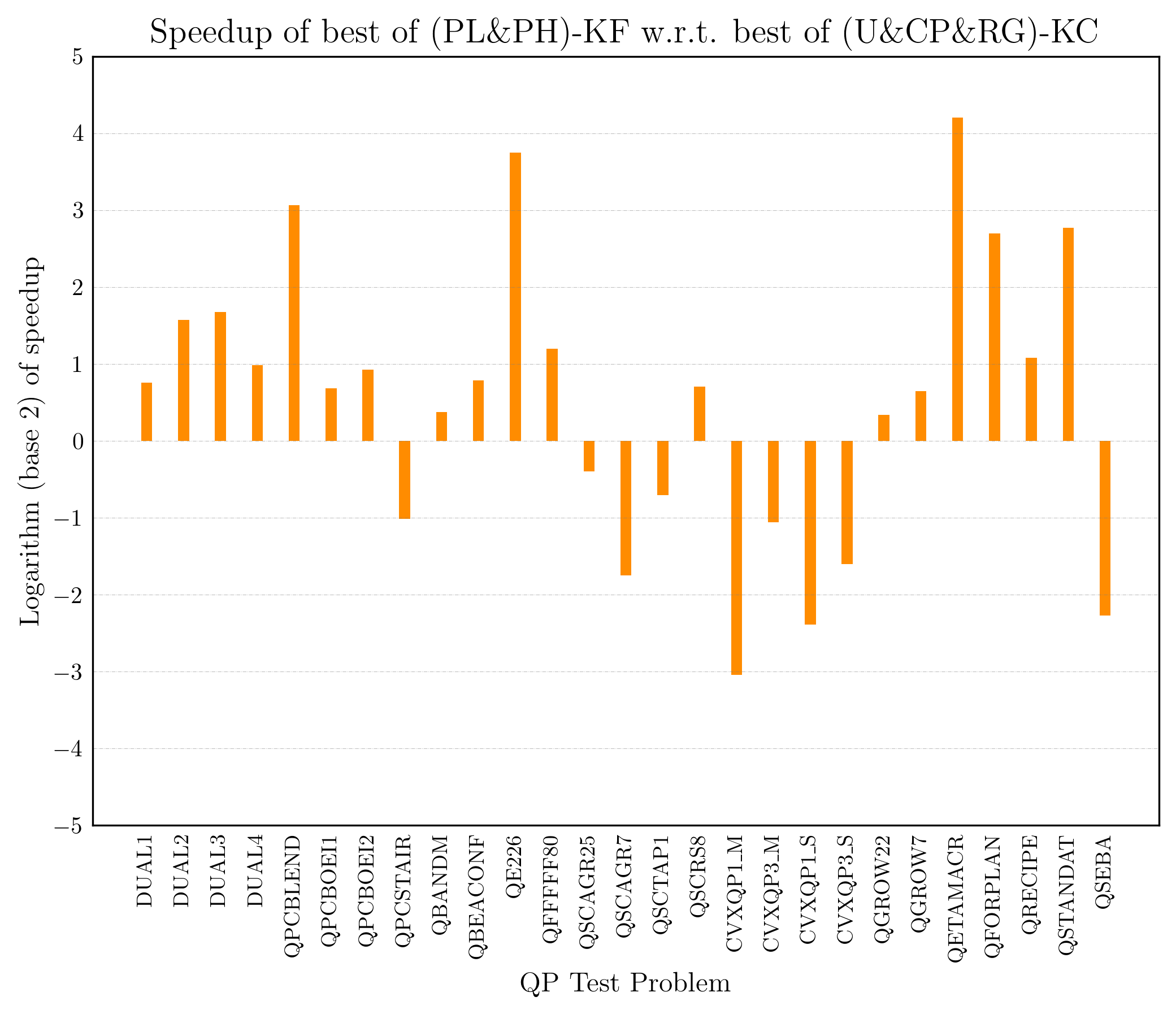}
  \caption{Speedups of PL-KF and PH-KF with respect to the fastest of the  iterative augmented system variants (U-KC, CP-KC and RG-KC) for the  $\text{MMQP}_\text{XD}$ problems for which the direct variant D-KF does not achieve the fastest overall. %\samah{median $\tilde{s}=1.628$, geomean $\hat{s}=1.432$.}
  }
  \label{speedup-barplot_vs_everything}
\end{figure}

Lastly, we directly compare our IPM variants PL-KF and PH-KF, to the alternative iterative variants namely U-KC, CP-KC and RG-KC by computing the following speedup for QP test problems, once more excluding the ones agreeable with the direct solve variant,
\begin{equation}
    s_\text{PL-PH}^\text{U-CP-RG}(q_j) = \frac{\min\bigl(c_\text{U-KC}(q_j),c_\text{CP-KC}(q_j),c_\text{RG-KC}(q_j)\bigr)}{\min\bigl(c_\text{PL-KF}(q_j),c_\text{PH-KF}(q_j)\bigr)},\;\;\;\;\;q_j\in {\{\text{MMQP}_\text{XD}\}}.
\end{equation}
%\sout{Figure~\ref{speedup-barplot_vs_everything} shows the maximum speedup among the (PL-KF, PH-KF) variants with respect to the maximum speedup among the iterative IPM variants (U-KC, CP-KC, RG-KC) for the QP test problems for which at least one of the iterative variants is faster than the direct variant (D-KC).} 
\Cref{speedup-barplot_vs_everything} presents the logarithm (base 2) of these speedup (and slowdown) results.  
Our preconditioned inexact IPM solvers (PL\&PH)-KF achieve a reduction in cost of $1.432$ (by geometric mean) relative to the best alternative inexact IPM variant. %\samahnew{, whether the linear solver is preconditioned or not}.
%median $\tilde{s}=1.628$, geomean $\hat{s}=1.432$

\section{Conclusion}\label{sec:conc}

Our theoretical analysis and experimental results demonstrate that the proposed new preconditioned iterative framework for IPM requires fewer iterations and minimal cost among alternatives in practical scenarios.
%The best choice of method is largely dependent on the sparsity structure of the Hessian and constraint matrices in the KKT system.
For QP problems arising in particular applications, the structure of the Hessian and constraint matrices is often known, so a choice among methods should be guided by how much fill is expected in factorizations employed by each approach.
%(constraint preconditioning also requires a factorization that is often very costly).
Due to factorizing $K_F$, our inexact IPM solvers require that the KKT subsystem $F$ be nonsingular and that the fill in this factorization is manageable. {In this paper, this invertibility is a enforced by a simple regularization, rendering the Hessian positive-definite. But an extension of our method to handle singular $H$ and/or $F$ is possible.}
%\sout{An extension of our method to work for semi-definite Hessian matrices would be of interest.While, in this paper, we enforce positive-definiteness by a simple regularization, extension of our approach to handle singular $H$ and/or $F$ are possible.}
{In particular, if $H$ is singular, but $F$ is full rank, our approach is applicable except for the high-d.o.f.\ preconditioner, which would need to become more sophisticated.
When $F$ is also singular, our method may be applied by taking the Schur complement of a smaller full-rank subsystem.}
%a more accurate solution of QPs with semi-definite systems may be possible with the use of a rank-revealing factorization of the Hessian or by using the regularized system as a preconditioner.

Our inexact IPM methods are a hybrid between direct and iterative solvers.
The benefit of our approach with respect to a direct method is in avoiding the need to compute a matrix factorization at each IPM step
and that the factorized matrix %\sout{does not include inequality constraints.} 
is independent of the inequality constraints. %and the \samah{re-}use of \sout{the factorization of a matrix smaller than the original KKT system.}\sout{a factorization of a simpler matrix than the whole KKT system.}
%Edgar: actually looking back, F has the same sparsity and dimension as K_C, unless C^TC has more nonzeros than H, which is rare, so we should not say the system is harder/bigger. It would be denser in some cases, so tried to make sentence more specific.
Relative to previously existing preconditioners, the preconditioned reduced system we employ is significantly better conditioned, which enables us to reduce iteration count and cost.
The positive results obtained in our cost evaluations and conditioning studies suggest that a high-performance implementation of our method might achieve favorable performance to existing high-performance linear solvers for IPM.

\pagebreak

\bibliographystyle{siamplain}
\bibliography{refs-ipmp.bib}

\begin{thebibliography}{10}

\bibitem{al2008preconditioning}
{\sc G.~Al-Jeiroudi, J.~Gondzio, and J.~Hall}, {\em Preconditioning indefinite
  systems in interior point methods for large scale linear optimisation},
  Optimisation Methods and Software, 23 (2008), pp.~345--363.

\bibitem{ali2018iterative}
{\sc F.~P. Ali~Beik and M.~Benzi}, {\em Iterative methods for double saddle
  point systems}, SIAM Journal on Matrix Analysis and Applications, 39 (2018),
  pp.~902--921.

\bibitem{altman1999regularized}
{\sc A.~Altman and J.~Gondzio}, {\em Regularized symmetric indefinite systems
  in interior point methods for linear and quadratic optimization},
  Optimization Methods and Software, 11 (1999), pp.~275--302.

\bibitem{arioli1989augmented}
{\sc M.~Arioli, I.~S. Duff, and P.~P. de~Rijk}, {\em On the augmented system
  approach to sparse least-squares problems}, Numerische Mathematik, 55 (1989),
  pp.~667--684.

\bibitem{axelsson2003preconditioning}
{\sc O.~Axelsson and M.~Neytcheva}, {\em Preconditioning methods for linear
  systems arising in constrained optimization problems}, Numerical linear
  algebra with applications, 10 (2003), pp.~3--31.

\bibitem{bank1989class}
{\sc R.~E. Bank, B.~D. Welfert, and H.~Yserentant}, {\em A class of iterative
  methods for solving saddle point problems}, Numerische Mathematik, 56 (1989),
  pp.~645--666.

\bibitem{bellavia1998inexact}
{\sc S.~Bellavia}, {\em Inexact interior-point method}, Journal of Optimization
  Theory and Applications, 96 (1998), pp.~109--121.

\bibitem{bellavia2015updating}
{\sc S.~Bellavia, V.~De~Simone, D.~di~Serafino, and B.~Morini}, {\em Updating
  constraint preconditioners for {KKT} systems in quadratic programming via
  low-rank corrections}, SIAM Journal on Optimization, 25 (2015),
  pp.~1787--1808.

\bibitem{benzi2019uzawa}
{\sc M.~Benzi and F.~P.~A. Beik}, {\em Uzawa-type and augmented lagrangian
  methods for double saddle point systems}, in Structured Matrices in Numerical
  Linear Algebra, Springer, 2019, pp.~215--236.

\bibitem{benzi2005numerical}
{\sc M.~Benzi, G.~H. Golub, and J.~Liesen}, {\em Numerical solution of saddle
  point problems}, Acta numerica, 14 (2005), pp.~1--137.

\bibitem{benzi2006augmented}
{\sc M.~Benzi and M.~A. Olshanskii}, {\em An augmented lagrangian-based
  approach to the oseen problem}, SIAM Journal on Scientific Computing, 28
  (2006), pp.~2095--2113.

\bibitem{benzi2008some}
{\sc M.~Benzi and A.~J. Wathen}, {\em Some preconditioning techniques for
  saddle point problems}, in Model Order Reduction: Theory, Research Aspects
  and Applications, Springer, 2008, pp.~195--211.

\bibitem{bergamaschi2021new}
{\sc L.~Bergamaschi, J.~Gondzio, {\'A}.~Mart{\'\i}nez, J.~W. Pearson, and
  S.~Pougkakiotis}, {\em A new preconditioning approach for an interior
  point-proximal method of multipliers for linear and convex quadratic
  programming}, Numerical Linear Algebra with Applications, 28 (2021),
  p.~e2361.

\bibitem{bergamaschi2004preconditioning}
{\sc L.~Bergamaschi, J.~Gondzio, and G.~Zilli}, {\em Preconditioning indefinite
  systems in interior point methods for optimization}, Computational
  Optimization and Applications, 28 (2004), pp.~149--171.

\bibitem{bocanegra2013improving}
{\sc S.~Bocanegra, J.~Castro, and A.~R. Oliveira}, {\em Improving an
  interior-point approach for large block-angular problems by hybrid
  preconditioners}, European Journal of Operational Research, 231 (2013),
  pp.~263--273.

\bibitem{cao2008augmentation}
{\sc Z.-H. Cao}, {\em Augmentation block preconditioners for saddle point-type
  matrices with singular (1, 1) blocks}, Numerical Linear Algebra with
  Applications, 15 (2008), pp.~515--533.

\bibitem{CASACIO2017129}
{\sc L.~Casacio, C.~Lyra, A.~R.~L. Oliveira, and C.~O. Castro}, {\em Improving
  the preconditioning of linear systems from interior point methods}, Computers
  and Operations Research, 85 (2017), pp.~129--138.

\bibitem{chow1997approximate}
{\sc E.~Chow and Y.~Saad}, {\em Approximate inverse techniques for
  block-partitioned matrices}, SIAM Journal on Scientific Computing, 18 (1997),
  pp.~1657--1675.

\bibitem{dollar2007using}
{\sc H.~S. Dollar, N.~I. Gould, W.~H. Schilders, and A.~J. Wathen}, {\em Using
  constraint preconditioners with regularized saddle-point problems},
  Computational optimization and applications, 36 (2007), pp.~249--270.

\bibitem{dravzic2015sparsity}
{\sc M.~D. Dra{\v{z}}i{\'c}, R.~P. Lazovi{\'c}, and V.~V.
  Kova{\v{c}}evi{\'c}-Vuj{\v{c}}i{\'c}}, {\em Sparsity preserving
  preconditioners for linear systems in interior-point methods}, Computational
  Optimization and Applications, 61 (2015), pp.~557--570.

\bibitem{duff2017direct}
{\sc I.~S. Duff, A.~M. Erisman, and J.~K. Reid}, {\em Direct methods for sparse
  matrices}, Oxford University Press, 2017.

\bibitem{durazzi2003indefinitely}
{\sc C.~Durazzi and V.~Ruggiero}, {\em Indefinitely preconditioned conjugate
  gradient method for large sparse equality and inequality constrained
  quadratic problems}, Numerical linear algebra with applications, 10 (2003),
  pp.~673--688.

\bibitem{ewing1990preconditioning}
{\sc R.~E. Ewing, R.~D. Lazarov, P.~Lu, and P.~S. Vassilevski}, {\em
  Preconditioning indefinite systems arising from mixed finite element
  discretization of second-order elliptic problems}, in Preconditioned
  conjugate gradient methods, Springer, 1990, pp.~28--43.

\bibitem{farrell2019augmented}
{\sc P.~E. Farrell, L.~Mitchell, and F.~Wechsung}, {\em An augmented lagrangian
  preconditioner for the 3d stationary incompressible navier--stokes equations
  at high reynolds number}, SIAM Journal on Scientific Computing, 41 (2019),
  pp.~A3073--A3096.

\bibitem{fletcher2013practical}
{\sc R.~Fletcher}, {\em Practical methods of optimization}, John Wiley \& Sons,
  2013.

\bibitem{forsgren2002inertia}
{\sc A.~Forsgren}, {\em Inertia-controlling factorizations for optimization
  algorithms}, Applied Numerical Mathematics, 43 (2002), pp.~91--107.

\bibitem{forsgren2007iterative}
{\sc A.~Forsgren, P.~E. Gill, and J.~D. Griffin}, {\em Iterative solution of
  augmented systems arising in interior methods}, SIAM Journal on Optimization,
  18 (2007), pp.~666--690.

\bibitem{forsgren2002interior}
{\sc A.~Forsgren, P.~E. Gill, and M.~H. Wright}, {\em Interior methods for
  nonlinear optimization}, SIAM review, 44 (2002), pp.~525--597.

\bibitem{fortin2000augmented}
{\sc M.~Fortin and R.~Glowinski}, {\em Augmented Lagrangian methods:
  applications to the numerical solution of boundary-value problems}, Elsevier,
  2000.

\bibitem{freund1999convergence}
{\sc R.~W. Freund, F.~Jarre, and S.~Mizuno}, {\em Convergence of a class of
  inexact interior-point algorithms for linear programs}, Mathematics of
  Operations Research, 24 (1999), pp.~50--71.

\bibitem{freund1995software}
{\sc R.~W. Freund and N.~M. Nachtigal}, {\em Software for simplified {Lanczos}
  and {QMR} algorithms}, Applied Numerical Mathematics, 19 (1995),
  pp.~319--341.

\bibitem{gansterer2003mathematical}
{\sc W.~N. Gansterer, J.~Schneid, and C.~W. Ueberhuber}, {\em Mathematical
  properties of equilibrium systems}, tech. report, University of Vienna, 2003.

\bibitem{gatica2000dual}
{\sc G.~N. Gatica and N.~Heuer}, {\em A dual-dual formulation for the coupling
  of mixed-fem and bem in hyperelasticity}, SIAM Journal on Numerical Analysis,
  38 (2000), pp.~380--400.

\bibitem{gertz2003object}
{\sc E.~M. Gertz and S.~J. Wright}, {\em Object-oriented software for quadratic
  programming}, ACM Transactions on Mathematical Software (TOMS), 29 (2003),
  pp.~58--81.

\bibitem{gill1991solving}
{\sc P.~E. Gill, W.~Murray, D.~B. Ponceleon, and M.~A. Saunders}, {\em Solving
  reduced kkt systems in barrier methods for linear and quadratic programming},
  tech. report, STANFORD UNIV CA SYSTEMS OPTIMIZATION LAB, 1991.

\bibitem{gill1992preconditioners}
{\sc P.~E. Gill, W.~Murray, D.~B. Poncele{\'o}n, and M.~A. Saunders}, {\em
  Preconditioners for indefinite systems arising in optimization}, SIAM journal
  on matrix analysis and applications, 13 (1992), pp.~292--311.

\bibitem{gill1993solving}
{\sc P.~E. Gill, W.~Murray, D.~B. Poncele{\'o}n, M.~A. Saunders, G.~Watson, and
  D.~Griffiths}, {\em Solving reduced kkt systems in barrier methods for linear
  programming}, Numerical Analysis,  (1993), pp.~89--104.

\bibitem{glowinski1989augmented}
{\sc R.~Glowinski and P.~Le~Tallec}, {\em Augmented Lagrangian and
  operator-splitting methods in nonlinear mechanics}, SIAM, 1989.

\bibitem{golub2003solving}
{\sc G.~H. Golub and C.~Greif}, {\em On solving block-structured indefinite
  linear systems}, SIAM Journal on Scientific Computing, 24 (2003),
  pp.~2076--2092.

\bibitem{golub2013matrix}
{\sc G.~H. Golub and C.~F. Van~Loan}, {\em Matrix computations}, vol.~3, JHU
  press, 2013.

\bibitem{gondzio2012matrix}
{\sc J.~Gondzio}, {\em Matrix-free interior point method}, Computational
  Optimization and Applications, 51 (2012), pp.~457--480.

\bibitem{gould2001solution}
{\sc N.~I. Gould, M.~E. Hribar, and J.~Nocedal}, {\em On the solution of
  equality constrained quadratic programming problems arising in optimization},
  SIAM Journal on Scientific Computing, 23 (2001), pp.~1376--1395.

\bibitem{gould2015cutest}
{\sc N.~I. Gould, D.~Orban, and P.~L. Toint}, {\em Cutest: a constrained and
  unconstrained testing environment with safe threads for mathematical
  optimization}, Computational optimization and applications, 60 (2015),
  pp.~545--557.

\bibitem{greif2004augmented}
{\sc C.~Greif, G.~H. Golub, and J.~M. Varah}, {\em Augmented lagrangian
  techniques for solving saddle point linear systems}, matrix, 500 (2004),
  pp.~1--1.

\bibitem{doi:10.1137/120890600}
{\sc C.~Greif, E.~Moulding, and D.~Orban}, {\em Bounds on eigenvalues of
  matrices arising from interior-point methods}, SIAM Journal on Optimization,
  24 (2014), pp.~49--83.

\bibitem{greif2006preconditioners}
{\sc C.~Greif and D.~Sch{\"o}tzau}, {\em Preconditioners for saddle point
  linear systems with highly singular (1, 1) blocks}, ETNA, Special Volume on
  Saddle Point Problems, 22 (2006), pp.~114--121.

\bibitem{greif2007preconditioners}
{\sc C.~Greif and D.~Sch{\"o}tzau}, {\em Preconditioners for the discretized
  time-harmonic {Maxwell} equations in mixed form}, Numerical Linear Algebra
  with Applications, 14 (2007), pp.~281--297.

\bibitem{haws2002preconditioning}
{\sc J.~C. Haws}, {\em Preconditioning {KKT} systems}, PhD thesis, North
  Carolina State University, 2002.

\bibitem{haynsworth1968inertia}
{\sc E.~V. Haynsworth and A.~M. Ostrowski}, {\em On the inertia of some classes
  of partitioned matrices}, Linear Algebra and its Applications, 1 (1968),
  pp.~299--316.

\bibitem{hestenes1952methods}
{\sc M.~R. Hestenes, E.~Stiefel, et~al.}, {\em Methods of conjugate gradients
  for solving linear systems}, Journal of research of the National Bureau of
  Standards, 49 (1952), pp.~409--436.

\bibitem{keller2000constraint}
{\sc C.~Keller, N.~I. Gould, and A.~J. Wathen}, {\em Constraint preconditioning
  for indefinite linear systems}, SIAM Journal on Matrix Analysis and
  Applications, 21 (2000), pp.~1300--1317.

\bibitem{krzyzanowski2001block}
{\sc P.~Krzyzanowski}, {\em On block preconditioners for nonsymmetric saddle
  point problems}, SIAM Journal on Scientific Computing, 23 (2001),
  pp.~157--169.

\bibitem{lukvsan1998indefinitely}
{\sc L.~Luk{\v{s}}an and J.~Vl{\v{c}}ek}, {\em Indefinitely preconditioned
  inexact newton method for large sparse equality constrained non-linear
  programming problems}, Numerical linear algebra with applications, 5 (1998),
  pp.~219--247.

\bibitem{lustig1992implementing}
{\sc I.~J. Lustig, R.~E. Marsten, and D.~F. Shanno}, {\em On implementing
  {Mehrotra’s} predictor--corrector interior-point method for linear
  programming}, SIAM Journal on Optimization, 2 (1992), pp.~435--449.

\bibitem{mardal2004uniform}
{\sc K.-A. Mardal and R.~Winther}, {\em Uniform preconditioners for the time
  dependent stokes problem}, Numerische Mathematik, 98 (2004), pp.~305--327.

\bibitem{maros1999repository}
{\sc I.~Maros and C.~M{\'e}sz{\'a}ros}, {\em A repository of convex quadratic
  programming problems}, Optimization Methods and Software, 11 (1999),
  pp.~671--681.

\bibitem{mihajlovic2004efficient}
{\sc M.~D. Mihajlovi{\'c} and D.~J. Silvester}, {\em Efficient parallel solvers
  for the biharmonic equation}, Parallel Computing, 30 (2004), pp.~35--55.

\bibitem{morini2016spectral}
{\sc B.~Morini, V.~Simoncini, and M.~Tani}, {\em Spectral estimates for
  unreduced symmetric kkt systems arising from interior point methods},
  Numerical Linear Algebra with Applications, 23 (2016), pp.~776--800.

\bibitem{morini2017comparison}
{\sc B.~Morini, V.~Simoncini, and M.~Tani}, {\em A comparison of reduced and
  unreduced {KKT} systems arising from interior point methods}, Computational
  Optimization and Applications, 68 (2017), pp.~1--27.

\bibitem{NoceWrig06}
{\sc J.~Nocedal and S.~J. Wright}, {\em Numerical Optimization}, Springer, New
  York, NY, USA, second~ed., 2006.

\bibitem{oliveira2005new}
{\sc A.~R. Oliveira and D.~C. Sorensen}, {\em A new class of preconditioners
  for large-scale linear systems from interior point methods for linear
  programming}, Linear Algebra and its applications, 394 (2005), pp.~1--24.

\bibitem{paige1975solution}
{\sc C.~C. Paige and M.~A. Saunders}, {\em Solution of sparse indefinite
  systems of linear equations}, SIAM journal on numerical analysis, 12 (1975),
  pp.~617--629.

\bibitem{pearson2017fast}
{\sc J.~W. Pearson and J.~Gondzio}, {\em Fast interior point solution of
  quadratic programming problems arising from pde-constrained optimization},
  Numerische Mathematik, 137 (2017), pp.~959--999.

\bibitem{pearson2020preconditioners}
{\sc J.~W. Pearson and J.~Pestana}, {\em Preconditioners for krylov subspace
  methods: An overview}, GAMM-Mitteilungen, 43 (2020), p.~e202000015.

\bibitem{perugia2000block}
{\sc I.~Perugia and V.~Simoncini}, {\em Block-diagonal and indefinite symmetric
  preconditioners for mixed finite element formulations}, Numerical linear
  algebra with applications, 7 (2000), pp.~585--616.

\bibitem{pougkakiotis2019dynamic}
{\sc S.~Pougkakiotis and J.~Gondzio}, {\em Dynamic non-diagonal regularization
  in interior point methods for linear and convex quadratic programming},
  Journal of Optimization Theory and Applications, 181 (2019), pp.~905--945.

\bibitem{pougkakiotis2021interior}
{\sc S.~Pougkakiotis and J.~Gondzio}, {\em An interior point-proximal method of
  multipliers for convex quadratic programming}, Computational Optimization and
  Applications, 78 (2021), pp.~307--351.

\bibitem{powell2003optimal}
{\sc C.~E. Powell and D.~Silvester}, {\em Optimal preconditioning for
  raviart--thomas mixed formulation of second-order elliptic problems}, SIAM
  journal on matrix analysis and applications, 25 (2003), pp.~718--738.

\bibitem{ramage2013preconditioned}
{\sc A.~Ramage and E.~C. Gartland~Jr}, {\em A preconditioned nullspace method
  for liquid crystal director modeling}, SIAM Journal on Scientific Computing,
  35 (2013), pp.~B226--B247.

\bibitem{rees2007preconditioner}
{\sc T.~Rees and C.~Greif}, {\em A preconditioner for linear systems arising
  from interior point optimization methods}, SIAM Journal on Scientific
  Computing, 29 (2007), pp.~1992--2007.

\bibitem{rozloznik2002krylov}
{\sc M.~Rozlozn{\'\i}k and V.~Simoncini}, {\em Krylov subspace methods for
  saddle point problems with indefinite preconditioning}, SIAM Journal on
  Matrix Analysis and Applications, 24 (2002), pp.~368--391.

\bibitem{saad2003iterative}
{\sc Y.~Saad}, {\em Iterative methods for sparse linear systems}, SIAM, 2003.

\bibitem{saad1986gmres}
{\sc Y.~Saad and M.~H. Schultz}, {\em {GMRES}: A generalized minimal residual
  algorithm for solving nonsymmetric linear systems}, SIAM Journal on
  scientific and statistical computing, 7 (1986), pp.~856--869.

\bibitem{saunders1996cholesky}
{\sc M.~A. Saunders et~al.}, {\em Cholesky-based methods for sparse least
  squares: The benefits of regularization}, Linear and nonlinear conjugate
  gradient-related methods, 100 (1996), pp.~92--100.

\bibitem{shen2012augmentation}
{\sc S.-Q. Shen, T.-Z. Huang, and J.-S. Zhang}, {\em Augmentation block
  triangular preconditioners for regularized saddle point problems}, SIAM
  Journal on Matrix Analysis and Applications, 33 (2012), pp.~721--741.

\bibitem{silvester1994fast}
{\sc D.~Silvester and A.~Wathen}, {\em Fast iterative solution of stabilised
  stokes systems part ii: Using general block preconditioners}, SIAM Journal on
  Numerical Analysis, 31 (1994), pp.~1352--1367.

\bibitem{sogn2019schur}
{\sc J.~Sogn and W.~Zulehner}, {\em Schur complement preconditioners for
  multiple saddle point problems of block tridiagonal form with application to
  optimization problems}, IMA Journal of Numerical Analysis, 39 (2019),
  pp.~1328--1359.

\bibitem{toh2004block}
{\sc K.-C. Toh, K.-K. Phoon, and S.-H. Chan}, {\em Block preconditioners for
  symmetric indefinite linear systems}, International Journal for Numerical
  Methods in Engineering, 60 (2004), pp.~1361--1381.

\bibitem{tong1998iterative}
{\sc Z.~Tong and A.~Sameh}, {\em On an iterative method for saddle point
  problems}, Numerische Mathematik, 79 (1998), pp.~643--646.

\bibitem{van1992bi}
{\sc H.~A. Van~der Vorst}, {\em {Bi-CGSTAB}: {A} fast and smoothly converging
  variant of {Bi-CG} for the solution of nonsymmetric linear systems}, SIAM
  Journal on scientific and Statistical Computing, 13 (1992), pp.~631--644.

\bibitem{vandenberghe2010cvxopt}
{\sc L.~Vandenberghe}, {\em The cvxopt linear and quadratic cone program
  solvers}, Online: http://cvxopt. org/documentation/coneprog. pdf,  (2010).

\bibitem{wachter2006implementation}
{\sc A.~W{\"a}chter and L.~T. Biegler}, {\em On the implementation of an
  interior-point filter line-search algorithm for large-scale nonlinear
  programming}, Mathematical programming, 106 (2006), pp.~25--57.

\bibitem{wang2000adaptive}
{\sc W.~Wang and D.~P. O'Leary}, {\em Adaptive use of iterative methods in
  predictor--corrector interior point methods for linear programming},
  Numerical Algorithms, 25 (2000), pp.~387--406.

\bibitem{wathen1995convergence}
{\sc A.~Wathen, B.~Fischer, and D.~Silvester}, {\em The convergence rate of the
  minimal residual method for the stokes problem}, Numerische Mathematik, 71
  (1995), pp.~121--134.

\bibitem{wright1992interior}
{\sc M.~H. Wright}, {\em Interior methods for constrained optimization}, Acta
  numerica, 1 (1992), pp.~341--407.

\bibitem{wright1997primal}
{\sc S.~J. Wright}, {\em Primal-dual interior-point methods}, SIAM, 1997.

\bibitem{zhang1998solving}
{\sc Y.~Zhang}, {\em Solving large-scale linear programs by interior-point
  methods under the {MATLAB} environment}, Optimization Methods and Software,
  10 (1998), pp.~1--31.

\end{thebibliography}

\begin{table}[htbp]
    \centering
    \small{
    \begin{tabular}{|c|c|c|c|}
    \hline
    QP Problem & $\min({\kappa}(\text U),{\kappa}(\text{CP}),{\kappa}(\text{RG}))$ &  $\min({\kappa}(\text{PL}),{\kappa}(\text{PH}))$& ratio\\ \hline
        DUAL1 & 7.89e+02 & 2.16e+01  & 3.66e+01 \\
\hline
DUAL2 & 2.14e+02 & 3.76e+00  & 5.70e+01 \\
\hline
DUAL3 & 4.64e+02 & 2.89e+00  & 1.61e+02 \\
\hline
DUAL4 & 6.48e+02 & 1.82e+00  & 3.56e+02 \\
\hline
DUALC1 & 1.36e+09 & 7.89e+03  & 1.73e+05 \\
\hline
DUALC5 & 1.46e+08 & 8.60e+02  & 1.70e+05 \\
\hline
QPCBLEND & 4.36e+11 & 5.60e+04  & 7.79e+06 \\
\hline
QPCBOEI1 & 5.16e+07 & 6.43e+03  & 8.03e+03 \\
\hline
QPCBOEI2 & 1.71e+08 & 1.05e+05  & 1.63e+03 \\
\hline
QPCSTAIR & 1.58e+07 & 4.81e+03  & 3.29e+03 \\
\hline
QBANDM & 3.64e+01 & 5.45e+06  & 6.68e-06 \\
\hline
QADLITTL & 3.42e+07 & 3.76e+06  & 9.11e+00 \\
\hline
QAFIRO & 3.35e+08 & 5.86e+03  & 5.71e+04 \\
\hline
QBEACONF & 3.70e+00 & 2.18e+09  & 1.70e-09 \\
\hline
QE226 & 1.78e+09 & 9.79e+04  & 1.82e+04 \\
\hline
QFFFFF80 & 1.30e+15 & 2.46e+08  & 5.28e+06 \\
\hline
QSC205 & 2.30e+10 & 6.71e+04  & 3.43e+05 \\
\hline
QSCAGR25 & 2.74e+07 & 2.32e+05  & 1.18e+02 \\
\hline
QSCAGR7 & 1.85e+07 & 2.60e+05  & 7.13e+01 \\
\hline
QSCFXM1 & 1.46e+11 & 1.41e+08  & 1.04e+03 \\
\hline
QSCFXM2 & 1.73e+12 & 3.33e+08  & 5.20e+03 \\
\hline
QSCTAP1 & 1.66e+08 & 5.39e+05  & 3.08e+02 \\
\hline
QSHARE1B & 4.79e+08 & 2.56e+07  & 1.87e+01 \\
\hline
QSHARE2B & 1.80e+09 & 2.56e+07  & 7.03e+01 \\
\hline
QSCRS8 & 6.49e+12 & 2.77e+08  & 2.34e+04 \\
\hline
CVXQP1\_M & 9.26e+05 & 6.89e+02  & 1.34e+03 \\
\hline
CVXQP3\_M & 9.25e+06 & 4.99e+02  & 1.85e+04 \\
\hline
CVXQP1\_S & 7.79e+05 & 8.99e+03  & 8.66e+01 \\
\hline
CVXQP3\_S & 3.44e+06 & 4.02e+04  & 8.55e+01 \\
\hline
QGROW15 & 4.68e+03 & 2.73e+03  & 1.71e+00 \\
\hline
QGROW22 & 9.04e+03 & 2.09e+03  & 4.33e+00 \\
\hline
QGROW7 & 9.50e+02 & 2.65e+03  & 3.59e-01 \\
\hline
VALUES & 2.39e+02 & 3.31e+03  & 7.22e-02 \\
\hline
DUALC2 & 4.25e+09 & 4.55e+04  & 9.35e+04 \\
\hline
DUALC8 & 1.65e+09 & 4.42e+03  & 3.73e+05 \\
\hline
QETAMACR & 4.28e+10 & 3.17e+07  & 1.35e+03 \\
\hline
QFORPLAN & 1.62e+13 & 1.98e+09  & 8.16e+03 \\
\hline
QRECIPE & 1.08e+08 & 1.35e+09  & 8.04e-02 \\
\hline
QSTAIR & 1.47e+07 & 7.76e+04  & 1.90e+02 \\
\hline
QSTANDAT & 1.14e+12 & 1.14e+10  & 9.96e+01 \\
\hline
QSEBA & 1.83e+06 & 8.99e+04  & 2.03e+01 \\
\hline
\hline
%Geometric mean & 6.38 & 4.23  & 1.51\\ \hline
\begin{tabular}[c]{@{}l@{}}Geometric mean\\over QP problems\end{tabular} & 3.13e+07 
 &  1.30e+05 
 &  2.40e+02 \\
\hline
    \end{tabular}}
    \caption{Conditioning of the preconditioned augmented system and the preconditioned inequality constraint reduced system (with best choices of preconditioners for each), as well as the ratio of improvement obtained by the preconditioners, defined as 
    %$ \min(\tilde{\kappa})$(U-CP-RG)$/\min(\hat{\ka%ppa})$(PL-PH).
    $ \min({\kappa}(\text{U}),{\kappa}(\text{CP}),{\kappa}(\text{RG}))/\min({\kappa}(\text{PL}),{\kappa}(\text{PH}))$.
    }
    \label{tab:cutest-cond-table}
\end{table}
\end{document}